\documentclass[final]{siamltex}
\usepackage{amsmath}
\usepackage{graphicx}
\usepackage{latexsym}
\usepackage{amssymb}
\usepackage{amsbsy}

\DeclareSymbolFont{msbm}{U}{msb}{m}{n}
\DeclareMathSymbol{\C}{\mathalpha}{msbm}{'103}
\DeclareMathSymbol{\R}{\mathalpha}{msbm}{'122}
\DeclareMathSymbol{\Z}{\mathalpha}{msbm}{'132}
\DeclareMathSymbol{\N}{\mathalpha}{msbm}{'116}

\newtheorem{remark}{Remark}
\newtheorem{lgrthm}{Algorithm}

\def\P{{\cal P}}
\def\R{{\cal R}}
\def\T{{\cal T}}
\def\RR{\mathbb R}
\def\e{\varepsilon}

\def\be{\begin{equation}}
\def\ee{\end{equation}}
\def\bea{\begin{eqnarray}}
\def\ba{\begin{array}{l}\displaystyle}
\def\eea{\end{eqnarray}}
\def\ea{\end{array}}

\def\ca{\\[+0.3cm]\displaystyle}

\def\IR{\mathop{\mbox{\rm Iround}}\nolimits}

\begin{document}
\title{Fluid Solver Independent Hybrid Methods for Multiscale Kinetic equations\thanks{This work was partially supported by
the INDAM project ``Kinetic Innovative Models for the Study of the
Behavior of Fluids in Micro/nano Electromechanical Systems''}}

\author{Giacomo Dimarco\thanks{
University of Ferrara, Department of Mathematics and CMCS,
Ferrara, Italy ({\tt Giacomo.Dimarco@unife.it}). }  \and Lorenzo
Pareschi\thanks{University of Ferrara, Department of Mathematics
and CMCS, Ferrara, Italy ({\tt Lorenzo.Pareschi@unife.it}).} }
\date{~}
\maketitle


\begin{abstract}
In some recent works \cite{dimarco1,dimarco3} we developed a
general framework for the construction of hybrid algorithms which
are able to face efficiently the multiscale nature of some
hyperbolic and kinetic problems. Here, at variance with respect to
the previous methods, we construct a method form-fitting to any
type of finite volume or finite difference scheme for the reduced
equilibrium system. Thanks to the coupling of Monte Carlo
techniques for the solution of the kinetic equations with
macroscopic methods for the limiting fluid equations, we show how
it is possible to solve multiscale fluid dynamic phenomena faster
with respect to traditional deterministic/stochastic methods for
the full kinetic equations. In addition, due to the hybrid nature
of the schemes, the numerical solution is affected by less
fluctuations when compared to standard Monte Carlo schemes.
Applications to the Boltzmann-BGK equation are presented to show
the performance of the new methods in comparison with classical
approaches used in the simulation of kinetic equations.
\end{abstract}

%
%
%
\maketitle




\section{Introduction}

The classical fluid dynamic models like the Navier-Stokes or the
Euler equations are not always satisfactory when dealing with
large temperature or very low densities, and a more detailed
analysis becomes often necessary to obtain correct values of the
macroscopic quantities. In such cases a kinetic approach based on
the Boltzmann equation \cite{Cer} is used. The introduction of
such a model is closely linked with the introduction of strong
difficulties from the numerical and computational point of view.
In fact, the system of equations to solve becomes very large,
especially in multidimensional situations, and even with computers
of the last generation the computational cost of a direct
discretization is often prohibitive. Moreover the Boltzmann
collision term that characterizes the kinetic equation, is very
hard to treat in practice due to its nonlinear nature and physical
properties. To these aims, probabilistic techniques such as Direct
Simulation Monte Carlo (DSMC) are extensively used in real
simulations for their great flexibility, capability of treating
different collision terms and low computational cost compared to
deterministic schemes for kinetic equations \cite{bird, Cf, sliu,
CPima}. On the other hand solutions are affected by large
fluctuations and, in non stationary situations, the difficulty to
compute averaged quantities leads to low accurate solutions or
very expensive simulations. However, even in extremely rarefied
regimes the fluid dynamic equations still furnish correct solution
in regions of the domain where the gas is not subjected to sharp
gradient. The direct consequence is that domain decomposition
methods \cite{BLPQ, Letallec, degond1} which consider the problem
at different scales, fluid or kinetic, in different part of the
computational domain, is a practical way to take advantage of the
physics without loosing accuracy. We quote also the possibility to
improve domain decomposition schemes through a moving boundary
\cite{dimarco2, tiwari_JCP}, in order to follow discontinuities
and sharp gradients inside the domain, these methods are
particularly important in the simulation of non stationary
problems. Clearly, the exact identification of the non equilibrium
zones remains an hard task to deal with and an open research
argument.

In some recent works we proposed an alternative approach to domain
decomposition methods based on the use of different numerical
methods on the whole computational domain
\cite{dimarco1,dimarco3}. We mention here that similar hybrid
approaches have been considered in \cite{degond, WE1, WE3, CPima,
RGV2}. Even if we develop our methods in the case of rarefied gas
dynamics (RGD) the formulation we proposed permits easily
generalizations to others fields in which kinetic and hyperbolic
multiscale phenomena are present. In order to use the hybrid
approach here described, it is essential to identify a local
equilibrium function, like the Maxwellian distribution for RGD,
either analytically or numerically. This local equilibrium
originates a model reduction from the microscale to the macroscale
formulation and allows to ignore the details of microscopic
interactions in terms of simplified equations which describe the
equilibrium system.

The schemes here presented represent an important improvement with
respect to the schemes developed in \cite{dimarco1, dimarco3}
where the limiting equilibrium method was limited by construction
to a kinetic scheme. In the present work we generalize the
approach to make the method independent of the fluid solver used.
We point out that this generalization is not trivial since in an
arbitrary fluid solver we miss the kinetic information on the
distribution function which is present in a standard kinetic
scheme. The main advantage is that the method in the fluid limit
degenerates into a standard fluid solver without any additional
cost of a kinetic simulation. To our knowledge this is the first
method which satisfy this property for RGD.

Although we will focus, in the construction of the schemes, on the
simplified BGK collision model, in principle the schemes can be
extended to the full Boltzmann collision operator through the use
of time relaxed Monte Carlo methods \cite{CPmc, CPima, Trazzi,
PR}. The basic idea consists in solving the kinetic model and the
macroscopic model in the entire domain, the first through Monte
Carlo techniques which are robust in the fluid limit and the
latter through a deterministic scheme and to consider as a
solution a suitable hybrid merging of the two. A remarkable
feature of the new method is the use of the hybrid moments to
correct the stochastic moments in the pure Monte Carlo scheme.
This is an important source of fluctuations reduction in the
method.

In addition we will show that it is not necessary to keep the
number of sample particles fixed in the Monte Carlo scheme,
instead it is sufficient to describe at the particles level only
the fraction of the solution which is far from the thermodynamic
equilibrium. The immediate consequence of the above observation is
a potential reduction of the number of samples used in the Monte
Carlo solution, and thereby, of computational time and
fluctuations. These improvements are directly linked with the
decrease of the local Knudsen number which is a measure of the
rarefaction of the gas. The implementation of such methodology
produces numerical schemes which, in general, are much faster than
deterministic kinetic schemes and, for flow regimes close to the
fluid limit, also to DSMC schemes. Moreover, thanks to the general
formulation of the algorithm, a domain decomposition technique can
be directly derived forcing the Knudsen number to zero (see
\cite{dimarco2b}) in some regions of the domain.







Finally, let us observe that the method here developed is based on
the classical operator splitting for the kinetic equation. This is
essential if one want to match the fluid scheme with standard DSMC
solvers for RGD, since the latter are based on such splitting.
Even if there are well-known limitations of such splitting when
dealing with Navier-Stokes asymptotics, namely the time step has
to be related to the Knudsen number in order to describe correctly
the Navier-Stokes level \cite{Xu}, here we don't aim at an
under-resolved method at all the scales but simply at developing a
method which is asymptotic preserving (AP) in the stiff limit
\cite{Jin}. On the other hand the advantages provided by our final
method in terms of fluctuations and computational cost reduction
are essentially independent of the small scales resolution but
depend only on the Knudsen number. Note that in principle one can
improve the method coupling the DSMC solver with a Navier-Stokes
fluid description instead of a Euler one. This coupling however is
not straightforward and we don't explore this direction here.

The rest of the article is organized as follows. In Section 2 we
introduce the BGK equations and his properties. In Section 3 we
recall the general structure of the hybrid methods derived in
\cite{dimarco1, dimarco3}. Next in Section 4 the fluid solver
independent hybrid scheme is described together with an
acceleration technique and two possible ways in which the
equilibrium fraction can be increased. Section 5 is devoted to
numerical results to compare performances respect traditional
Monte Carlo and kinetic schemes. Some final considerations and
future developments are discussed in the last Section.

\section{The Boltzmann-BGK Model}

We consider the following Boltzmann-BGK kinetic model \be
\partial_t f(x,v,t) + v\cdot\nabla_{x}f(x,v,t) = \frac1{\tau(x,t)} (M_{f}(x,v,t)-f(x,v,t)),
\label{eq:B} \ee with the initial condition \be
 f(x,v,t=0)=f_{0}(x,v).
\label{eq:B1} \ee In (\ref{eq:B}) the function $f(x,v,t)$ is non
negative and describes the time evolution of the distribution of
particles with velocity $v \in \RR^{d}$ and position $x \in \Omega
\subset \RR^{D}$, with $d$ and $D$ representing the dimension in
velocity and physical space respectively. In this simplified model
the Boltzmann collision term is substituted by a relaxation
towards equilibrium. In the sequel we will work with
nondimensional quantities, in that case $\tau$, the relaxation
frequency, can be written as \be
\tau(x,t)^{-1}=\frac{C}{\varepsilon(x,t)},
 \label{tau} \ee
where $\varepsilon(x,t)$ is the Knudsen number. Here we assume
$C=1$ \cite{Coron, Pieraccini}. Others choices of the relaxation
time do not change the hybrid algorithm we will describe in next
section. Observe anyway that the ratio of deterministic and
stochastic component will be a function of the relaxation time,
being linked to the ratio of the distribution function with
respect to the Maxwellian equilibrium function as explained in
details in the next Section. In the following, for simplicity, we
will skip the space and time dependency of the Knudsen number thus
$\varepsilon(x,t)=\varepsilon$.

The local Maxwellian function, representing the local equilibrium,
is defined by \be
 M_{f}(\varrho,u,T)(x,v,t)=\frac{\varrho}{(2\pi T)^{d/2}}\exp\left(\frac{-|u-v|^{2}}{2T}\right) ,
\label{eq:M} \ee where $\varrho$, $u$, $T$ are the density, mean
velocity and temperature of the gas in the x-position and at time
$t$\be
 \varrho=\int_{\RR^d} fdv, \quad u=\frac{1}{\varrho}\int_{\RR^d} vfdv, \quad T=\frac{1}{d\varrho} \int_{\RR^d}|v-u|^{2}fdv,
\label{eq:Mo} \ee while the energy $E$ is defined as \be
E=\frac{1}{2} \int_{\RR^d}|v|^{2}fdv.
 \label{eq:E} \ee
Consider now the BGK equation (\ref{eq:B}) and multiply it for
$1$, $v$, $\frac{1}{2}|v|^{2}$, the so-called collision invariant.
By integrating in $v$ the above quantities, the equations for the
first three moments of the distribution function $f$ are obtained.
They describe respectively the conservations laws for mass,
momentum and energy. Unfortunately, the system obtained through
the above average in velocity space is not closed since it
involves higher order moments of the distribution function.

Note that, formally from (\ref{eq:B}) as $\varepsilon \rightarrow
0$, the function $f$ approaches the local Maxwellian. In this case
it is possible to compute analytically the higher moments of $f$
from $\varrho$, $u$ and $T$. Carrying on this computation we
obtain the set of compressible Euler equations (see
\cite{cercignani} for details)

\be \ba \frac{\partial \varrho}{\partial t} + \nabla_{x}
\cdot(\varrho u) = 0\ca \frac{\partial \varrho u}{\partial t} +
\nabla_{x} \cdot (\varrho u \otimes u+p) =  0, \ca \frac{\partial
E}{\partial t} +\nabla_x\cdot(Eu+pu) =  0 \ca p=\varrho T, \ \
E=\frac{d}{2}\varrho T +\frac{1}{2} \varrho |u|^{2} \ea
\label{eq:sys1} \ee where $p$ is the thermodynamical pressure
while $\otimes$ represents a tensor product. Higher order fluid
model, like Navier-Stokes, can be derived similarly
\cite{cercignani}.

\section{Hybrid Methods}

The schemes derived in this paper are based on the same hybrid
representation defined in \cite{dimarco3}. Here we recall only the
key points of the previous method, for details we remind to
\cite{dimarco3}.

For a fixed space point $x$ we can interpret the distribution
function as a probability density in the velocity space (the
x-dependence is omitted)
\begin{equation}
f(v,t)\geq 0, \quad \varrho=\int_{-\infty}^{+\infty} f(v,t)dv=1.
\end{equation}
Next we recall the following definition of hybrid representation
\cite{dimarco3}
\begin{definition}
Given a probability density $f(v,t)$, and a probability density
$M_f(v,t)$, called equilibrium density, we define $w(v,t) \in
[0,1]$ and $\tilde f(v,t)\geq 0$ in the following way
\begin{equation}
w(v,t)=\left\lbrace
\begin{array}{ll}
\displaystyle\frac{f(v,t)}{M_f(v,t)}, & f(v,t) \leq M_f(v,t)\neq 0\\
1, & f(v,t) \geq M_f(v,t)\\
\end{array}
\right.
\end{equation}
and
\begin{equation}
\tilde f(v,t)=f(v,t)-w(v,t) M_f(v,t).
\end{equation}
Thus $f(v,t)$ can be represented as 
\begin{equation}
f(v,t)=\tilde f(v,t)+w(v,t)M_f(v,t). \label{eq:r1}
\end{equation}
\label{def:1}
\end{definition}
If we take now $\beta(t)=\min_{v}\{w(v,t)\}$ and $\tilde
f(v,t)=f(v,t)-\beta(t) M_f(v,t)$ we have
\[
\int_{v} \tilde f(v,t) dv =1-\beta(t).
\]
Let us define for $\beta(t)\neq 1$ the probability density
\[
f^p(v,t) =\frac{\tilde f(v,t)}{1-\beta(t)}.
\]
The case $\beta(t)=1$ is trivial since it implies
$f(v,t)=M_f(v,t)$. Thus the probability density $f(v,t)$, can be
written as a convex combination of two probability densities in
the form \cite{CPmc,
CPima} 
\begin{equation}
f(v,t)=(1-\beta(t)) f^p(v,t)+\beta(t) M_f(v,t). \label{eq:r2}
\end{equation}
Clearly the above representation is a particular case of
(\ref{eq:r1}).

Now we consider the following general representation, including
space dependence \be f(x,v,t)=\tilde f(x,v,t)+w(x,v,t)M_f(x,v,t),
\label{eq:pj} \ee where $w(x,v,t)\geq 0$ is a function that
characterizes the equilibrium fraction and $\tilde f(x,v,t)$ the
non equilibrium part of the distribution function. This
representation in general can be obtained for the initial data of
the kinetic equation using directly Definition 1.

The starting point of the method is the classical operator
splitting which consists in solving first a homogeneous relaxation
step
\be
\partial_t f^r(x,v,t)  =-\frac1{\e}(f^r(x,v,t)-M^{r}(x,v,t))
\label{eq:coll}\ee and then a free transport equation \be
\partial_t f^c(x,v,t) + v\cdot \nabla_{x} f^c(x,v,t)=0.
\label{eq:trasp} \ee

In a single time step $\Delta t$ the computation of the hybrid
method derived in \cite{dimarco3} can be summarized as follows
\begin{itemize}
  \item Starting from a function $f^r(x,v,t)=f(x,v,t)$ in the form (\ref{eq:pj}) solve the
  relaxation step (\ref{eq:coll}) either analytically or with
  a suitable numerical time integrator for stiff ODEs, like
  backward
  Euler. This originates the decomposition
  \begin{eqnarray*}
  f^r(x,v,t+\Delta t)&=&\lambda f^r(x,v,t)+ (1-\lambda)
  M_f(x,v,t)\\
  &=&\lambda \tilde f(x,v,t)+(1-\lambda+\lambda
  w(x,v,t))M_f(x,v,t),
  \end{eqnarray*}
  with $0\leq \lambda=\lambda(\Delta t/\e) \leq 1$ a scheme dependent constant such that $\lambda\to 0$ as $\Delta t/\e \to \infty$. This decomposition
  can be cast again in the form (\ref{eq:pj}) taking $\tilde f^r(x,v,t+\Delta
  t)=\lambda \tilde f(x,v,t)$ and $w^r(x,v,t+\Delta t)=1-\lambda+\lambda
  w(x,v,t)$.
  \begin{enumerate}
  \item The new value $w^r(x,v,t+\Delta t)$ follows directly from
  the choice of $\lambda$, so from the time solver used for (\ref{eq:coll}).
  \item The new value $\tilde f^r(x,v,t+\Delta t)$ is computed by a Monte Carlo method simply
  discarding a fraction of the samples since $0\leq \lambda \leq 1$ and so $w^r(x,v,t+\Delta t)\geq
  w(x,v,t)$.
  \end{enumerate}
  \item Starting from the function $f^c(x,v,t)=f^r(x,v,t+\Delta t)$ in the form (\ref{eq:pj}) computed above
  solve the transport step (\ref{eq:trasp}).
  \begin{enumerate}
  \item Transport the particle fraction $\tilde f^c(x,v,t)$ by simple particles shifts.
  \item Transport the deterministic fraction $w^c(x,v,t)M_f(x,v,t)$ by a deterministic
  scheme.
  \item Project the computed hybrid solution $f(x,v,t+\Delta t)$ to the form
  (\ref{eq:pj}) using Definition 1.
  \end{enumerate}
\end{itemize}
Clearly point 3 of the transport step is crucial for the details
of the hybrid method. Note that point 2 of the transport step
involves the solution of a so-called kinetic scheme for the Euler
equations\cite{Deshpande, perthame}.

In the sequel we will describe the Fluid Solver Independent (FSI)
schemes which remove the limitations given by the use of a kinetic
scheme. One major difference with respect to the hybrid scheme
described above is that a common value for the equilibrium
fraction in velocity space has to be chosen
$\beta(x,t)=\min_{v}\{w(x,v,t)\}$.

\section{Fluid Solver Independent Hybrid Methods}
The key feature of FSI methods is to take advantage from the
solution of the equilibrium part of the distribution function
through a macroscopic scheme instead of a kinetic scheme. Besides
its generality, this new feature, could, in principle, lead to a
strong reduction of the computational time with respect to any
kinetic scheme for the fluid equation.


In order to describe the FSI method we introduce the projection
operator $\P$, and, in a time step $\Delta t$, the relaxation
operator $\R_{\Delta t}$ and the transport operator $\T_{\Delta
t}$. The projection operator computes from the kinetic variable
$f$ (or $M_{f}$) the macroscopic averages
$U(x,t)=(\varrho(x,t),\rho u(x,t),E(x,t))$, thus
\be\P(f(x,v,t))=U(x,t), \quad \P(M_{f}(x,v,t))=U(x,t),\ee since
the local Maxwellian $M_f$ has the same moments of the
distribution function $f$. The relaxation and transport operators
solve the relaxation and transport steps. The first has the form
\be \R_{\Delta t}(f(x,v,t))=\lambda f(x,v,t)+(1-\lambda)
M_f(x,v,t), \ee where $\lambda=\exp(-\Delta t/\varepsilon)$,
whereas the second reads \be \T_{\Delta t}(f(x,v,t))=f(x-v\Delta
t,v,t). \ee Similarly we have the approximated relaxation and
transport operators $\R_a$ and $\T_a$. For simplicity, since their
particular structure does not play any role in the general
derivation of the method, we assume in the sequel that $\T_a=\T$
and $\R_a=\R$. Note that by definition $\R_{\Delta t}(M_f)=M_f$
and so $\P(\R_{\Delta t}(M_f))=\P(M_f)$.

\subsection{A simple FSI method}

Let us start from an hybrid solution in the form
\begin{equation}
f(x,v,t)=(1-\beta(x,t)) f^p(x,v,t)+\beta(x,t) M_f(x,v,t),
\label{eq:id}
\end{equation}
where $f^p(x,v,t)$ is represented by samples so that \be
(1-\beta(x,t))f^p(x,v,t)= m^p\sum_{j=1}^{N(t)}
\delta(x-p_j(t))\delta(v-\nu_j(t)), \ee where $p_j(t)$ and
$\nu_j(t)$ represent the particles position and velocity, and
\[
m^p=\frac1{N(0)}\int_{\RR^d}\int_{\RR^D} f(x,v,0)dxdv
\]
is the mass of a single particle, while $M_f(x,v,t)$ is
represented analytically. Note that $N(t)$, namely the total
number of samples is a function of time since we keep $m^p$
constant during the simulation. This is a crucial feature of the
method since if we increase $\beta(x,t)$ in the representation
above we must decrease the number of samples $N(t)$. In practice
this implies that we will not be able to represent exactly the
fraction $\beta(x,t)$ but only its approximation corresponding to
integer sums of particles.

Since, as described in the previous section, the first relaxation
step has only the consequence of a change of $\beta(x,t)$ in
(\ref{eq:id}) we derive the method starting from the transport
step.

The transport step produces the solution \be \nonumber \T_{\Delta
t}(f(x,v,t))=(1-\beta(x-v\Delta t,t)) f^p(x-v\Delta
t,v,t)+\beta(x-v\Delta t,t) M_f(x-v\Delta t,v,t).\ee From a
practical viewpoint $(1-\beta(x-v\Delta t,t)) f^p(x-v\Delta
t,v,t)$ corresponds to solve a simple particle shift for the Monte
Carlo samples. At variance the term $\beta(x-v\Delta t,t)
M_f(x-v\Delta t,v,t)$ corresponds to a Maxwellian shift analogous
to that usually performed in the so called kinetic or Boltzmann
schemes for the Euler equations \cite{Deshpande, perthame}. The
hybrid solution for the moments $U^H(x,t+\Delta t)$ is then
recovered as \bea\nonumber U^H(x,t+\Delta t)&=&\P(\T_{\Delta
t}(f(x,v,t)))\\
\nonumber &=&\P(\T_{\Delta t}((1-\beta(x,t))
f^p(x,v,t)))+\P(\T_{\Delta t}(\beta(x,t) M_f(x,v,t)))\\
\nonumber &=& U^p(x,t+\Delta t)+U^K(x,t+\Delta t). \eea In
particular $U^K(x,t+\Delta t)$ corresponds exactly to the
approximation of the Euler solution provided by a
kinetic/Boltzmann scheme.\\
We can state the following result (see also \cite{perthame2})

\begin{theorem}
If we denote with $U^E(x,t+\Delta t)$ the solution of the Euler
equations (\ref{eq:sys1}) with initial data
$U^E(x,t)=\P(\beta(x,t)M_f(x,v,t))$ we have \be U^E(x,t+\Delta
t)=U^K(x,t+\Delta t)+O(\Delta t^2). \label{eq:est}\ee \label{th:1}
\end{theorem}
{\bf Proof.} In a time step $\Delta t$ we can write for the Euler
solution
\[
U^E(x,t+\Delta t)=U^E(x,t)+\Delta t \partial_t
U^E(x,t)+\frac12(\Delta t)^2
\partial_{tt} U^E(x,t)+O(\Delta t^3)
\]
and similarly
\[
U^K(x,t+\Delta t)=U^K(x,t)+\Delta t \partial_t
U^K(x,t)+\frac12(\Delta t)^2
\partial_{tt} U^K(x,t)+O(\Delta t^3).
\]
Clearly the zero order terms in the expansions are the same since
the initial data of the Euler equations is simply the projection
of the initial data for the transport equation
\[
U^K(x,t)=\P(\beta(x,t)M_f(x,v,t))=U^E(x,t).
\]
Now let us consider the first order terms. We have
\[
\partial_t
U^E =-(\nabla_{x} \cdot(\varrho u), \nabla_{x} \cdot (\varrho u
\otimes u+p), \nabla_x\cdot(Eu+pu))^T
\]
and
\[
\partial_t
U^K =\P(-v\cdot\nabla_x f).\] Note that this last equation is not
closed since the right hand side involves third order moments of
$f$. Again, however, the two terms evaluated at the initial time
$t$ coincide since the initial data for the transport step is the
Maxwellian fraction $\beta(x,t)M_f(x,v,t)$ and so we have the
usual Euler closure in the kinetic term. By similar arguments one
can verify that the second order terms evaluated at the initial
time are different because of the fourth order moments appearing
in $\partial_{tt}U^K(x,t)=\P(v\cdot\nabla_x(v\cdot\nabla_x f))$.
This proves (\ref{eq:est}).\\ $\Box$

By virtue of the above result we can replace the hybrid solution for
the moments after the transport with \be {\tilde U}^H(x,t+\Delta
t)=U^p(x,t+\Delta t)+U^E(x,t+\Delta t), \label{eq:fin}\ee without
affecting the overall first order accuracy of the splitting method.

This hybrid values for the moments are used to compute the new
Maxwellian $M^H_f(x,v,t+\Delta t)$ and advance the computation. To
this goal we note that the next relaxation step takes the form
\bea \nonumber \R_{\Delta t}(\T_{\Delta
t}(f(x,v,t)))&=&\lambda\T_{\Delta t}(f(x,v,t))
+(1-\lambda) M^H_f(x,v,t+\Delta t)\\
\nonumber &=& \lambda(\T_{\Delta t}((1-\beta(x,t))
f^p(x,v,t))+\T_{\Delta
t}(\beta(x,t) M_f(x,v,t)))\\
\nonumber &+&(1-\lambda) M^H_f(x,v,t+\Delta t)\\
\nonumber &=&(1-\beta(x,t+\Delta t)) f^p(x,v,t+\Delta
t)+\beta(x,t+\Delta t) M^H_f(x,v,t+\Delta t),\eea where we set \be
\beta(x,t+\Delta t)=1-\lambda, \qquad f^p(x,v,t+\Delta t)=\T_{\Delta
t}(f(x,v,t))\label{eq:rel}\ee with $\lambda=e^{-\Delta
t/\varepsilon}$. This shows that in order to compute the new
particle fraction we need to sample particles from $\T_{\Delta
t}(\beta(x,t)M_f(x,v,t))$. In practice this can be realized in a
simple way transforming initially the equilibrium Maxwellian part
$\beta(x,t)M_f(x,v,t)$ into samples and then advecting the whole set
of samples.

Let us denote with $\T_{\Delta t}(\beta(x,t)M^p_f(x,v,t))$ this set
of advected equilibrium samples. Computationally this means that at
each time step me must solve the full BGK model with a Monte Carlo
scheme \cite{dimarco3} together with a suitable deterministic solver
for the Euler equation. We can improve the efficiency of the above
algorithm observing that we do not need to transform into samples
the whole Maxwellian part but only a fraction $\bar \lambda$ of it,
where $$\bar \lambda\geq \max_{x}\{\lambda(x,t+\Delta t)\}.$$ As
discussed before, the reason for this is that we know that at the
subsequent relaxation step a fraction $\beta(x,t+\Delta t)$ of
samples will be discarded in each cell. Thus we need only
$$(1-\beta(x,t+\Delta t))\T_{\Delta t}(\beta(x,t)
M^p_f(x,v,t))=\lambda(x,t+\Delta t)\T_{\Delta t}(\beta(x,t)
M^p_f(x,v,t))$$ advected Maxwellian particles, which is guaranteed
if in any cell before advection we have at least $\bar
\lambda\beta(x,t) M^p_f(x,v,t)$ particles since
\[
\T_{\Delta t}(\bar \lambda\beta(x,t) M^p_f(x,v,t))=\bar \lambda
\T_{\Delta t}(\beta(x,t) M^p_f(x,v,t)).\]

This is of paramount importance since $\bar \lambda$ vanishes as
$\e/\Delta t\to 0$ and so the number of samples effectively used
by the hybrid method is a decreasing function of the ration
between the Knudsen number and the time step.

Starting from an initial data represented by particles a simple FSI
hybrid scheme for the solution of the BGK equation with $\lambda$
constant in space and time is described in the following algorithm
\begin{lgrthm}[FSI Hybrid Scheme]{~}
\label{al:H M/mm}
\begin{enumerate}
\item Compute the initial velocity and position of the particles $\{\nu_{j}^{0},j=1,..,N\}$ $\{p_{j}^{0},j=1,\ldots,N\}$
by sampling them from initial density $f_{0}(x,v)$. Set
$m^p=\int\int f_0(x,v)dx dv /N$.
\item Given a mesh $x_i$, $i=1,\ldots,L$ with grid size $\Delta x$, and an estimate of
the larger sample velocity $\nu_{max}=4\sqrt{2T_{max}}$, with
$T_{max}$ the maximum temperature, set $\Delta t^p=\Delta
x/\nu_{max}$.
\item Compute the initial values of the moments of the distribution
function in each cell $\varrho_i$, $(\varrho u)_i$, $E_i$,
$i=1,\ldots,L$.
\item Compute the larger time step allowed by the
deterministic macroscopic scheme $\Delta t_{D}$.
\item Set $\Delta t=min(\Delta t^p,\Delta t_{D})$.

\item Set $n=0$, $t=0$, $\lambda=e^{-\Delta t /\e}$, $\bar \lambda=\lambda$, $\beta_i=1-\bar\lambda$, $i=1,\ldots,L$.

\item While $t\leq t_f$ with $t_f$ the final chosen time.
\begin{enumerate}
\item Estimate the number of Maxwellian samples we need from $\bar \lambda \beta(x,t) M^p_f(x,v,t)$.
\begin{enumerate}
\item In each cell set $N_i^M=\IR(\bar \lambda\beta_i \rho_i^n/(m^p/\Delta x))$ and sample $N_i^M$
equilibrium particles from the Maxwellian with moments $(\rho
u)_i^n, E_i^n$.
\end{enumerate}
\item Perform the transport step keeping track of the particles that come from
the above sampling.
\begin{enumerate}
\item Transport all particles  \be
p_{j}^{n+1}=p_{j}^{n}+\nu^{n}_{j}\Delta t, \quad \forall\, j. \ee
\item Compute the moments $U_i^{p,n+1}$ and the number of particles $N_i^p$ in each
cell using only the advected particles not sampled from the
Maxwellian.
\item Solve the Euler equations for $U_i^{E,n}=\beta_i U_i^n$ and find
$U_i^{E,n+1}$.
\item Compute the new hybrid moments $U_i^{n+1}=U_i^{p,n+1}+U_i^{E,n+1}$
\end{enumerate}
\item Perform the relaxation step.
\begin{enumerate}
\item In each cell set $N_i^{k}=\IR(\lambda N_i^p)$ and discard
$N_i^p-N_i^k$ particles.
\item Compute the new number of particles in non equilibrium regime, in each cell
$N_i^{p}=N_i^k+M_i^p$, $i=1,\ldots,L$, with $M_i^p$ the
transported Maxwellian particles ($N_i^M$ transported). 
\item Compute the effective value $\lambda_i^p=(N_i^k+M_i^p)/(\varrho_i^{n+1}\Delta x/m^p)$.
\item Set $\beta_i=1-\lambda_i^p$.


\end{enumerate}
\item Set $t=t+\Delta t$, $n=n+1$ and compute the updated value of $\Delta
t$.
\end{enumerate}
end while
\end{enumerate}
\label{H M/m}
\end{lgrthm}

\begin{remark}~\rm
\begin{itemize}
\item In this simple version of the FSI hybrid method the value of
the equilibrium fraction fluctuates in each cell around the
constant value $\beta=1-\lambda$, thus it depends on $\Delta
t/\e$. We will see how to remove this limitation and make the
equilibrium fraction essentially independent of $\Delta t$ in the
optimized version of the FSI scheme. Note that fluctuations are
due to the fact that to have mass conservation during the
relaxation step we compute the effective value $\lambda_i^p$ and
set $\beta_i=1-\lambda_i^p$.
\item
 To avoid bias in the algorithm we used a
stochastic rounding $\IR(x)$ of a positive real number $x$ defined
as
\[
   \IR(x) = \left\{\begin{array}{lll}
                     {[x]}     & \mbox{with probability} & {[x]}+1-x, \\
                     {[x]} + 1 & \mbox{with probability} & x-{[x]},
                   \end{array}
            \right.
\]
where $[x]$ denotes the integer part of $x$.

\item In the fluid limit the numerical method is characterized
by the particular solver adopted to compute the solution of the
Euler equation. Thus the order of accuracy of the limiting scheme
is completely independent from the first order splitting procedure
used to solve the kinetic equation. This is an advantage compared
to the classical approach based on kinetic schemes which gives
limited accuracy in time. Extensions to higher order in the non
fluid regime are not trivial since we are limited to first order
accuracy in time by Theorem \ref{th:1}.
\end{itemize}

\end{remark}

\subsubsection{Matching moments}
In order to have a conservative scheme it is desirable that the
set of advected equilibrium samples satisfies \be \P(\T_{\Delta
t}(\beta(x,t)M^p_f(x,v,t)))=U^E(x,t+\Delta t),\label{eq:cons} \ee
namely the kinetic particles solution to the fluid equations in
one time step should match the direct solution to the limiting
fluid equations. Moreover, since the right hand side is not
affected by statistical sampling error, imposing (\ref{eq:cons})
will decrease the variance of the samples.

To this goal it is natural to use a moment matching approach
\cite{Cf}. This can be done by simple transformations of the
sample points. Given a set of samples $\nu_1,\ldots,\nu_{J}$ with
first two moments $\mu_1$ and $\mu_2$ and a better estimate $m_1$
and $m_2$ of the same moments we can apply the
transformation\cite{Cf, Trazzi}
\[
\nu_j^*=(\nu_j-\mu_1)/c+m_1\quad
c=\sqrt{\frac{\mu_2^2-\mu_1^2}{m_2-m_1^2}},\quad i=1,\ldots,J
\]
to get
\[
\frac1{J}\sum_{j=1}^J \nu_j^*=m_1,\qquad \frac1{J}\sum_{j=1}^J
(\nu_j^*)^2=m_2.
\]
Of course this renormalization is not possible for the mass density.
In fact, to keep the algorithm simple, we take the weight of each
particle $m^p$ equal and constant during the simulation and this
implies that we can have only integer multiples of such weights as
mass density values in each cell.

However thanks to the particular structure of the algorithm we can
perform also a matching procedure for the mass using the following
trick. After the transport of Maxwellian particles we need in each
cell, in order to perform the moment matching of order zero, a
number of particles given by
\[
M_i^p=\IR(\lambda \rho^{E}(x_i,t+\Delta t)/m^p).
\]
This can be done if we take $\bar\lambda$ large enough before
transport which guarantees that we have enough particles in each
cell. In this way the difference between the particles mass and
the Euler mass is below the mass of one single particle.

Next, to have exactly mass conservation
 we compute the effective values
$$\lambda_i^p=(M_i^p+N_i^k)/
(\rho^E(x_i,t+\Delta t)\Delta x/m^p+N_i^p), \quad
\beta^p(x_i,t+\Delta t)=1-\lambda_i^p,$$ used in the method. After
this we renormalize the transported equilibrium samples in each
space cell as described above so that they have the same momentum
and energy of the Euler solution.

Similarly one can apply a moment matching strategy when sampling
from the Maxwellian during the relaxation step. In this case, as
an alternative to the moment matching technique described above,
one can use the algorithm developed by Pullin \cite{Pullin79}.

Note that the whole method can be seen as a Monte Carlo scheme for
the BGK equation in which we try to reduce fluctuations
substituting the moments of the transported Maxwellian, computed
with particles, with the moments given by the solution of the
compressible Euler equation obtained with a deterministic
macroscopic scheme. Moreover, as described above, if we force the
equilibrium particles to follow the moments given by the fluid
equations we can reinterpret the algorithm as a fluid-dynamic
guided Monte Carlo scheme.

\subsection{Optimal FSI Methods}
\label{increasingb} The method just described does not take into
account the possibility to optimize the equilibrium fraction by
increasing its value in time and make it independent on the choice
of the time step. In fact at each time step the equilibrium
structure is entirely lost and the new fraction of equilibrium is
only given by the relaxation step (see \ref{eq:rel}). However, in
principle, it is possible to recover some information from the
transported local Maxwellian although we know it through samples
rather than analytically. We recall, in fact, that we do not get
any microscopic information from $U^E(x,t)$ which corresponds to
the solution of the Euler equation with a macroscopic numerical
scheme. In the sequel, we will propose a method to optimize the
equilibrium fraction $\beta(x,t)$ after the transport step. We
start describing the generalization of the hybrid method once this
optimization has been achieved.

Thanks to Definition 1 we can define the velocity dependent
optimal equilibrium fraction as the ratio of the transported
Maxwellian at time $t$ respect to the new local Maxwellian at time
$n+1$
\begin{equation*}w^c(x,v,t+\Delta t)=\left\lbrace
\begin{array}{ll}
\displaystyle\frac{\T_{\Delta t}(\beta(x,t)
M_{f}(x,v,t))}{M^H_{f}(x,v,t+\Delta t)}, & \T_{\Delta t}(\beta(x,t)
M_{f}(x,v,t)) \leq M^H_{f}(x,v,t+\Delta t),\\
1, & \T_{\Delta t}(\beta(x,t)
M_{f}(x,v,t)) \geq M^H_{f}(x,v,t+\Delta t),\\
\end{array}
\right.
\end{equation*}
and the optimal equilibrium fraction as \be \beta^c(x,t+\Delta
t)=\min_v\{w^c(x,v,t+\Delta t)\}. \label{eq:def11}\ee
This value can be considered optimal, in the sense that it is the
maximum allowed value for which we have a decomposition like \be
\T_{\Delta t}(\beta(x,t) M_f(x,v,t))= \tilde M_f(x,v,t+\Delta
t)+\beta^c(x,t+\Delta t)M_f^H(x,v,t+\Delta t)\label{eq:Mtilde} \ee
with $\tilde M_f(x,v,t+\Delta t)\geq 0$. Clearly similar
decompositions hold true for any fraction of equilibrium below the
optimal one.

Suppose, for simplicity, that $\beta^c(x,t)=0$ at the beginning of
our computation, it follows that the method reads in the same way
from equation (\ref{eq:id}) to equation (\ref{eq:fin}). Now, given
an estimation for $\beta^c(x,t+\Delta t)$ the next relaxation step
reads as \bea \nonumber \R_{\Delta t}(\T_{\Delta
t}(f(x,v,t)))&=&\lambda\T_{\Delta t}(f(x,v,t))
+(1-\lambda) M^H_f(x,v,t+\Delta t)\\
\nonumber &=& \lambda(\T_{\Delta t}((1-\beta(x,t))
f^p(x,v,t))+\T_{\Delta t}(\beta(x,t) M_f(x,v,t)))\\\nonumber &+& (1-\lambda)M^H_f(x,v,t+\Delta t)\\
\nonumber &=& \lambda(\T_{\Delta t}((1-\beta(x,t))
f^p(x,v,t))+\beta^c(x,t+\Delta t) M^H_f(x,v,t+\Delta t)\\
\nonumber &+&\tilde M_f(x,v,t+\Delta t))+(1-\lambda) M^H_f(x,v,t+\Delta t)\\
\nonumber &=&(1-\beta(x,t+\Delta t)) f^p(x,v,t+\Delta
t)+\beta(x,t+\Delta t) M^H_f(x,v,t+\Delta t),\eea with \be
\beta(x,t+\Delta t)=1-\lambda(1-\beta^c(x,t+\Delta t))\ee and \be
f^p(x,v,t+\Delta t)=\frac{\T_{\Delta
t}((1-\beta(x,t))f^p(x,v,t))+\tilde M_f(x,v,t+\Delta
t)}{1-\beta^c(x,t+\Delta t)}.\ee



In order to sample from the distribution $\tilde M_f(x,v,t+\Delta
t)$, which is obtained as a difference of two distribution
functions, see (\ref{eq:Mtilde}), we can sample particles, exactly
as in the previous section, from the transported Maxwellian and
then apply an acceptance rejection technique that reads
\begin{lgrthm}[Acceptance-Rejection Sampling]{~}\\ \newline
do $i=1,N$ with $N$ number of particles to be sampled
\begin{enumerate}
\item Select randomly one particle from the distribution $\T_{\Delta t}(\beta(x,t)
M_{f}(x,v,t))$;
\item with probability $\displaystyle 1-\frac{\beta^c(x,t+\Delta t)M_f^H(x,v,t+\Delta
t)}{\T_{\Delta t}(\beta(x,t) M_{f}(x,v,t))}$ keep the particle.
\end{enumerate}
\label{a/r}
\end{lgrthm}
In the above algorithm particles can be taken more than once, in
other words the sampling is not exclusive. Finally a moment
matching strategy similar to the one described in the previous
section can be used in such a way that the equation \be \P(\tilde
M_f(x,v,t+\Delta t))=U^E(x,t+\Delta t)-\beta^c(x,t+\Delta t)\tilde
U^H(x,t+\Delta t),\label{eq:mm2}\ee is satisfied exactly.

The major problem we have to face when practically evaluating
$\beta^c(x,t+\Delta t)$ is that $M_f^H(x,v,t+\Delta t)$ is known
analytically while $\T_{\Delta t}(\beta(x,t) M_{f}(x,v,t))$ is known
only through samples. From a numerical viewpoint when approximating
$\beta^c(x,t+\Delta t)$ we want to avoid overestimates since these
may produce unphysical solutions.
In the following description, to simplify notations, we restrict
to 1-D in velocity and physical space, extensions of the methods
to multidimensional cases are straightforward. Our goal is to find
an estimation of $\beta^c(x,t)$ given by (\ref{eq:def11}). Without
loss of generality we assume that at each point $x$ there exist a
velocity $v$ such that
$$\T_{\Delta t}(\beta(x,t)M_{f}(x,v,t))\leq M^H_{f}(x,v,t+\Delta
t).$$ In fact, for those space points $x$ where the above
assumption is not satisfied we simply have $\beta^c(x,t+\Delta
t)=1$.

The first and simplest method consists in measuring the departure
from equilibrium reconstructing the transported Maxwellian from
samples. In order to do that we need a grid in velocity space and
a loop over the particles inside each spatial cell.
We omit here the details of the different reconstruction methods
that can be used, we refer to \cite{pa-inria} (and the references
therein) for the technical aspects.

Once we have reconstructed $\T_{\Delta t}(\beta(x_i,t)
M_{f}(x_i,v,t))$, with $\{x_i\}_{i\in I}$ a mesh of the physical
space, we are able to determine the quantity
\begin{equation}\beta^c(x_i,t+\Delta t)=\min_{v}\left\lbrace
\displaystyle\frac{\T_{\Delta t}(\beta(x_i,t)
M_{f}(x_i,v,t))}{M^H_{f}(x_i,v,t+\Delta t)}\right\rbrace.
\label{eq:def2}
\end{equation}
This method presents several drawbacks. The reconstruction of the
distribution function from samples increase the computational
cost, moreover a small number of particles inside a cell, which is
quite common in applications, gives large fluctuations and this
turns in an imprecise estimate of $\beta^c(x,t)$.

A better way to estimate the equilibrium fraction $\beta^c(x,t)$
after the transport is based on the analysis of a deterministic
transport of the Maxwellian part. Again we introduce a grid in
space.
We consider the following scheme for the transport of the
Maxwellian fraction \bea \nonumber \frac{\hat
M_{f,i}^{n+1}(v)-M_{f,i}^{n}(v)}{\Delta
t}+v\frac{M^{n}_{f,i}(v)-M_{f,i-1}^{n}(v)}{\Delta x}&=&0, \quad v\geq0,\\
\\
\label{eq:upwind} \nonumber \frac{\hat
M_{f,i}(v)^{n+1}-M_{f,i}^{n}(v)}{\Delta t}+v
\frac{M^{n}_{f,i+1}(v)-M_{f,i}^{n}(v)}{\Delta x}&=&0, \quad v<0,
\eea where $M_{f,i}^{n}(v)\approx M_f(x_i,v,t^n)$, $\hat
M_{f,i}^{n+1}(v)\approx M_f(x_i,v,t^{n+1})$ and $\Delta x$ is the
mesh size in space. We put an hat on the transported Maxwellian to
distinguish it with respect to the local Maxwellian at time
$t+\Delta t$, which is accordingly to the notations,
$M_{f,i}^{n+1}(v)$. The scheme described above is a simple first
order upwind for the Maxwellian transport. Of course to
effectively perform the computation it is necessary to truncate
the Maxwellian in order to obtain finite values for the velocity
and time step larger than zero. Typically this truncation leads to
several problems which are common in numerical methods for kinetic
equations (see for example \cite{Mieussens, Russo1}). Here we are
only interested to estimate the departure from the equilibrium of
the transported Maxwellian and for that scope we choose a bound
for the velocity space in such a way that no additional time step
restrictions are imposed. Solving Eqs. (\ref{eq:upwind}) we obtain
\bea \nonumber\hat M_{f,i}^{n+1}(v)&=&\left(1-\frac{v\Delta
t}{\Delta x}\right)M_{f,i}^{n}(v)+\frac{v\Delta t}{\Delta
x}M_{f,i-1}^{n}(v),
\quad v\geq0,\\
\\ \nonumber\hat M_{f,i}^{n+1}(v)&=&\left(1+\frac{v\Delta t}{\Delta
x}\right)M_{f,i}^{n}(v)-\frac{v\Delta t}{\Delta
x}M_{f,i+1}^{n}(v),\quad v<0.\eea Note that, since $|v|\Delta
t\leq \Delta x$, the updated function $\hat M_{f,i}^{n+1}(v)$ is a
convex combination of the local Maxwellian in the cells $i$ and
$i-1$ for positive velocities and in the cells $i$ and $i+1$ for
negative velocities.

Now, ignoring the error introduced by the truncation in velocity,
in each cell the equilibrium fraction $\beta^c(x_i,t+\Delta t)$
satisfies \be \beta^c(x_i,t+\Delta
t)=\min_{v}\left\{\frac{\T_{\Delta
t}(\beta(x_i,t)M_{f}(x_i,v,t))}{M^H_{f}(x_i,v,t+\Delta
t)}\right\}=\min_{v}\left\{\frac{\hat
M_{f,i}^{n+1}(v)}{M^{H,n+1}_{f,i}(v)}\right\}+O(\Delta
t^2).\label{eq:ratio}\ee In the general case a numerical method is
required to compute the minimum on the right hand side. This
operation can be expensive since it has to be done at each time
step and in each spatial cell. However we can restrict ourselves
to a lower estimate of $\beta^c(x_i,t+\Delta t)$ (to avoid an
overestimate of the equilibrium fraction) and we can choose
instead of the minimum a lower bound for (\ref{eq:ratio}) using
the convexity property of the scheme. That value can be estimated
observing that

\begin{eqnarray}
\nonumber   \min_{v\geq 0}\left\{\frac{\hat
M_{f,i}^{n+1}(v)}{M^{H,n+1}_{f,i}(v)}\right\}&\geq&\min\left\{\min_{v\geq
0}\left\{\frac{M_{f,i}^{n}(v)}{M^{H,n+1}_{f,i}(v)}\right\},
\min_{v\geq 0}\left\{\frac{M_{f,i-1}^{n}(v)}{M^{H,n+1}_{f,i}(v)}\right\}\right\}  \\
&=& \beta_R^c(x,t+\Delta t),\label{eq:min}
\end{eqnarray}

\begin{eqnarray}
\nonumber
   \min_{v<0}\left\{\frac{\hat
M_{f,i}^{n+1}(v)}{M^{H,n+1}_{f,i}(v)}\right\}&\geq&\min\left\{\min_{v<0}\left\{\frac{M_{f,i}^{n}(v)}{M^{H,n+1}_{f,i}(v)}\right\},
\min_{v<0}\left\{\frac{M_{f,i+1}^{n}(v)}{M^{H,n+1}_{f,i}(v)}\right\}\right\} \\
   &=&\beta^c_L(x,t+\Delta
t),\label{eq:min1}
\end{eqnarray}
and setting \be \beta^c(x,t+\Delta t)=\min\{\beta^c_R(x,t+\Delta
t),\beta^c_L(x,t+\Delta t)\}\label{eq:min2}\ee where the minimum
of the ratios in (\ref{eq:min})-(\ref{eq:min1}) can be computed
exactly, being the ratios of Maxwellian functions.

An algorithm that can be used to implement the optimized FSI
method for the solution of the BGK equation, in which for
simplicity $\lambda$ is constant, is the following

\begin{lgrthm}[Optimized FSI Hybrid Scheme]{~}
\begin{enumerate}
\item Compute the initial velocity and position of the particles $\{\nu_{j}^{0},j=1,..,N\}$ $\{p_{j}^{0},j=1,\ldots,N\}$
by sampling them from initial density $f_{0}(x,v)$. Set
$m^p=\int\int f_0(x,v)dx dv /N$.
\item Given a mesh $x_i$, $i=1,\ldots,L$ with grid size $\Delta x$, and an estimate of
the larger sample velocity $\nu_{max}=4\sqrt{2T_{max}}$, with
$T_{max}$ the maximum temperature, set $\Delta t^p=\Delta
x/\nu_{max}$.
\item Compute the initial values of the moments of the distribution
function in each cell $\varrho_i$, $(\varrho u)_i$, $E_i$,
$i=1,\ldots,L$.
\item Compute the larger time step allowed by the
deterministic macroscopic scheme $\Delta t_{D}$.
\item Set $\Delta t=min(\Delta t^p,\Delta t_{D})$.

\item Set $n=0$, $t=0$, $\lambda=e^{-\Delta t /\e}$, $\bar \lambda=\lambda$, $\beta_i=1-\bar\lambda$, $i=1,\ldots,L$.

\item While $t\leq t_f$ with $t_f$ the final chosen time.
\begin{enumerate}
\item Estimate the number of Maxwellian samples we need from $\bar \lambda \beta(x,t) M^p_f(x,v,t)$.
\begin{enumerate}
\item In each cell set $N_i^M=\IR(\beta_i \rho_i^n/(m^p/\Delta x))$ and sample $N_i^M$
equilibrium particles from the Maxwellian with moments $(\rho
u)_i^n, E_i^n$.
\end{enumerate}
\item Perform the transport step keeping track of the particles that come from
the above sampling.
\begin{enumerate}
\item Transport all particles  \be
p_{j}^{n+1}=p_{j}^{n}+\nu^{n}_{j}\Delta t, \quad \forall\, j. \ee
\item Compute the moments $U_i^{p,n+1}$ and the number of particles $N_i^p$ in each
cell using only the advected particles not sampled from the
Maxwellian.
\item Solve the Euler equations for $U_i^{E,n}=\beta_i U_i^n$ and find
$U_i^{E,n+1}$.
\item Compute the new hybrid moments $U_i^{n+1}=U_i^{p,n+1}+U_i^{E,n+1}$
\end{enumerate}
\item Compute the optimal equilibrium fraction $\beta^{c,n+1}_{i}$ as
described in (\ref{eq:min})-(\ref{eq:min2}).
\item Perform the relaxation step.
\begin{enumerate}
\item In each cell set $N_i^{k}=\IR(\lambda N_i^p)$ and discard
$N_i^p-N_i^k$ particles.
\item In each cell sample $\tilde N_i^M$ particles from the distribution $\tilde M^{n+1}_{f,i}(v)$with
the acceptance-rejection technique described in Algorithm 2.
\item Apply the moment matching technique to the $\tilde N_i^M$ particles in order to satisfy (\ref{eq:mm2}).
\item Compute the new number of particles in non equilibrium regime, in each cell
$N_i^{p}=N_i^k+\tilde N_i^M$, $i=1,\ldots,L$. 
\item Compute the effective equilibrium fraction $\beta_i=1-(N_i^k+\tilde N_i^M)/(\varrho_i^{n+1}\Delta x/m^p)$.

\end{enumerate}
\item Set $t=t+\Delta t$, $n=n+1$ and compute the updated value of $\Delta
t$.
\end{enumerate}
end while
\end{enumerate}

\end{lgrthm}

\begin{remark}\rm
\begin{itemize}
\item We emphasize that the first order upwind method to compute
the deterministic transport of Maxwellians is never used in
practice. It serves us only as an approximation strategy in order
to compute a lower bound for the optimal equilibrium fraction
$\beta^c(x,t+\Delta t)$. In this sense it is worth to notice that
the additional first order dissipation introduced by this
upwinding produces additional smearing and in principle, close to
discontinuities, can produce overestimates of the equilibrium
fraction when computed from (\ref{eq:ratio}). Besides
computational efficiency this is an additional motivation to use a
lower bound for that value.
\item The hybrid composition of the solution in the final method
does not depend on the time step $\Delta t$ but only on $\e$.
Note, however, that small time steps, below the CFL condition of
the deterministic Euler solver, may increase the computational
cost. To reduce this effect one can use different time steps in
the kinetic and the Euler solver and perform the hybridization and
matching only at intermediate steps. This strategy can be used,
for example, where there is the need to resolve small scales at
the Navier-Stokes level or in boundary layer effects.
\end{itemize}

\end{remark}

\section{Implementation and numerical tests}
In principle any finite volume or finite difference numerical
scheme can be used to solve the compressible Euler equations in
our hybrid method. In the sequel we will use a second order finite
volume MUSCL type relaxed scheme (see \cite{JX} for details).

In the next paragraphs we analyze the performances of the fluid
solver independent hybrid schemes in comparison with a classical
Monte Carlo method (MCM) for several one-dimensional problems with
different Knudsen numbers ranging from $\varepsilon=10^{-2}$ to
$\varepsilon=10^{-5}$.

As a reference solution we used a deterministic discrete velocity
model (DVM) for the BGK equation for all tests (see
\cite{Mieussens} for details). We use the shorthand FSI, FSI1 and
TVD to denote the simple FSI method, the optimal FSI method and
the second order in space MUSCL Euler solver respectively.

\subsection{Accuracy test}
Because our aim is to compare the differences in accuracy and
computational time between the different methods first we have
considered an accuracy test with a a periodic smooth solution. We
compare the results of the kinetic-solver based hybrid methods
developed in \cite{dimarco3} to the new independent fluid solver
schemes.

We report the total $L_1$ norm of the errors for the conserved
quantity $\varrho$, $u$, and $T$ as the computational times by
considering a problem with the following initial data \be
\nonumber \varrho(x,0)=1+a_\varrho \sin\frac{2\pi x}{L} \ee \be
u(x,0)=1.5+a_u \sin\frac{2\pi x}{L} \ee \be \nonumber
E(x,0)=2.5+a_T \sin\frac{2\pi x}{L} \ee where we set \be \nonumber
a_\varrho=0.3 \ \ a_u=0.1 \ \ a_E=1. \ee The equations are
integrated for $t \in [0,5\times 10^{-2}]$ using $200$ space cell.

\begin{figure}
\begin{center}
\includegraphics[scale=0.42]{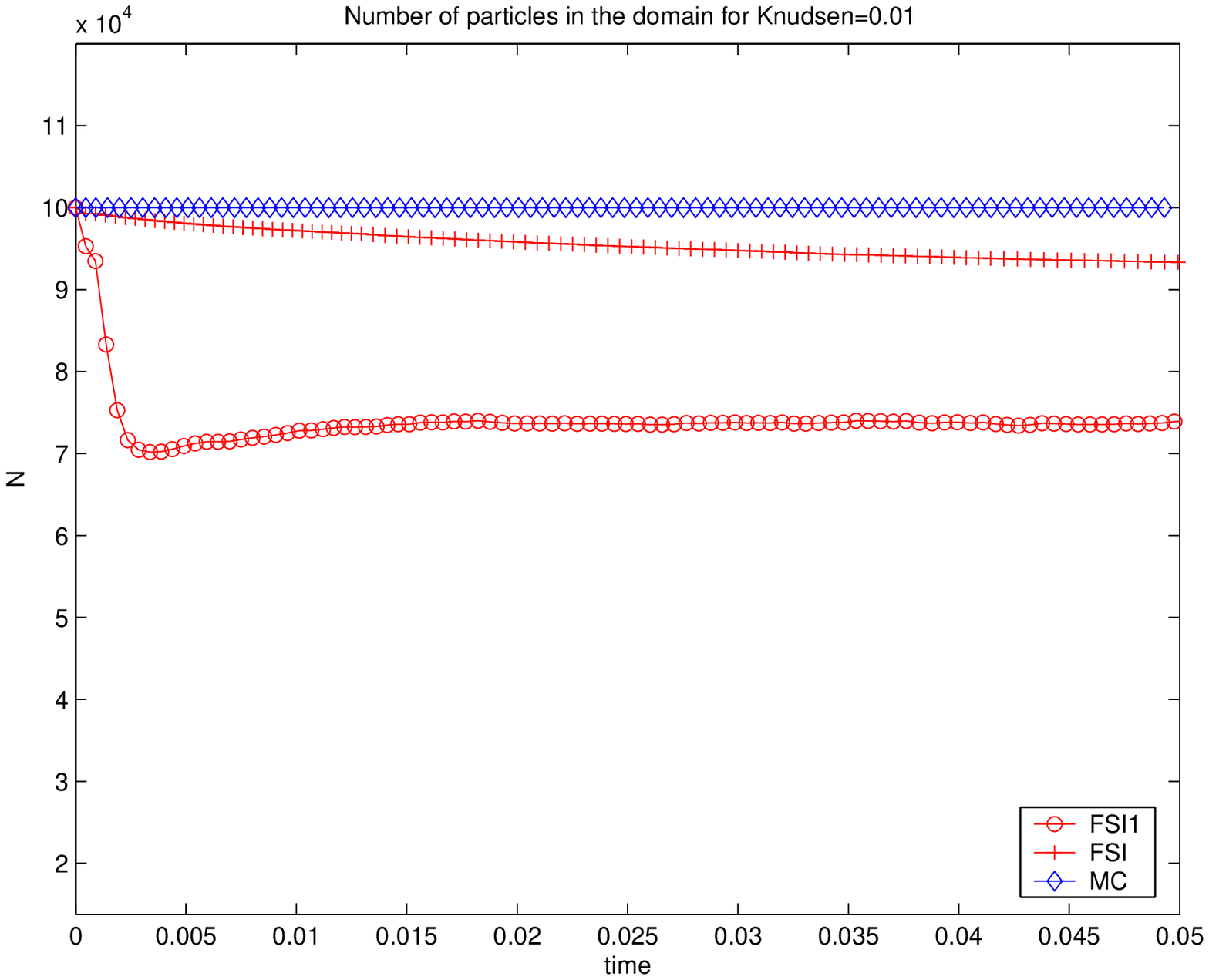}
\includegraphics[scale=0.42]{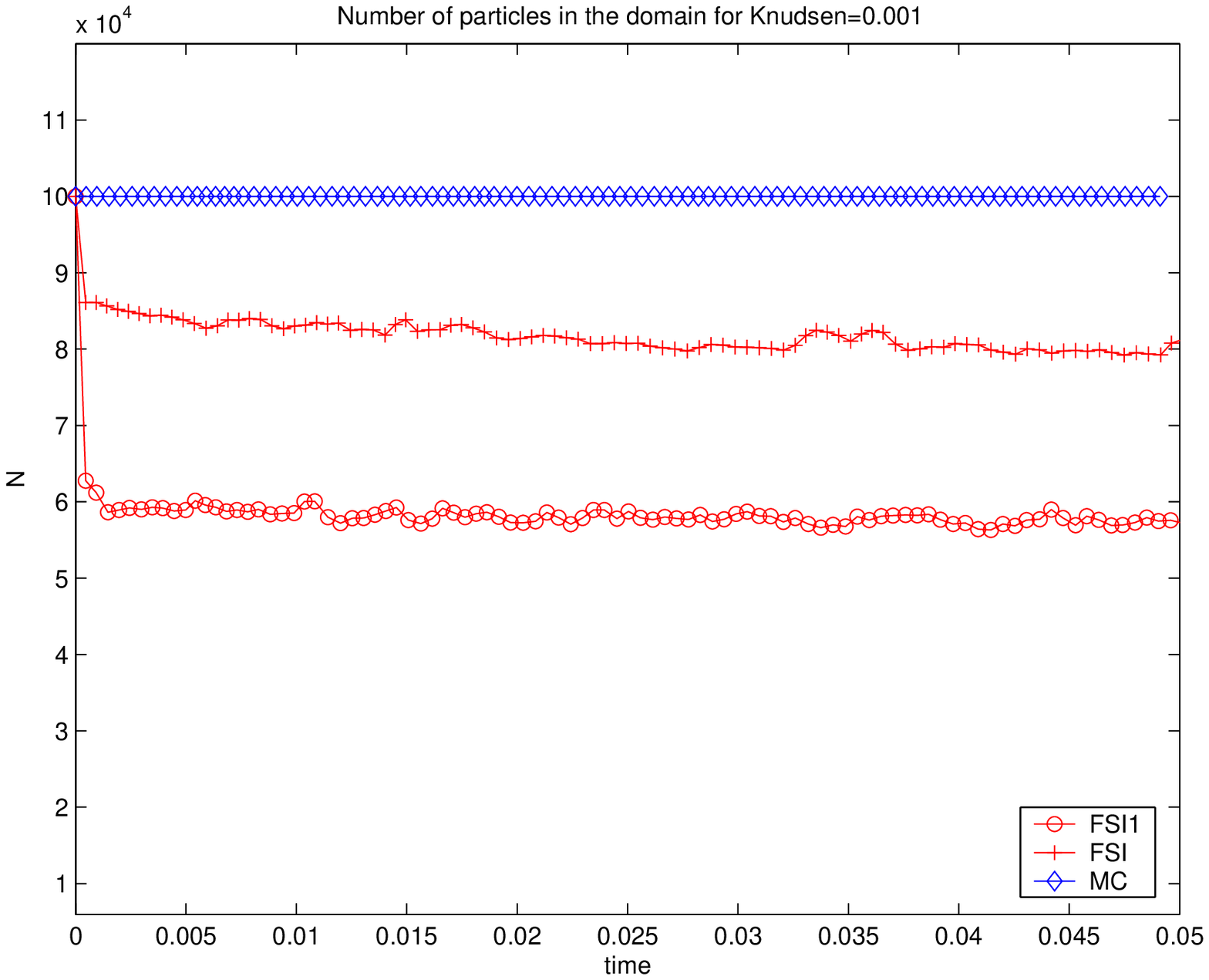}
\includegraphics[scale=0.42]{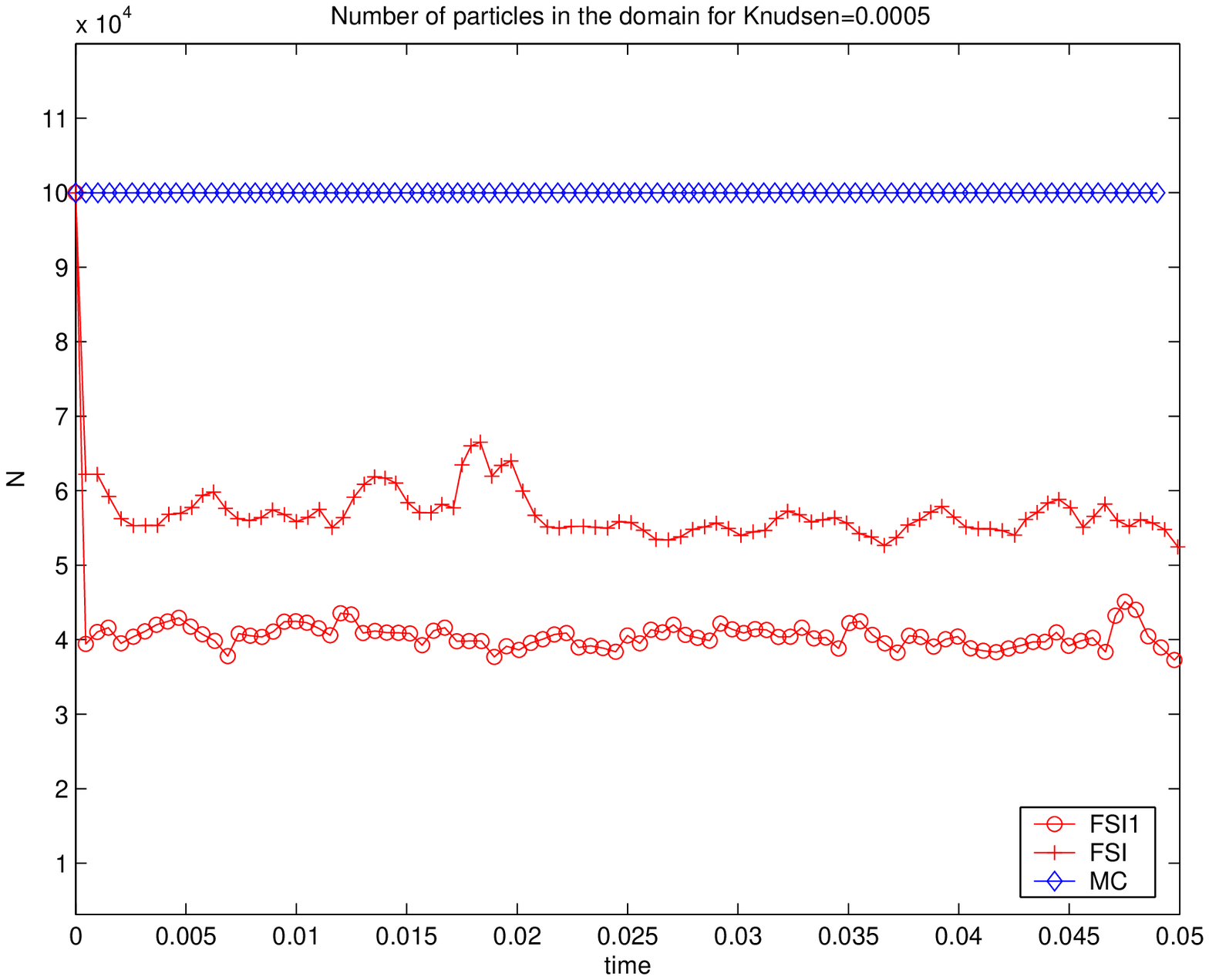}
\includegraphics[scale=0.42]{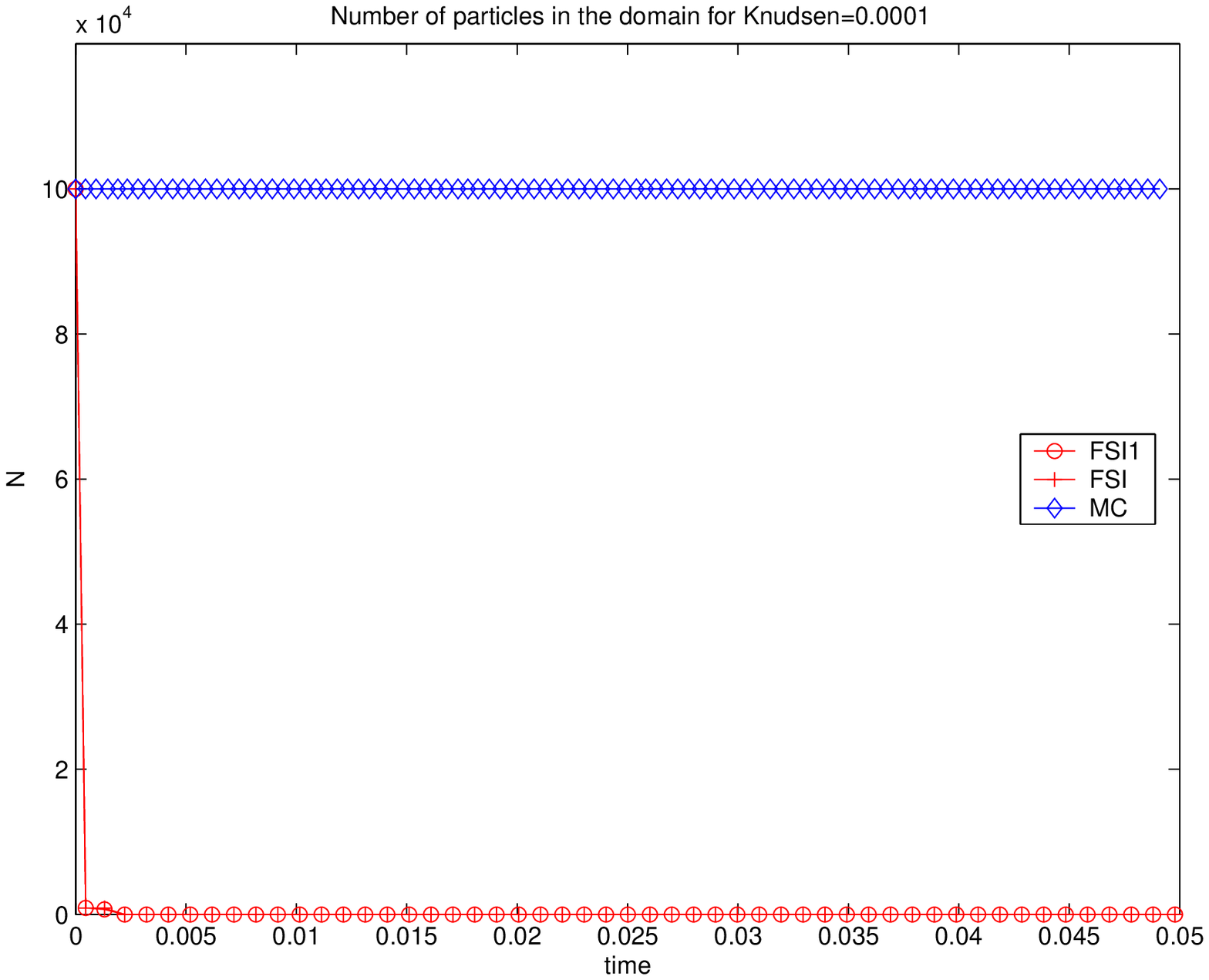}
\caption{Accuracy test. Number of particles in time inside the
computational domain for FSI1, FSI and MC schemes. Knudsen numbers
$\varepsilon=10^{-2}$ (top left) $\varepsilon=10^{-3}$ (top right)
$\varepsilon=5\times 10^{-4}$ (bottom left) and
$\varepsilon=10^{-4}$ (bottom right).} \label{US7}
\end{center}
\end{figure}

In order to make a fair comparison with the previous schemes named
BHM, BHM1, BCHM (see \cite{dimarco3} for details and parameters
settings of the methods) we use at the beginning $1500$ particles
for cell with the same time step of the Boltzmann-BGK schemes.
Then, because in general FSI type schemes allow larger time steps
and the moment matching techniques produce lower fluctuations we
repeat the computation with $N=500$ and the time step prescribed
by the Monte Carlo method. We remark that, while estimating the
time differences between the FSI schemes and MCM is quite easy (in
fact the FSI methods are based on the MCM for the kinetic part),
the same comparison with another kinetic solver such as DVM or
BHM, is not straightforward due to the several possible choices
involved in such schemes (for example the way the velocity domain
is truncated and the type of solver chosen for the space
derivatives). For this reason we stress that the simulations
times, reported in Table 1, are just indicative. It is clear that
the hybrid methods here developed represent a strong improvement
with respect to the previous schemes as well as to the classical
Monte Carlo method in terms of computational time. In Figure
\ref{US7} we report the total number of particles used by the
different algorithms. Note how computational time and fluctuations
reduce dramatically when the Knudsen number diminishes.

\begin{table}
\begin{center}
\begin{tabular}{|l|r|c|c|c|}
\hline
 & $\varepsilon=10^{-2}$ & $\varepsilon=10^{-3}$ & $\varepsilon=10^{-4}$ & $\varepsilon=10^{-5}$ \\
\hline
MCM N=1500 & 23 sec & 25 sec & 27 sec & 26 sec  \\
\hline
BHM N=1500 & 35 sec & 25 sec & 22 sec & 22 sec  \\
\hline
BHM1 N=1500 & 34 sec & 20 sec & 19 sec & 20 sec  \\
\hline
BCHM N=1500 & 15 sec & 11 sec & 17 sec & 21 sec  \\
\hline
FSI N=1500 & 25 sec & 22 sec & 3 sec & 0.6 sec  \\
\hline
FSI1 N=1500 & 18 sec & 17 sec & 2 sec & 0.6 sec  \\
\hline
FSI N=500 & 9 sec & 8 sec & 0.4 sec & 0.3 sec  \\
\hline
FSI1 N=500 & 7 sec & 6 sec & 0.4 sec & 0.3 sec  \\
\hline
\end{tabular}
\end{center} \caption{Accuracy test. Computational times for FSI and FSI1 with
two different initial numbers of particles $N=1500$ and $N_1=500$
compared to the previous hybrid methods (see \cite{dimarco3}) for
different values of the Knudsen number.}
\end{table}

\begin{table}
\begin{center}
\begin{tabular}{|l|r|c|c|c|} \hline
 & $\varepsilon=10^{-2}$ & $\varepsilon=10^{-3}$ & $\varepsilon=5\times 10^{-4}$ & $\varepsilon=10^{-4}$ \\
\hline
MCM & 5.494 & 5.786  & 5.153 &  5.184 \\
\hline
FSI &  5.545 & 3.926 & 3.067 & 0.268  \\
\hline
FSI1 & 4.588 & 3.406 & 2.451 & 0.243 \\
\hline
\end{tabular}
\end{center} \caption{Accuracy test. $L_1$ norm of the error (in units
of $10^{-2}$) for the density with respect to different values of
the Knudsen number $\varepsilon$.}\label{tab1}
\end{table}

The results for the relative $L_1$ errors are reported in Tables
\ref{tab1}, \ref{tab2}, \ref{tab3}, respectively for density, mean
velocity and temperature for the FSI, FSI1 and the Monte Carlo
schemes. In all the methods we used $N=200$ particles for cell and
the moment matching techniques. We notice that the hybrid methods
have approximately the same accuracy of the Euler solver for small
Knudsen numbers and the same accuracy of the Monte Carlo method
for large Knudsen while to intermediate values correspond
intermediate behaviors. From the above results it is clear how the
performances of the hybrid schemes are better than the ones of a
Monte Carlo method, moreover FSI1 gives in general better results
respect to FSI both in term of computational time and accuracy for
almost all regimes.

\begin{table}
\begin{center}
\begin{tabular}{|l|r|c|c|c|}
\hline
 & $\varepsilon=10^{-2}$ & $\varepsilon=10^{-3}$ & $\varepsilon=5\times 10^{-4}$ & $\varepsilon=10^{-4}$ \\
\hline
MCM & 6.565 & 5.437 & 5.338 & 6.035 \\
\hline
FSI &  4.802 & 4.401 & 3.264 & 0.641 \\
\hline
FSI1 & 5.135 & 4.102 & 2.848 & 0.610 \\
\hline
\end{tabular}
\end{center} \caption{Accuracy test. $L_1$ norm of the error (in units
of $10^{-2}$) for the mean velocity with respect to different
values of the Knudsen number $\varepsilon$.} \label{tab2}
\end{table}

\begin{table}
\begin{center}
\begin{tabular}{|l|r|c|c|c|}
\hline
 & $\varepsilon=10^{-2}$ & $\varepsilon=10^{-3}$ & $\varepsilon=5\times 10^{-4}$ & $\varepsilon=10^{-4}$ \\
\hline
MCM & 6.762 & 7.611 & 7.578 & 7.316 \\
\hline
FSI &  7.007 & 6.022 & 4.500 & 0.641  \\
\hline
FSI1 & 6.662 & 4.939 & 3.773 & 0.598 \\
\hline
\end{tabular}
\end{center} \caption{Accuracy test. $L_1$ norm of the error (in units
of $10^{-2}$) for the temperature with respect to different values
of the Knudsen number $\varepsilon$.}\label{tab3}
\end{table}

\subsection{1-D Unsteady shock}
\begin{figure}
\begin{center}
\includegraphics[scale=0.40]{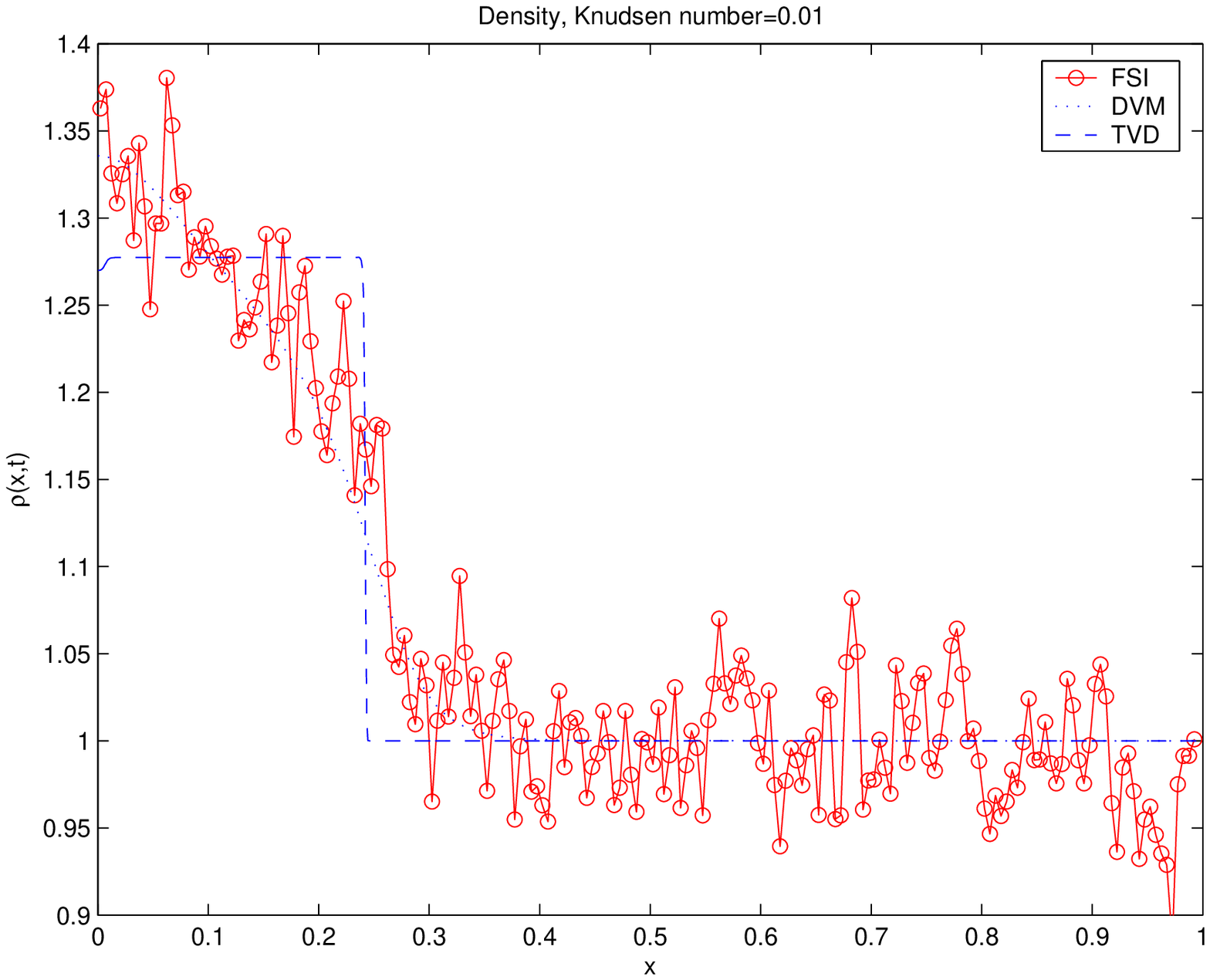}
\includegraphics[scale=0.40]{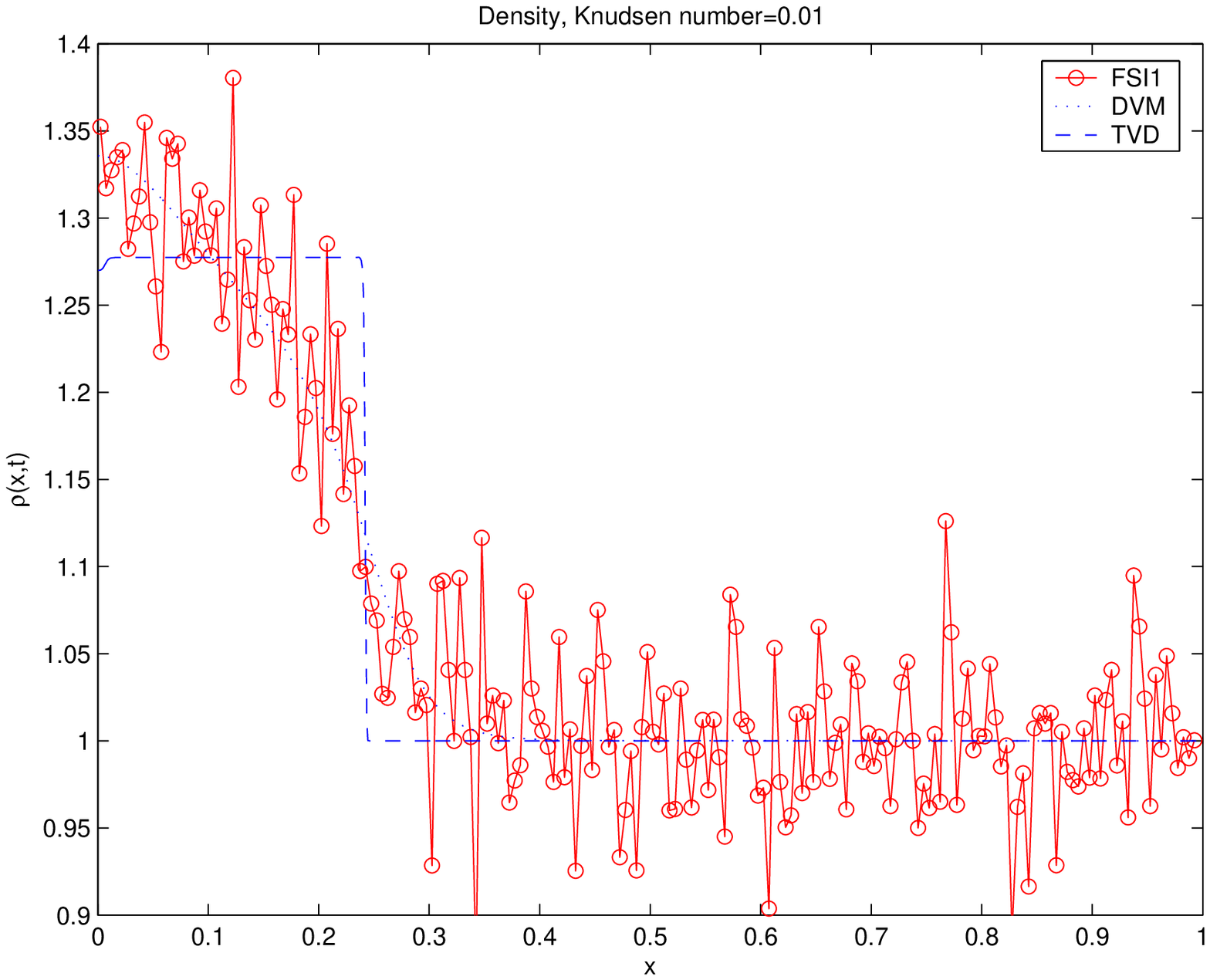}
\includegraphics[scale=0.40]{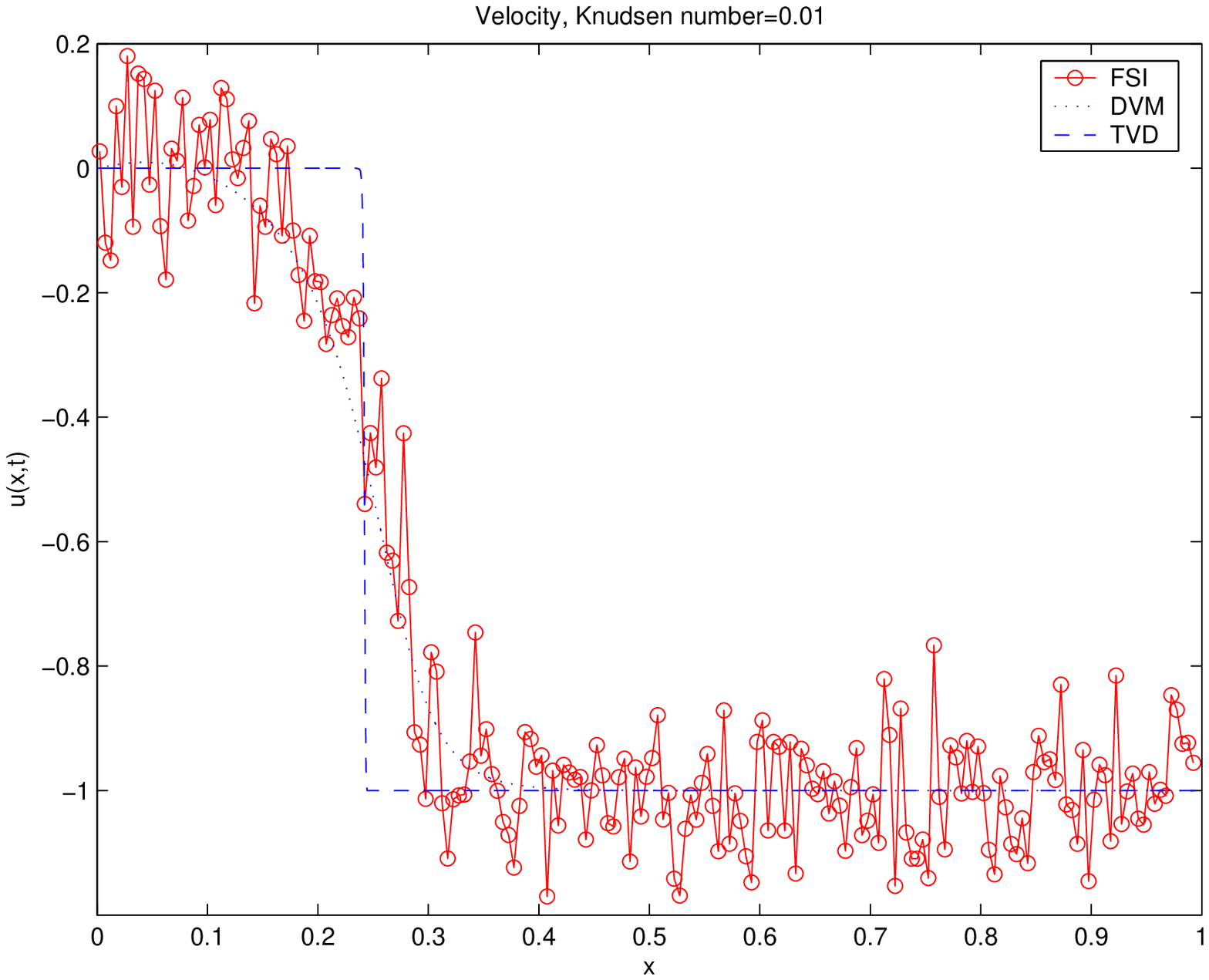}
\includegraphics[scale=0.40]{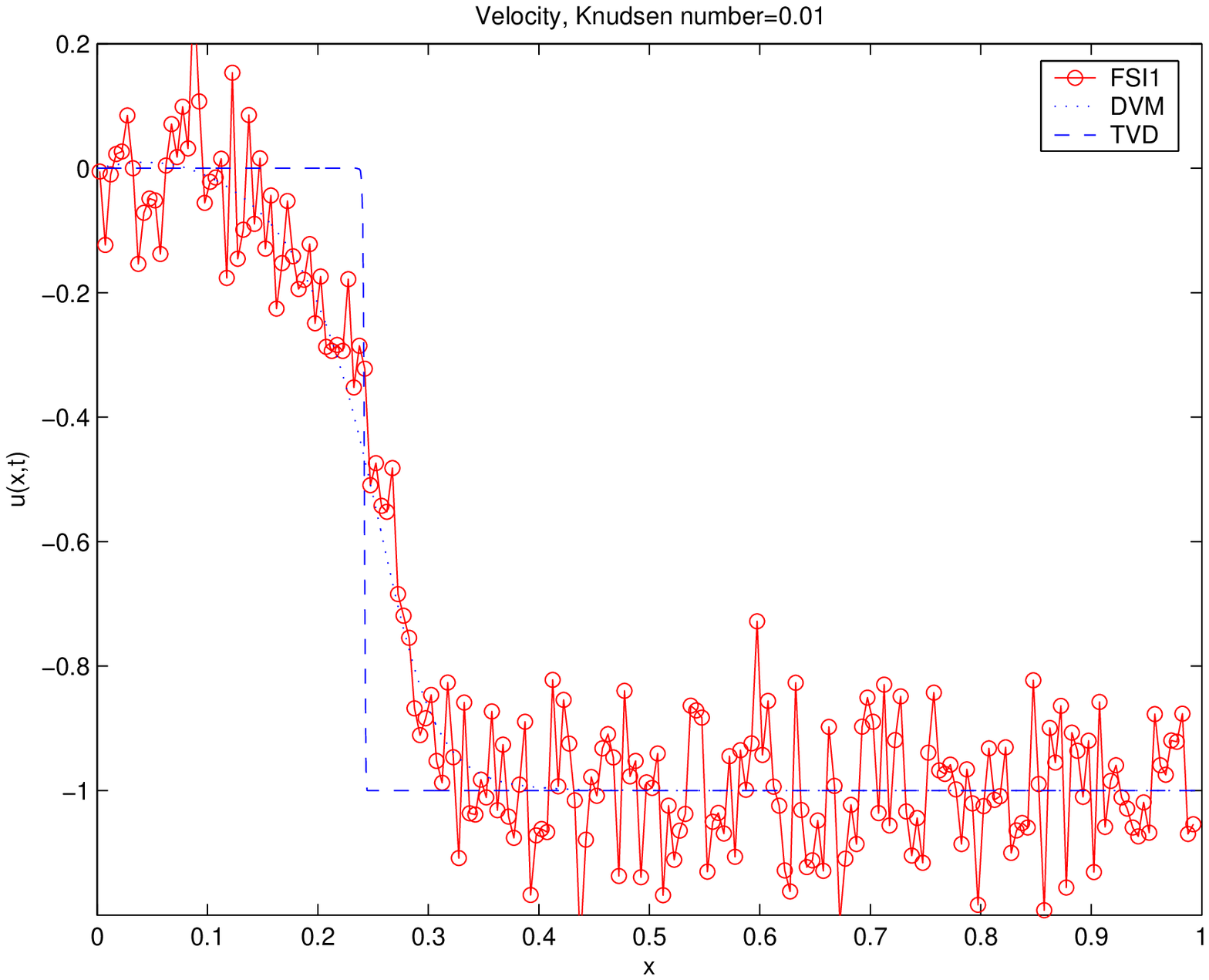}
\includegraphics[scale=0.40]{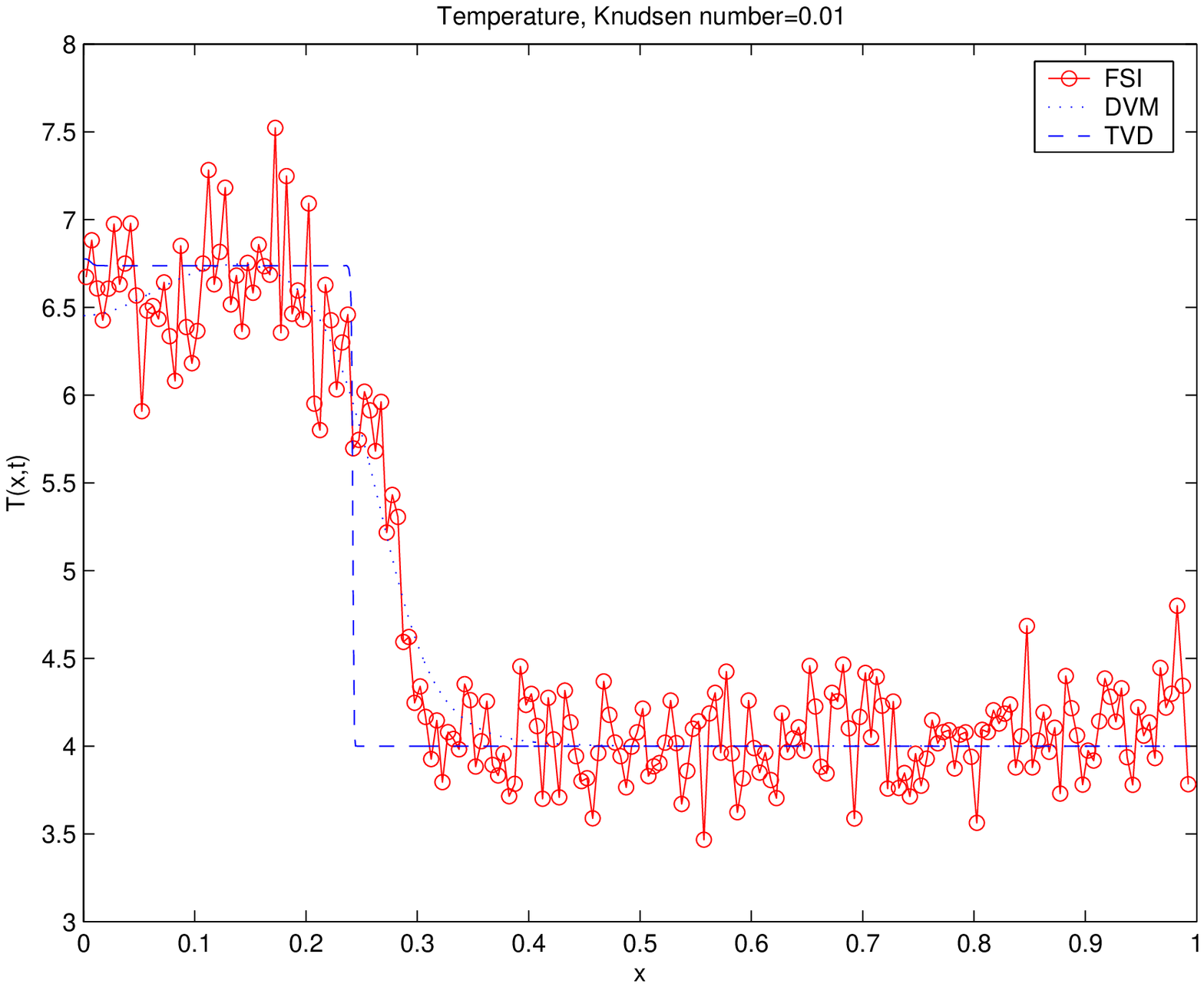}
\includegraphics[scale=0.40]{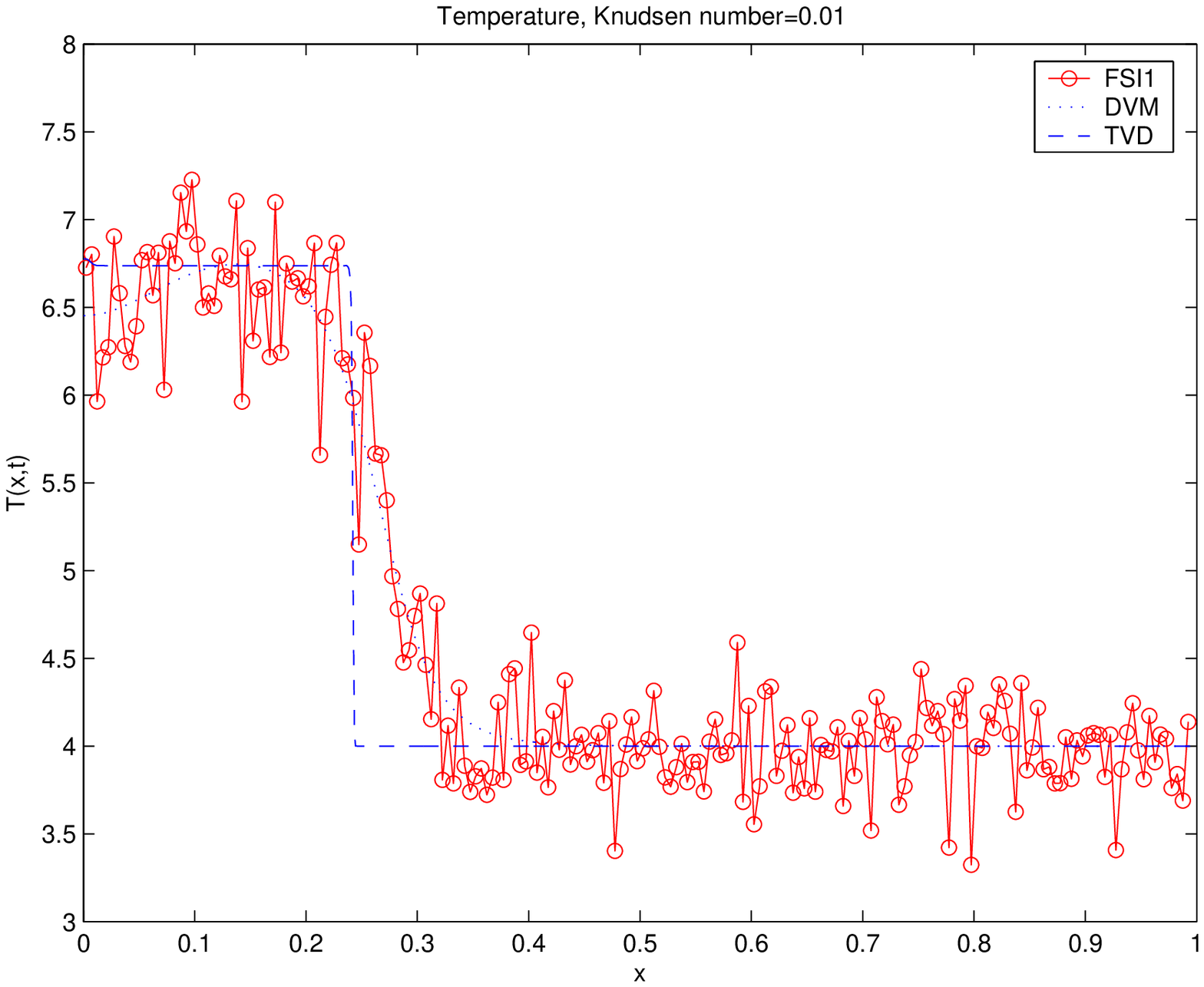}
\caption{Unsteady Shock: $\varepsilon=10^{-2}$. Solution at
$t=0.065$ for FSI (left) FSI1 (right). From top to bottom density,
mean velocity and temperature.} \label{US1}
\end{center}
\end{figure}

\begin{figure}
\begin{center}
\includegraphics[scale=0.40]{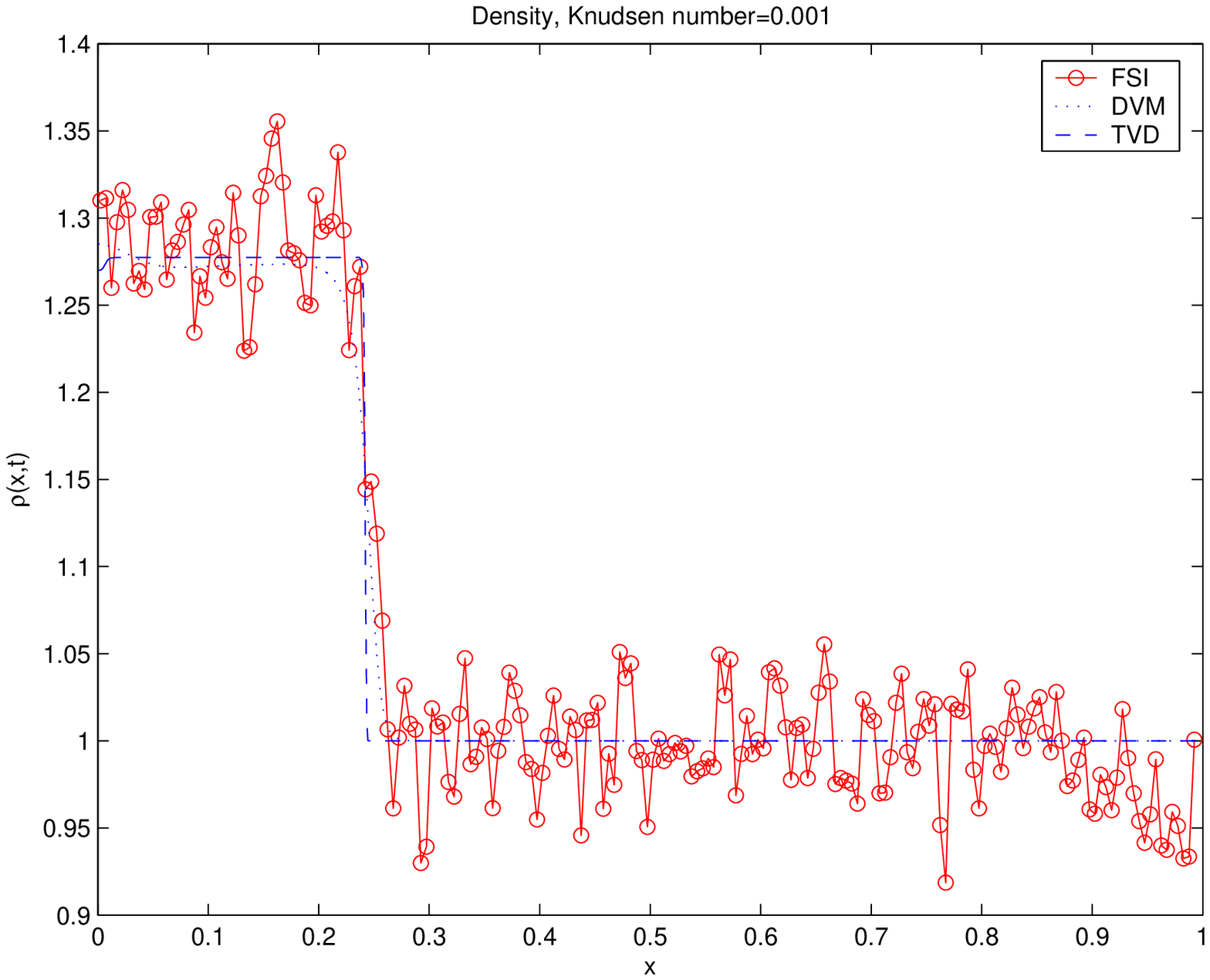}
\includegraphics[scale=0.40]{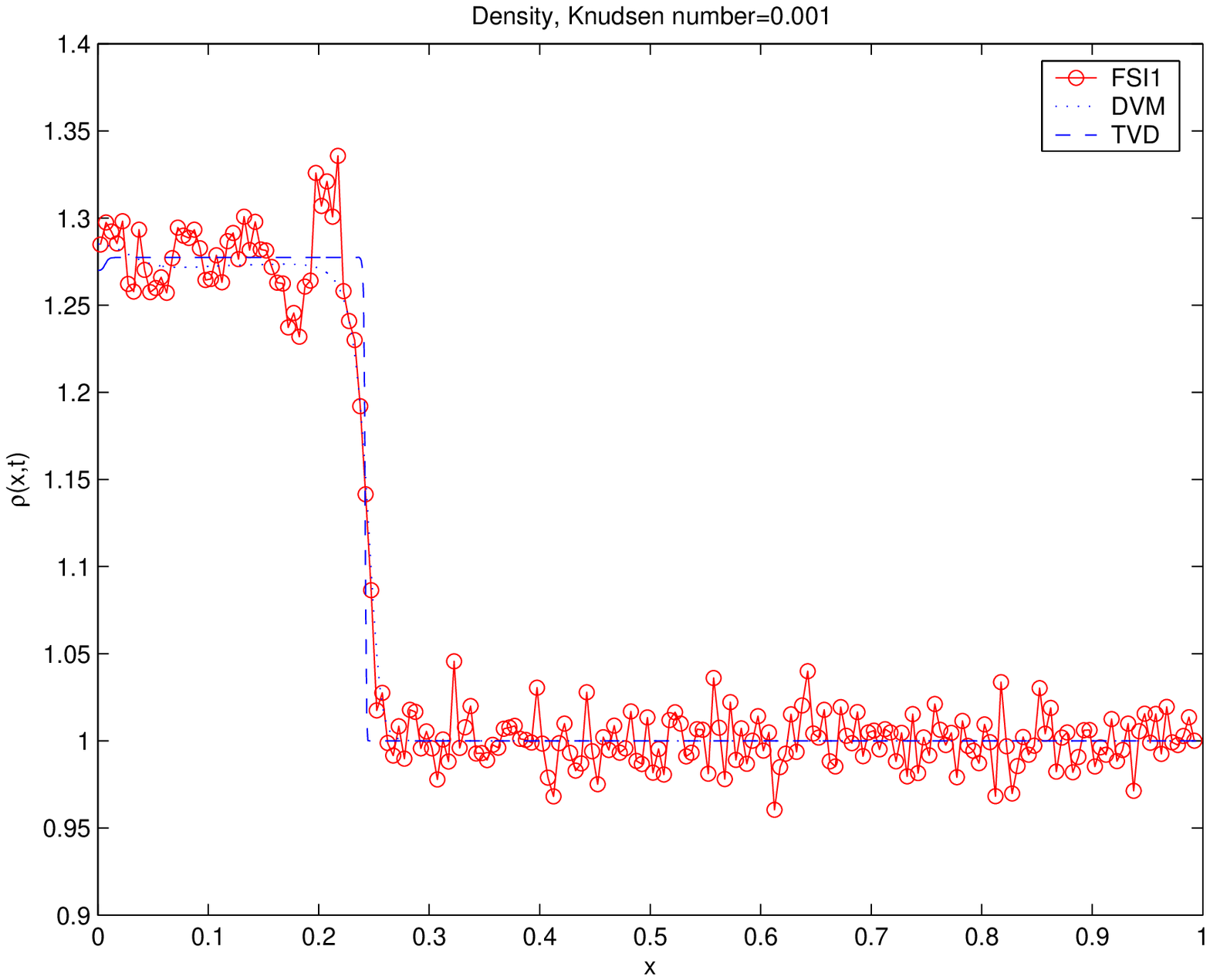}
\includegraphics[scale=0.40]{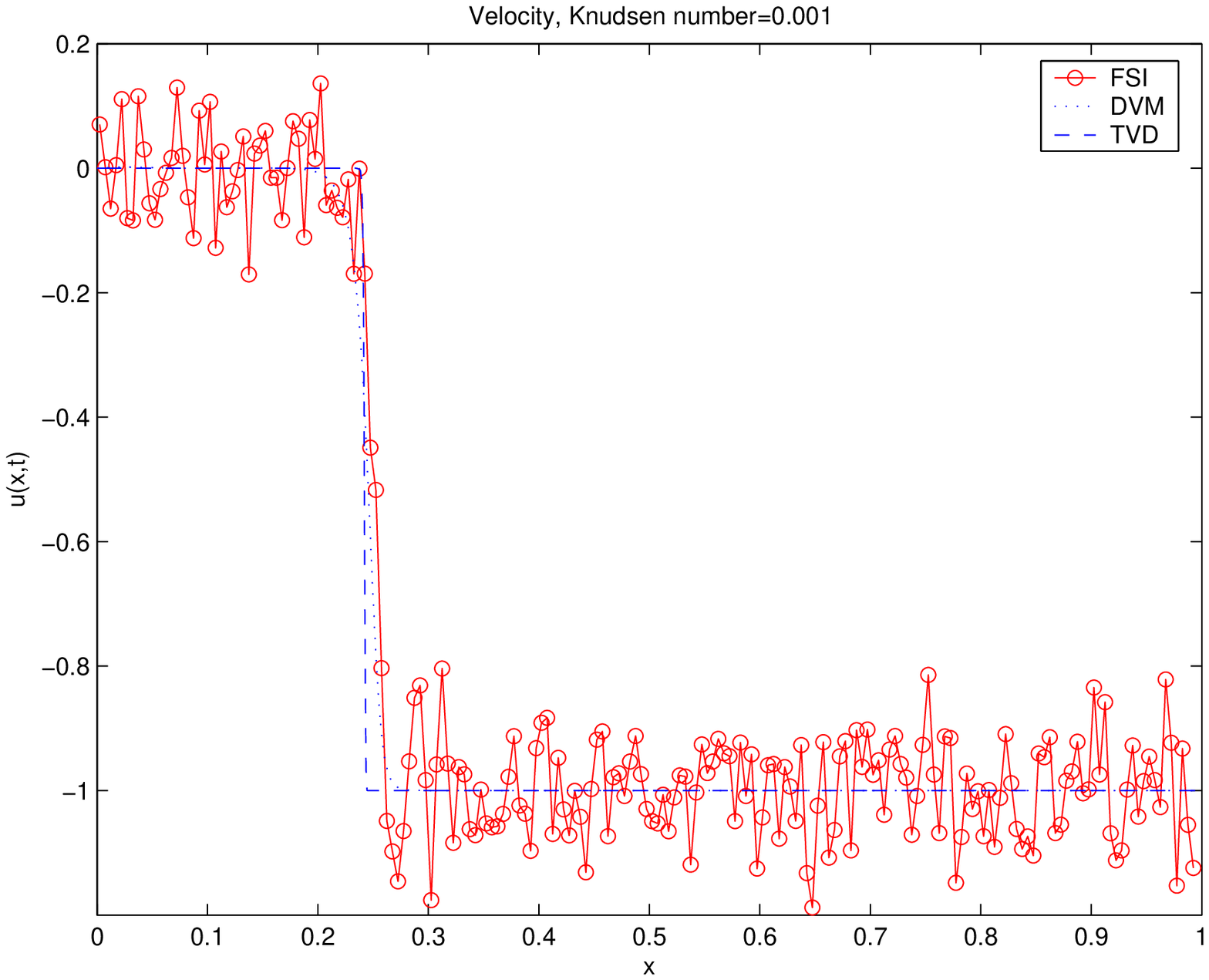}
\includegraphics[scale=0.40]{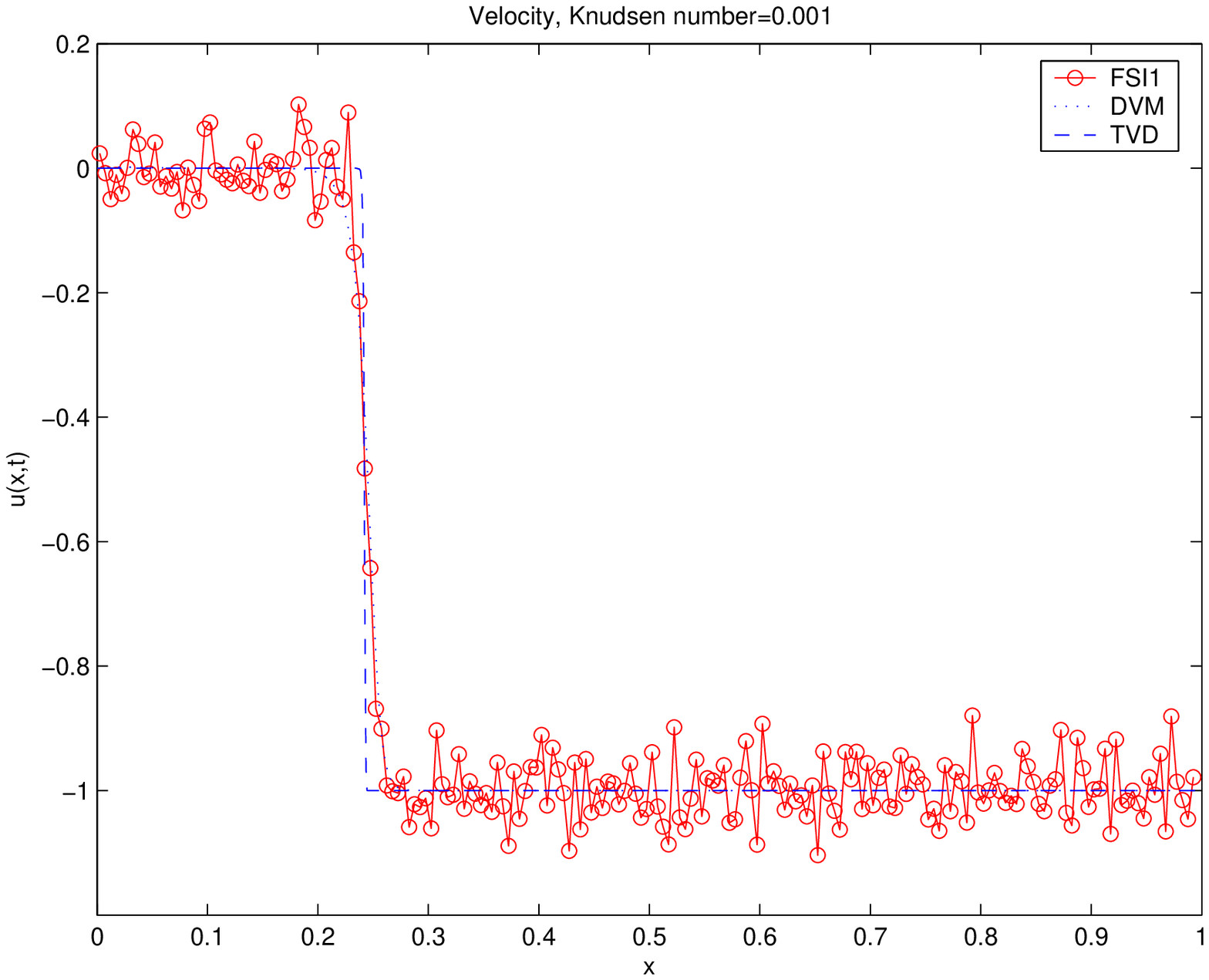}
\includegraphics[scale=0.40]{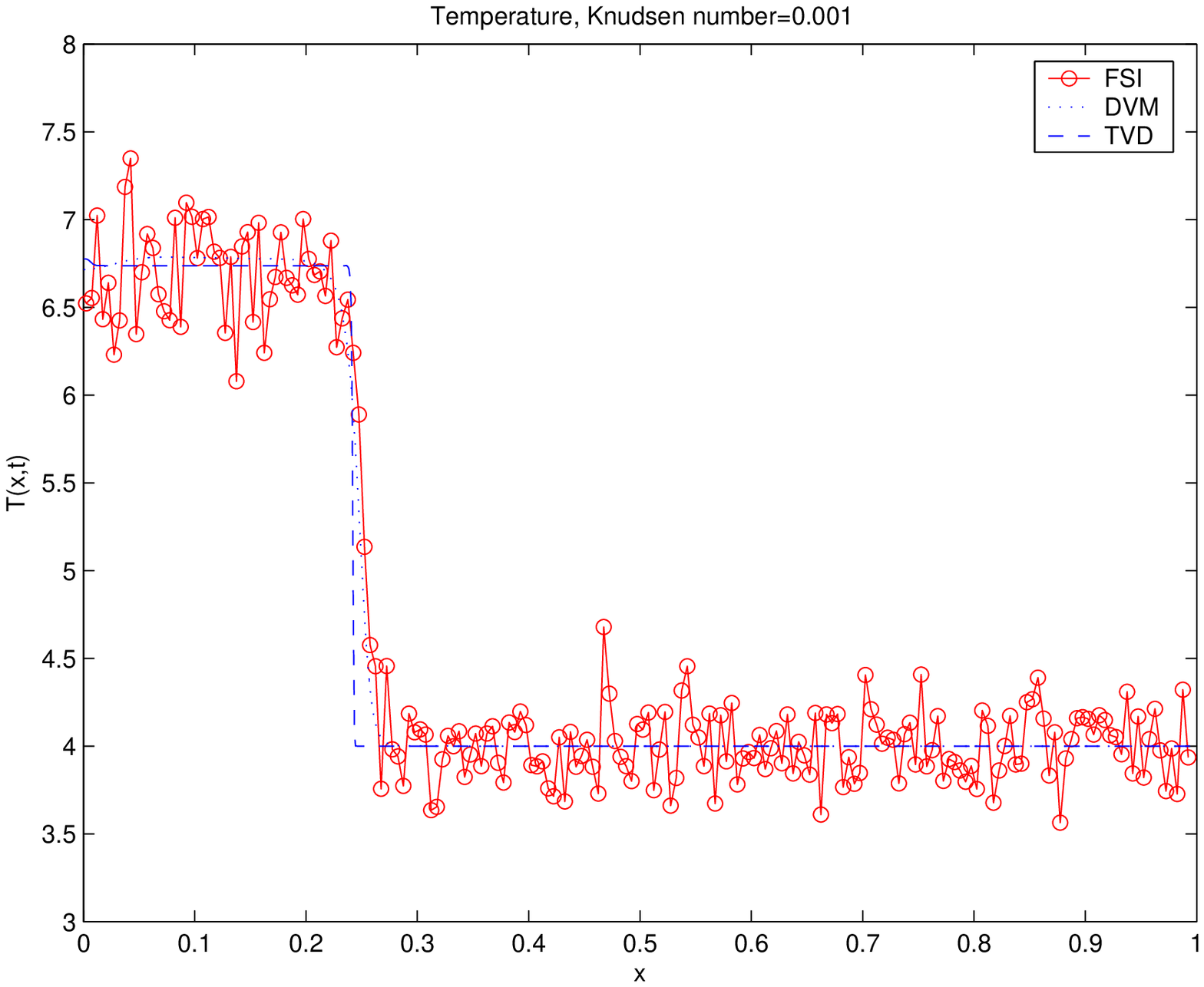}
\includegraphics[scale=0.40]{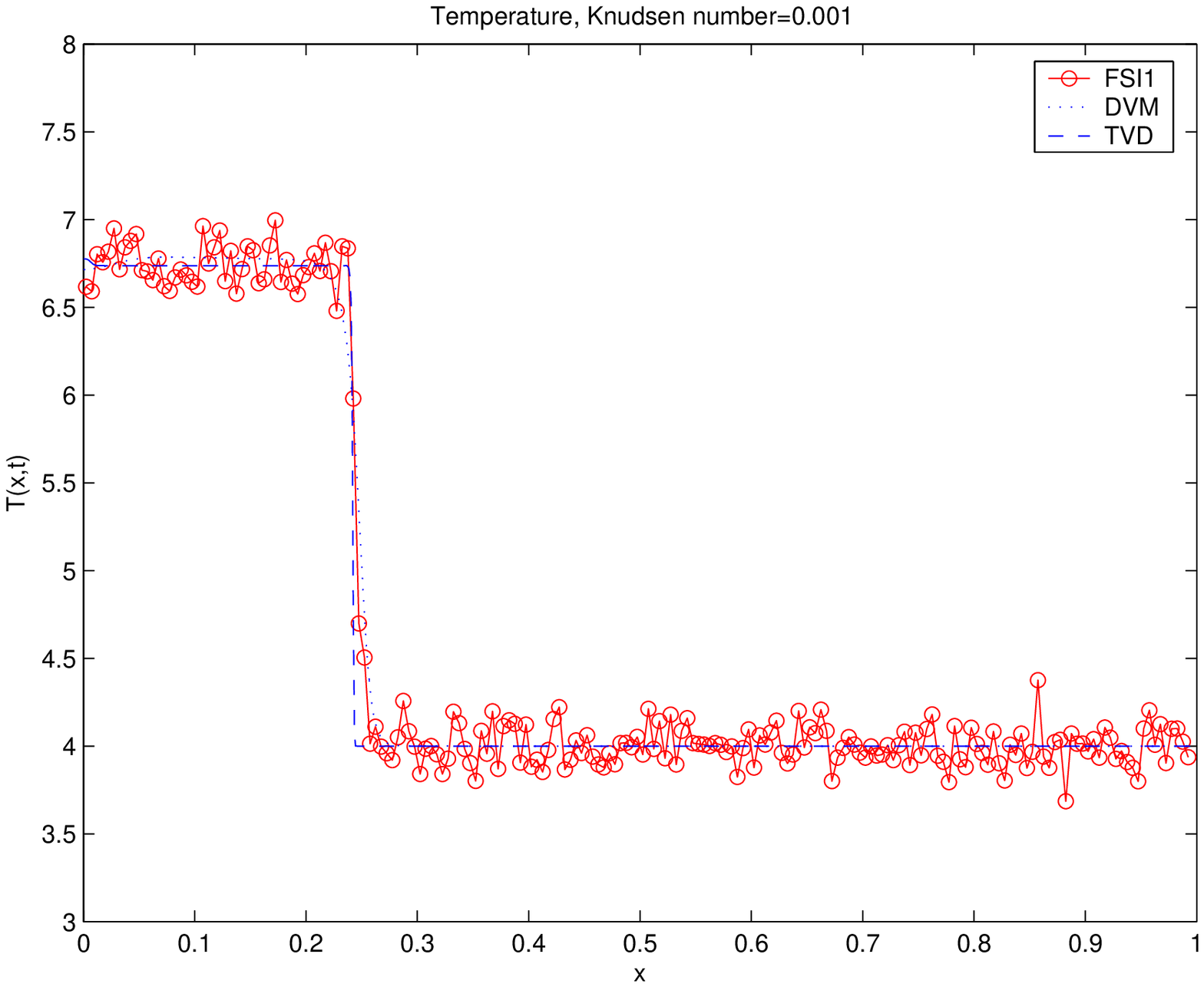}
\caption{Unsteady Shock: $\varepsilon=10^{-3}$. Solution at
$t=0.065$ for FSI (left) FSI1 (right). From top to bottom density,
mean velocity and temperature.} \label{US2}
\end{center}
\end{figure}

\begin{figure}
\begin{center}
\includegraphics[scale=0.40]{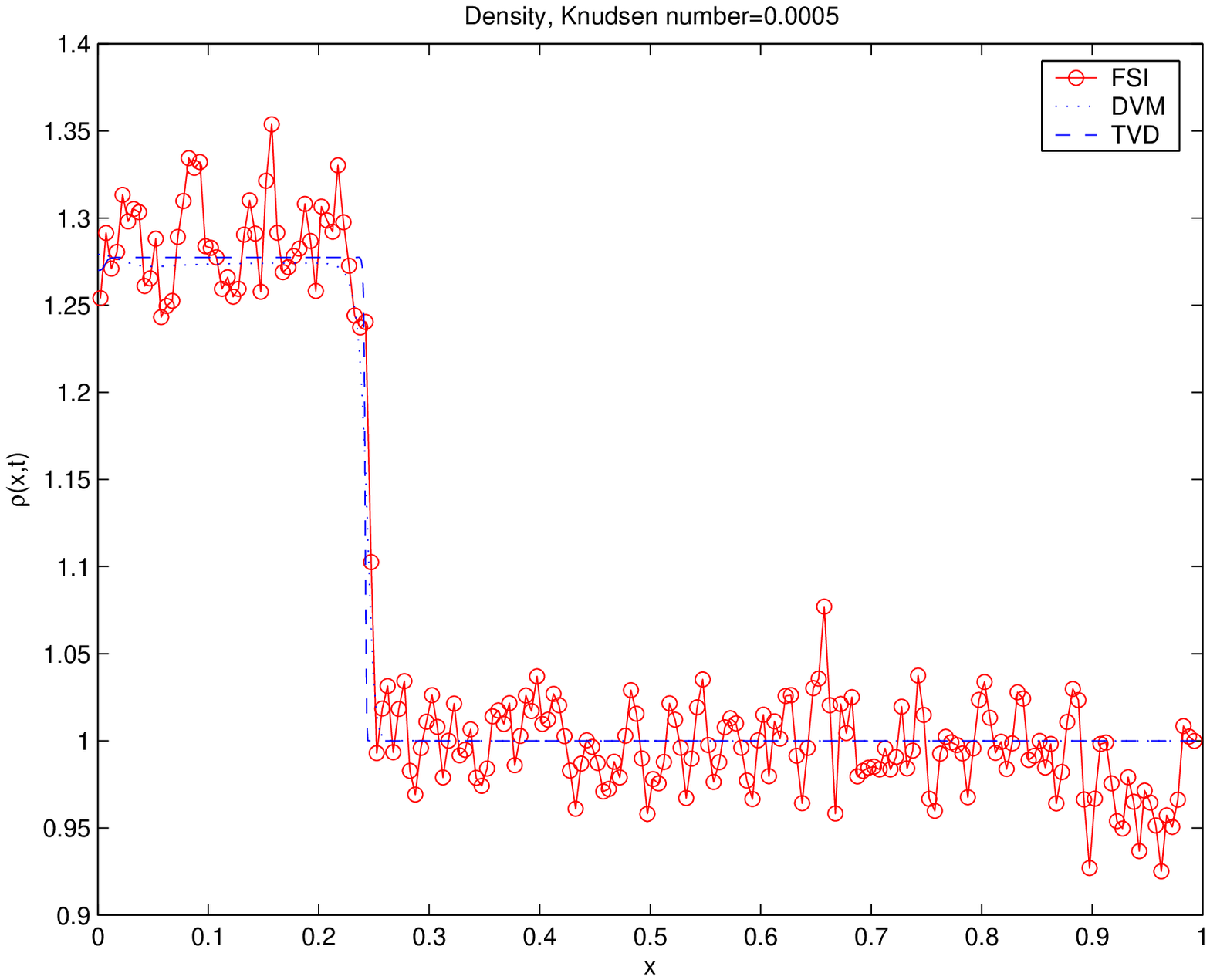}
\includegraphics[scale=0.40]{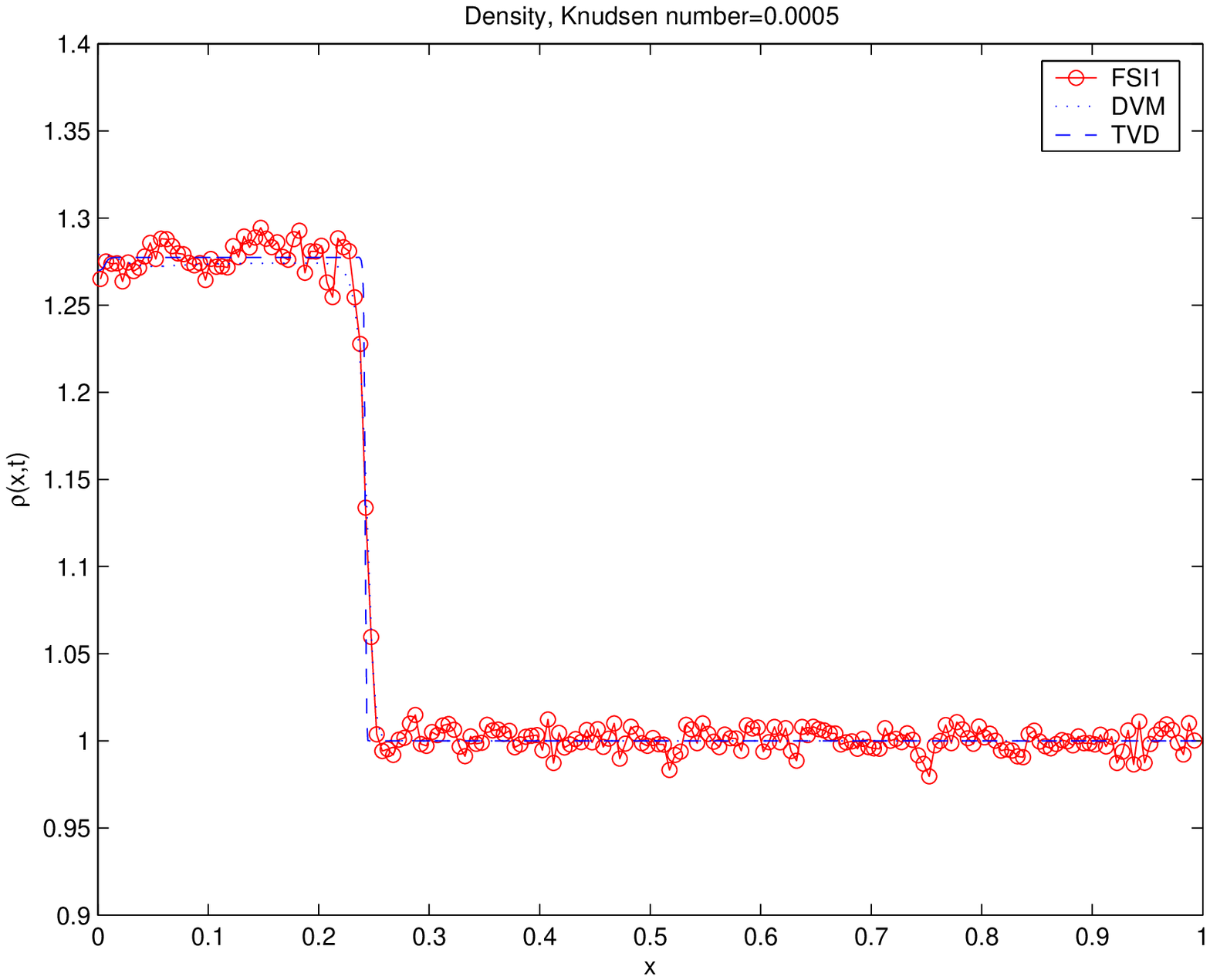}
\includegraphics[scale=0.40]{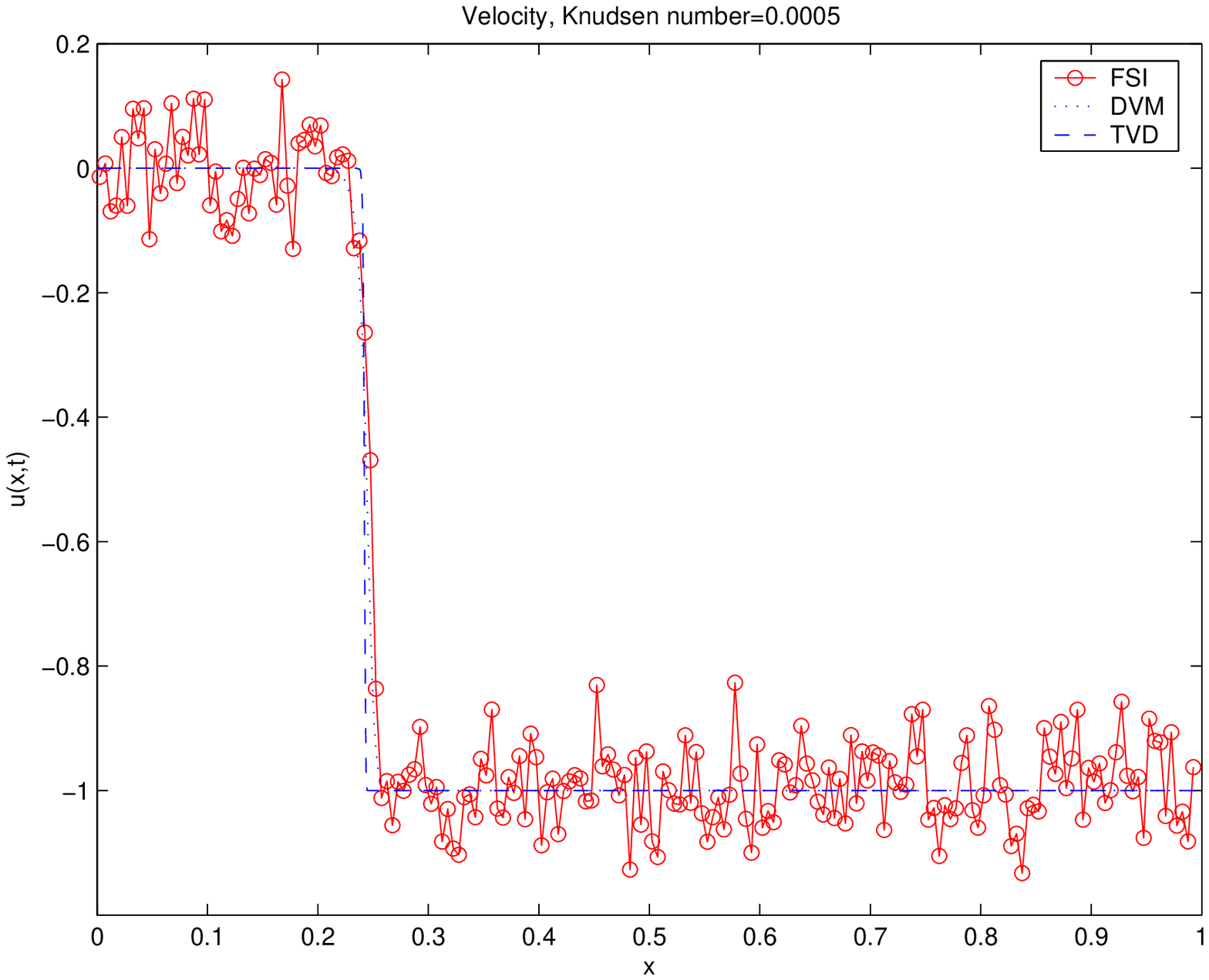}
\includegraphics[scale=0.40]{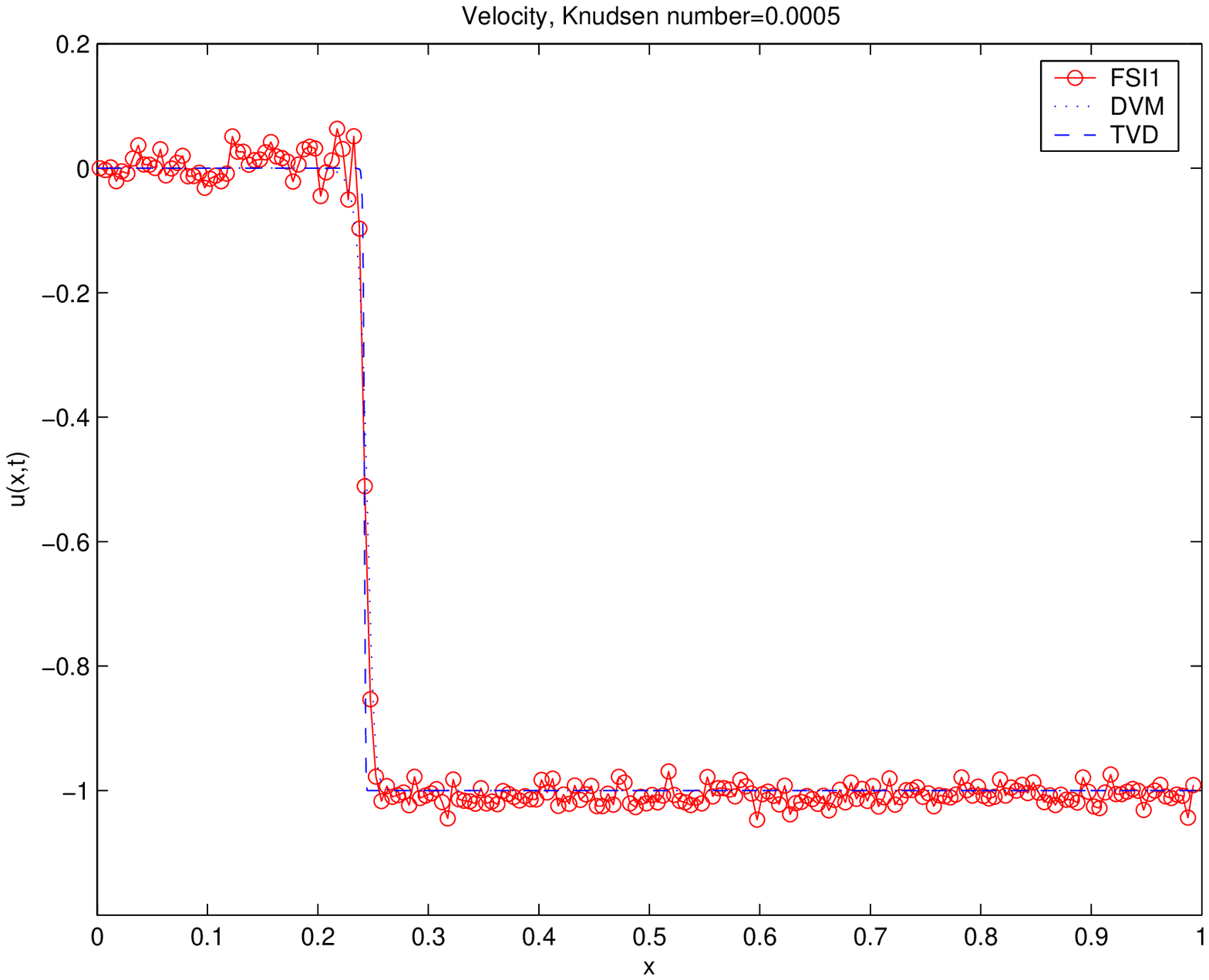}
\includegraphics[scale=0.40]{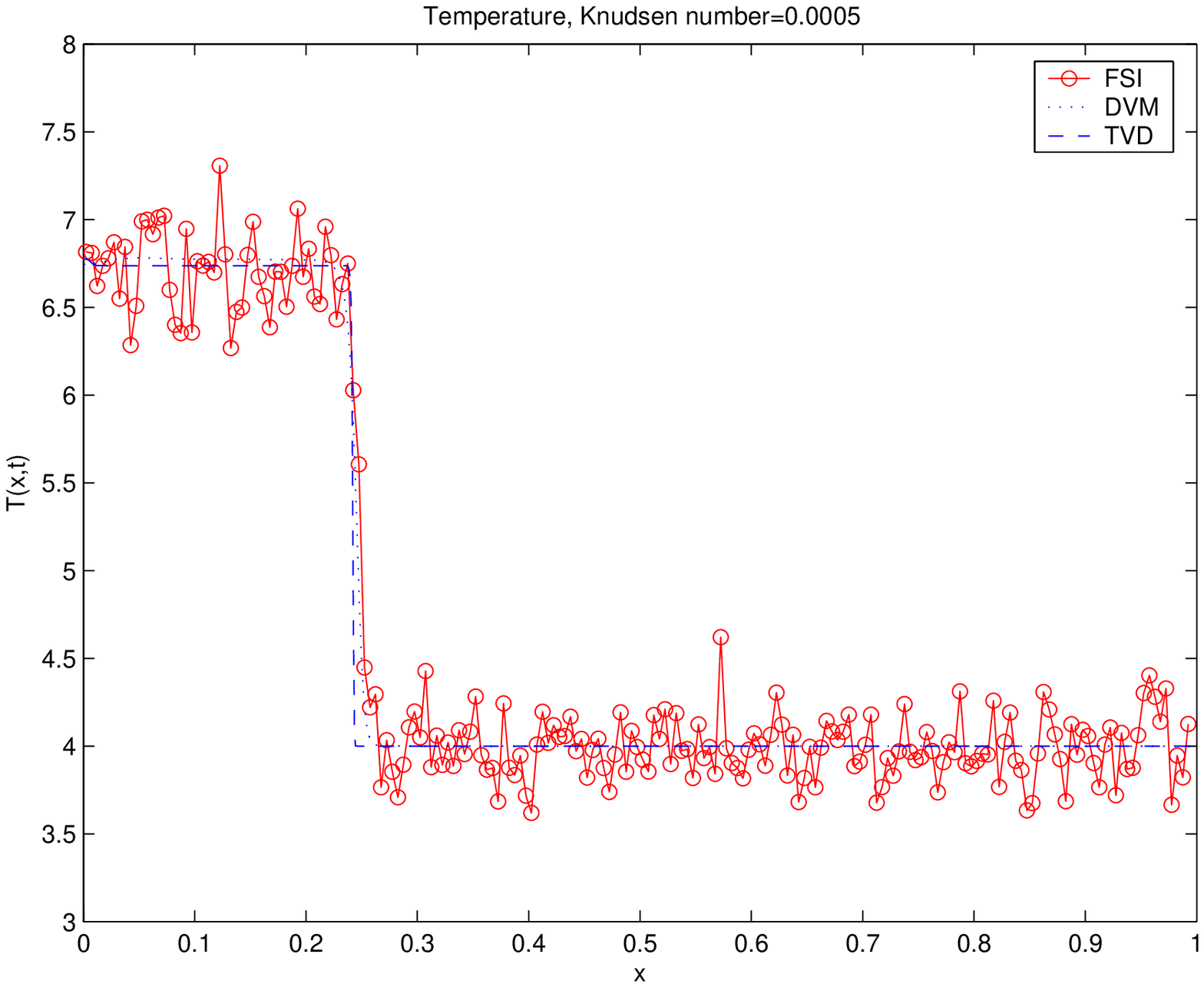}
\includegraphics[scale=0.40]{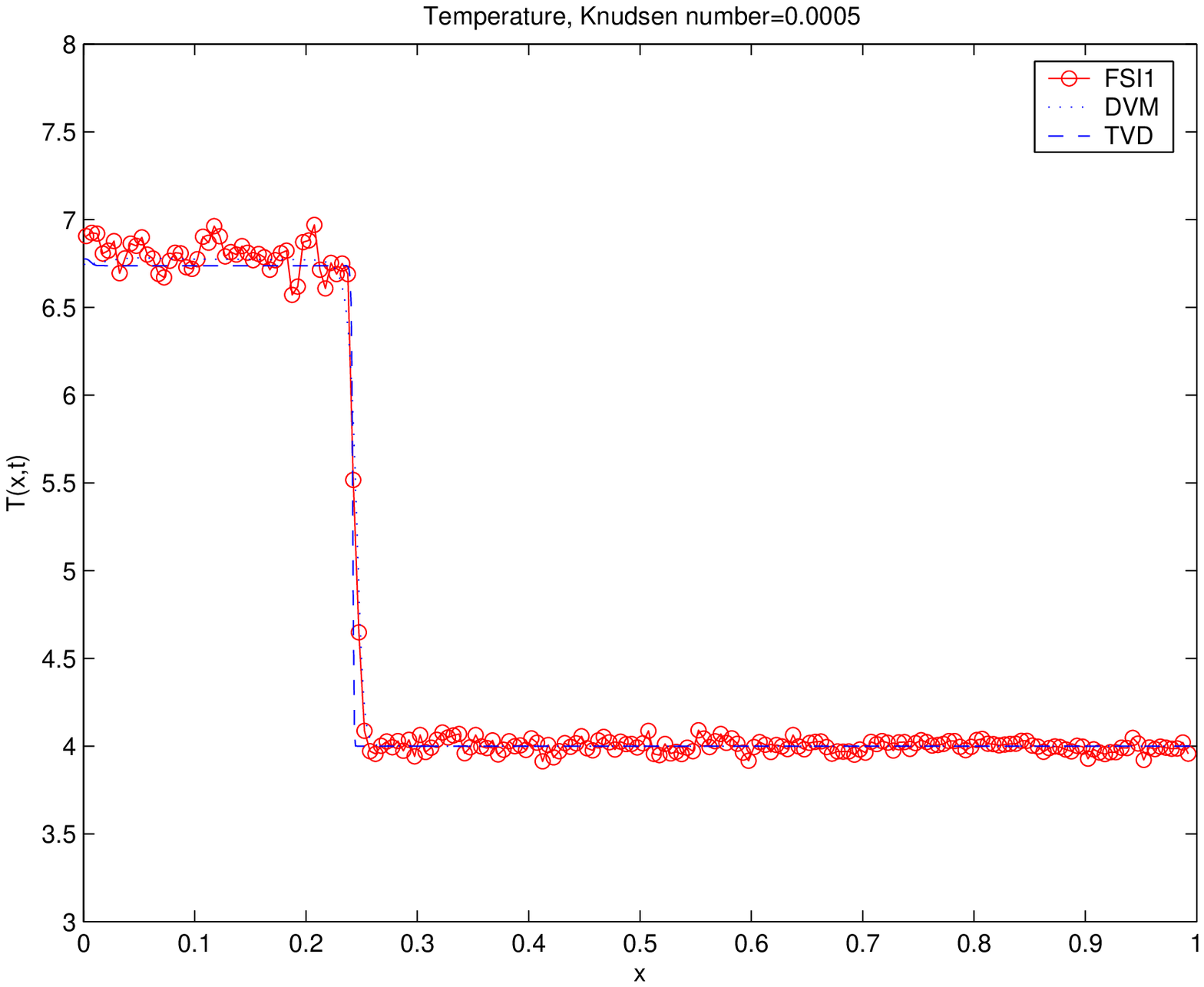}
\caption{Unsteady Shock: $\varepsilon=5\times 10^{-4}$. Solution
at $t=0.065$ for FSI (left) FSI1 (right). From top to bottom
density, mean velocity and temperature.} \label{US3}
\end{center}
\end{figure}

\begin{figure}
\begin{center}
\includegraphics[scale=0.40]{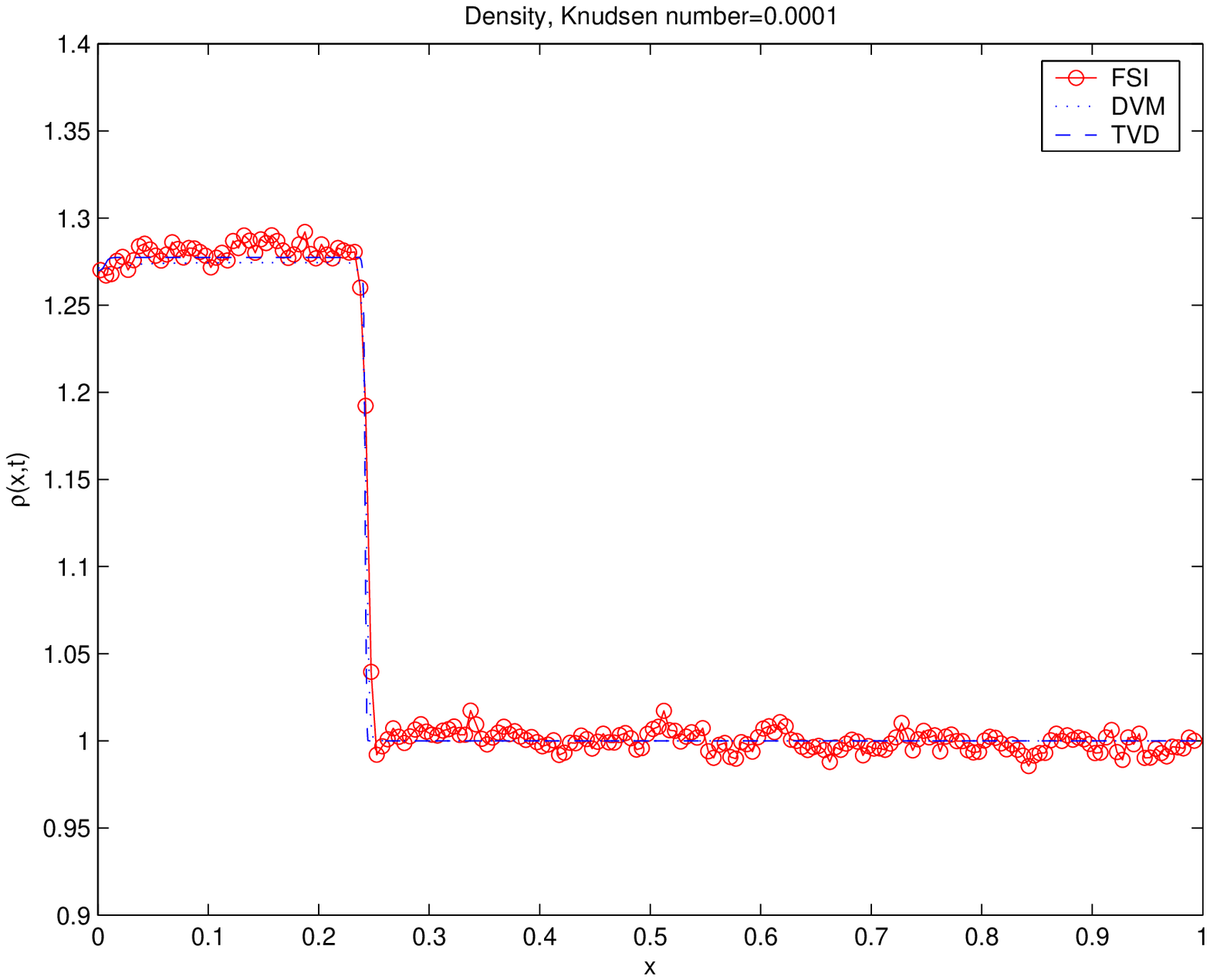}
\includegraphics[scale=0.40]{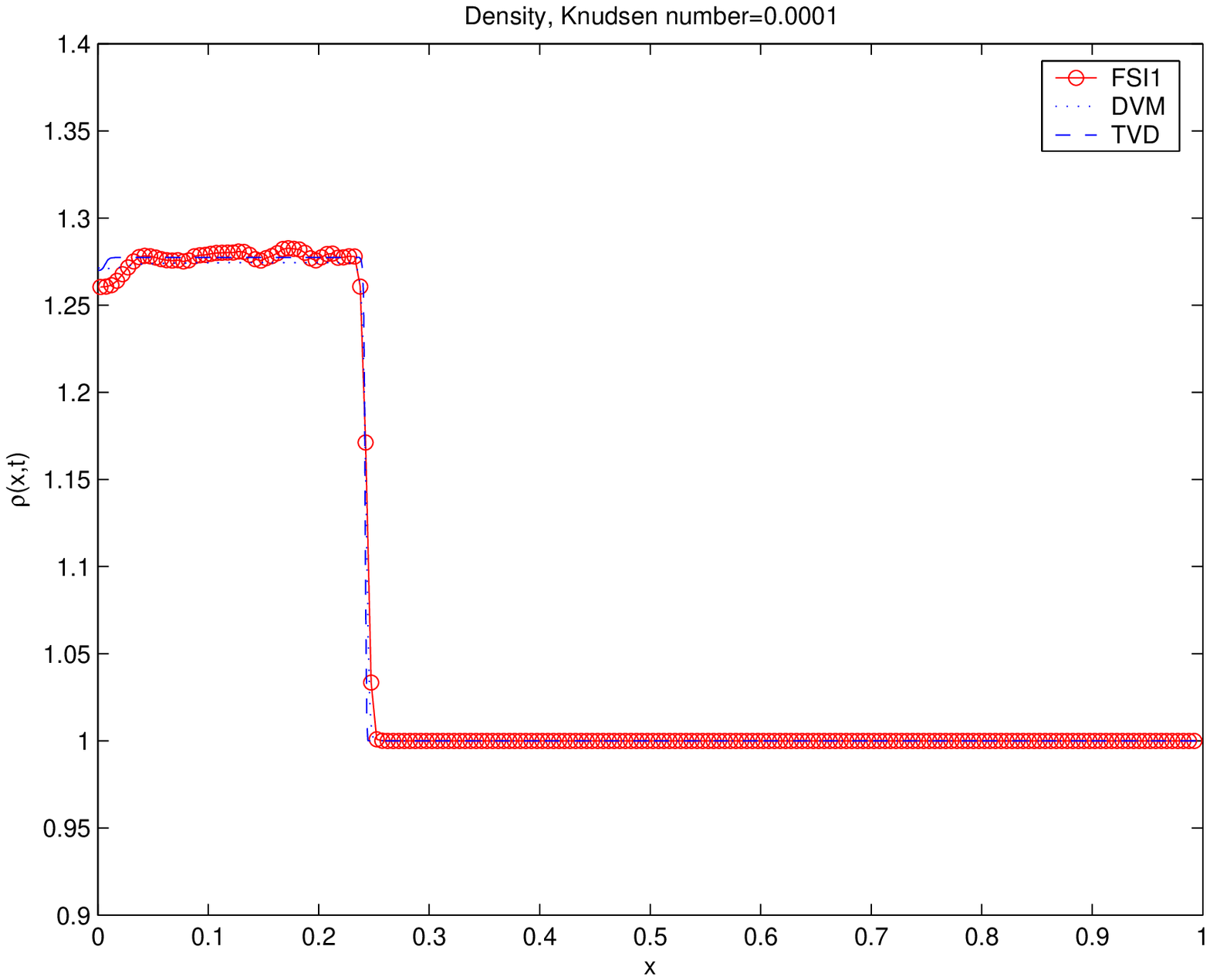}
\includegraphics[scale=0.40]{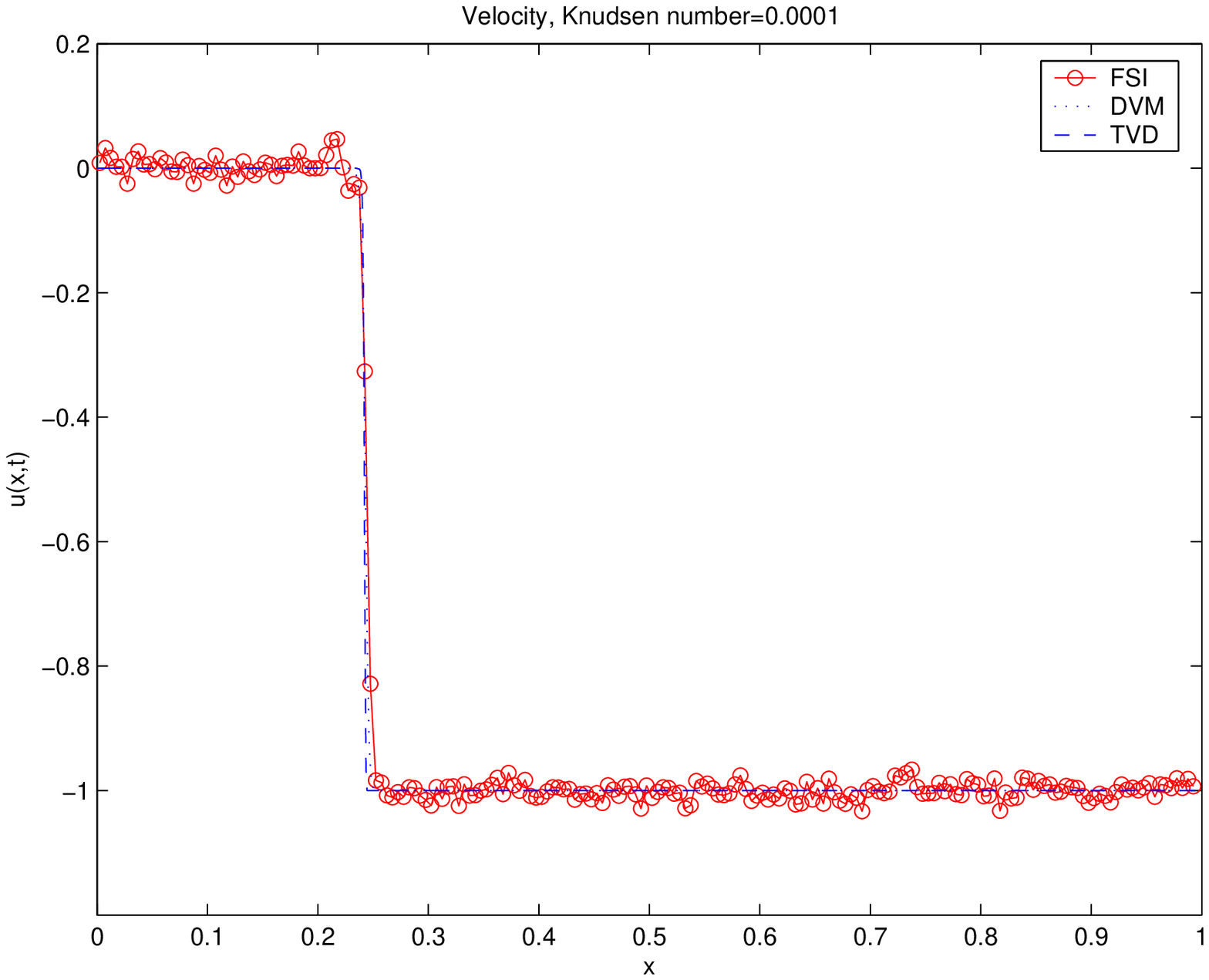}
\includegraphics[scale=0.40]{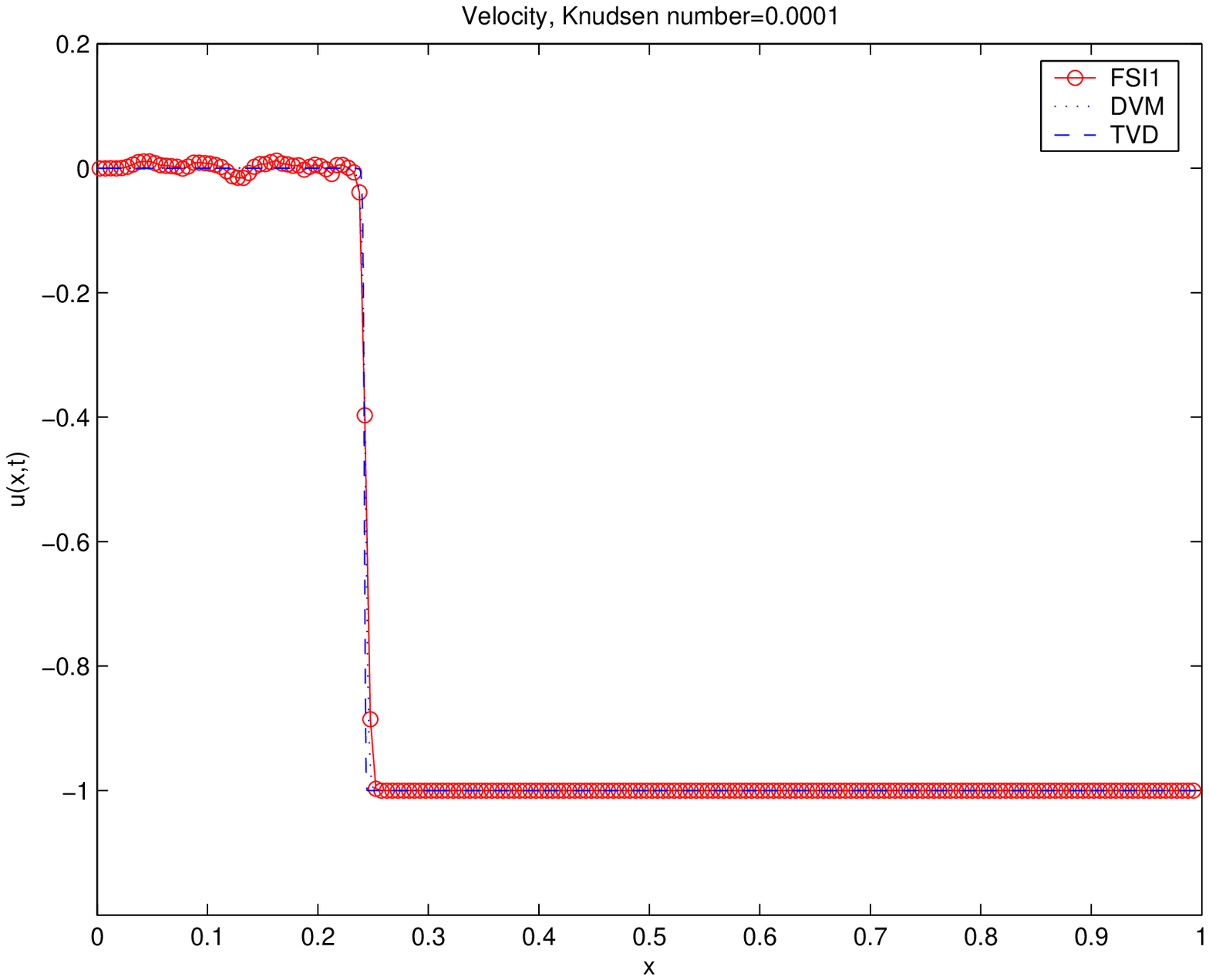}
\includegraphics[scale=0.40]{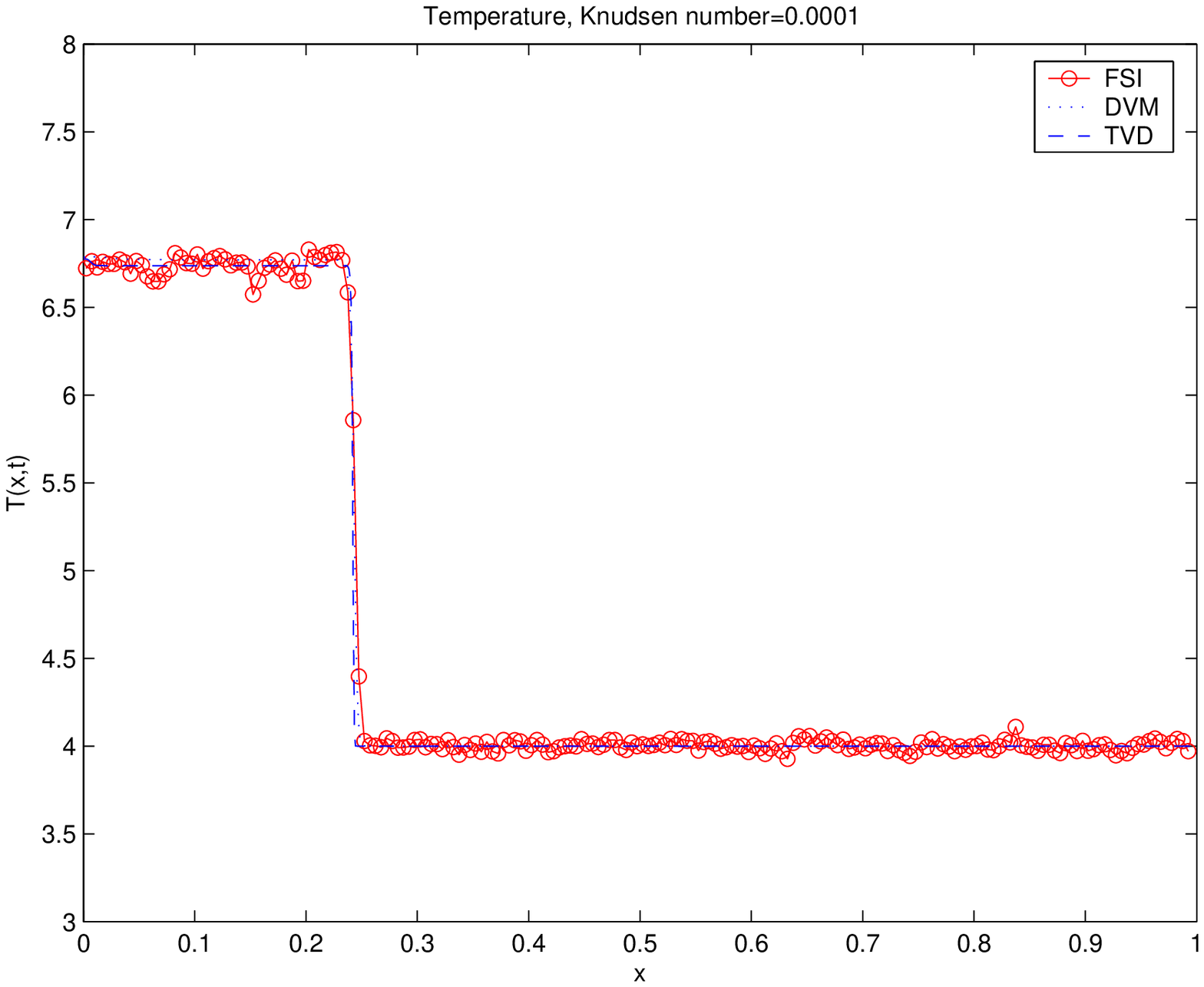}
\includegraphics[scale=0.40]{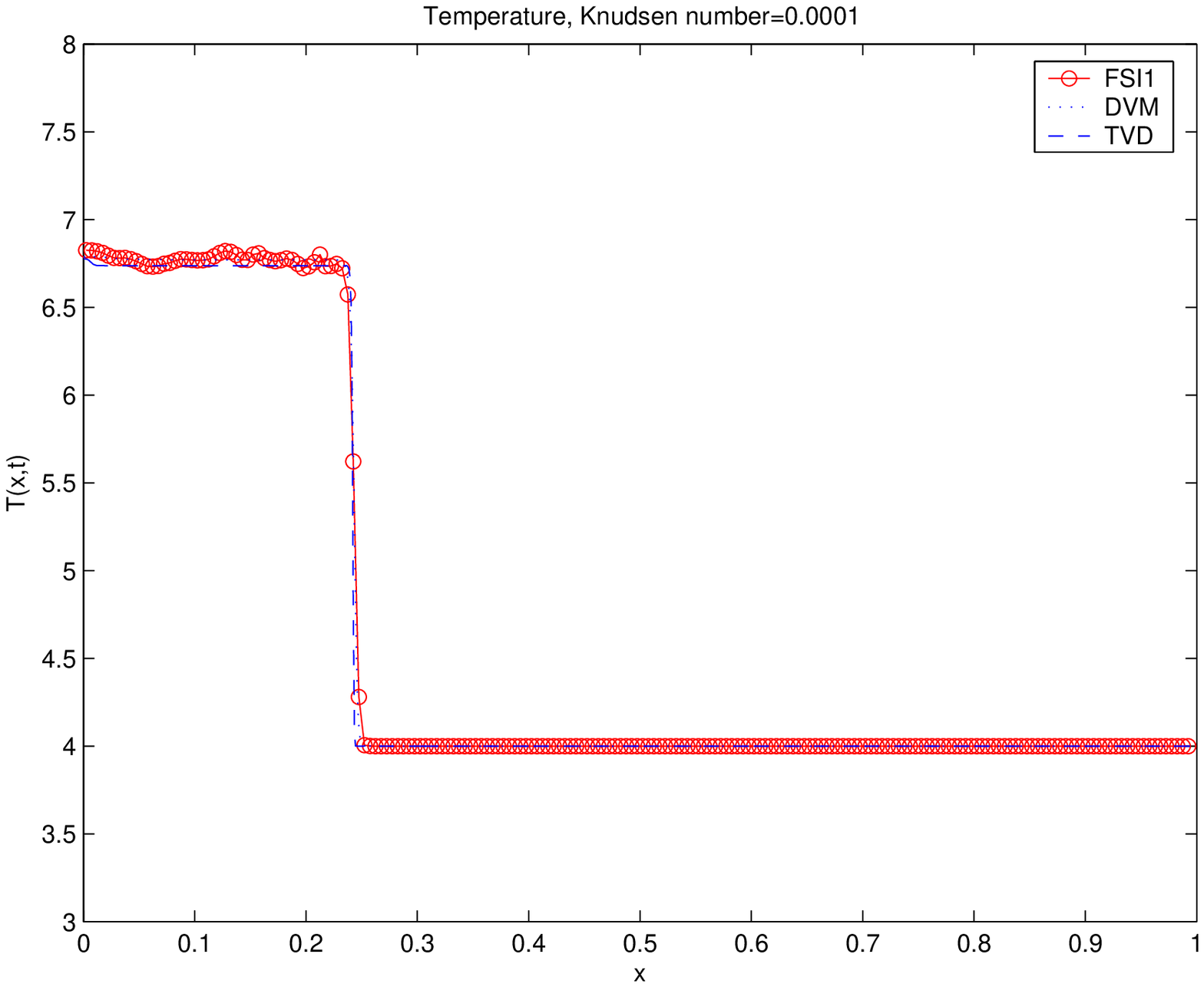}
\caption{Unsteady Shock: $\varepsilon=10^{-4}$. Solution at
$t=0.065$ for FSI (left) FSI1 (right). From top to bottom density,
mean velocity and temperature.} \label{US4}
\end{center}
\end{figure}

\begin{figure}
\begin{center}
\includegraphics[scale=0.40]{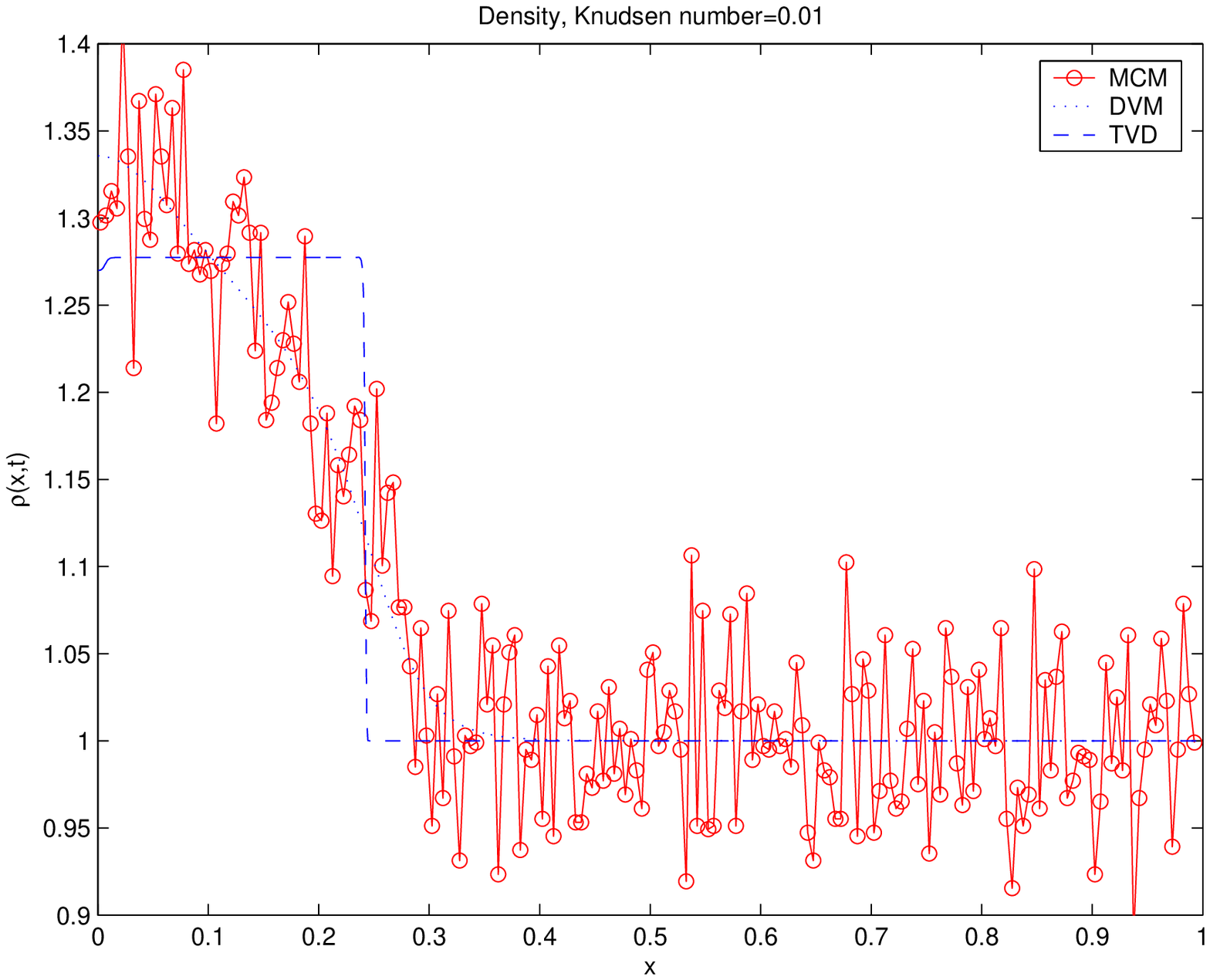}
\includegraphics[scale=0.40]{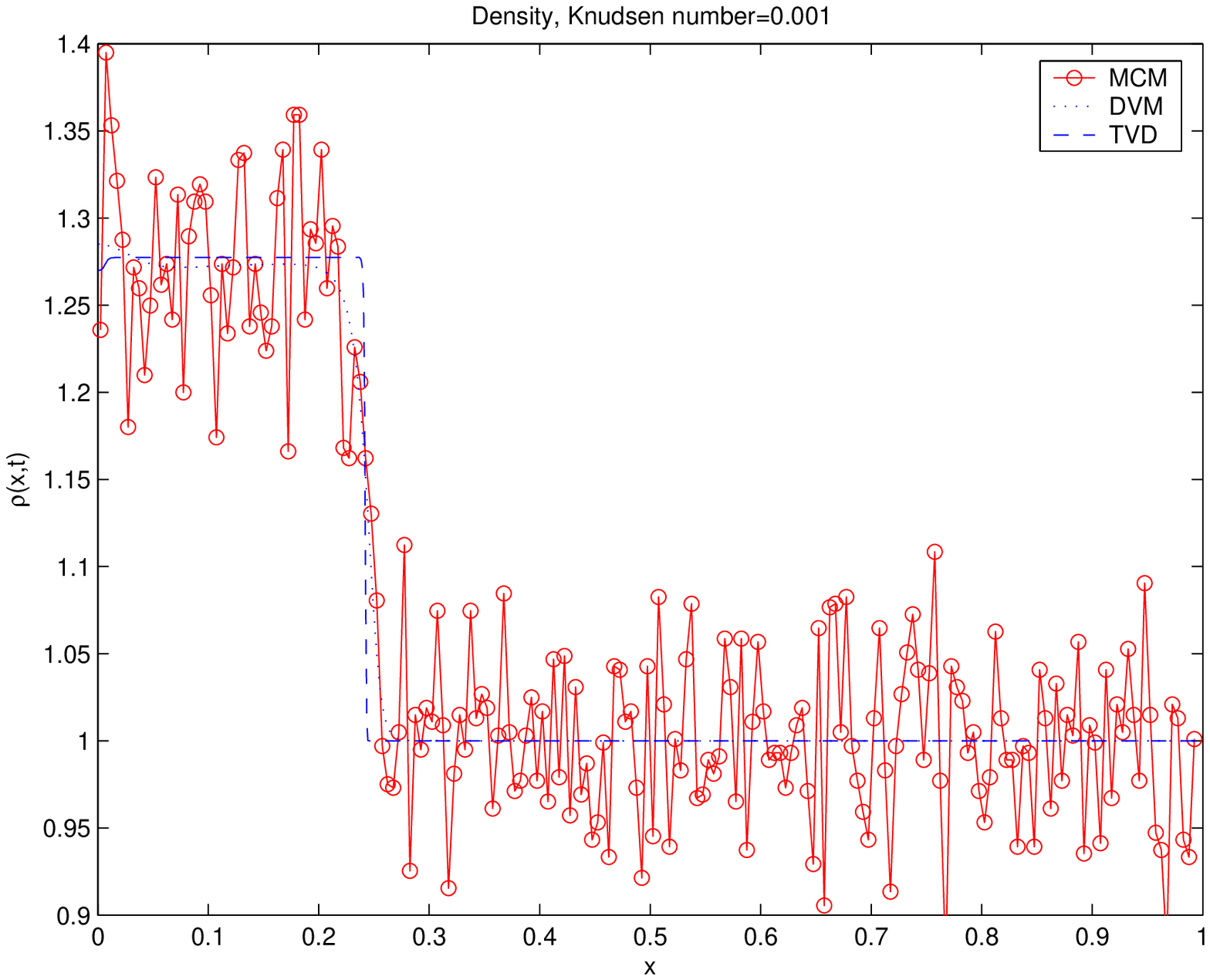}
\includegraphics[scale=0.40]{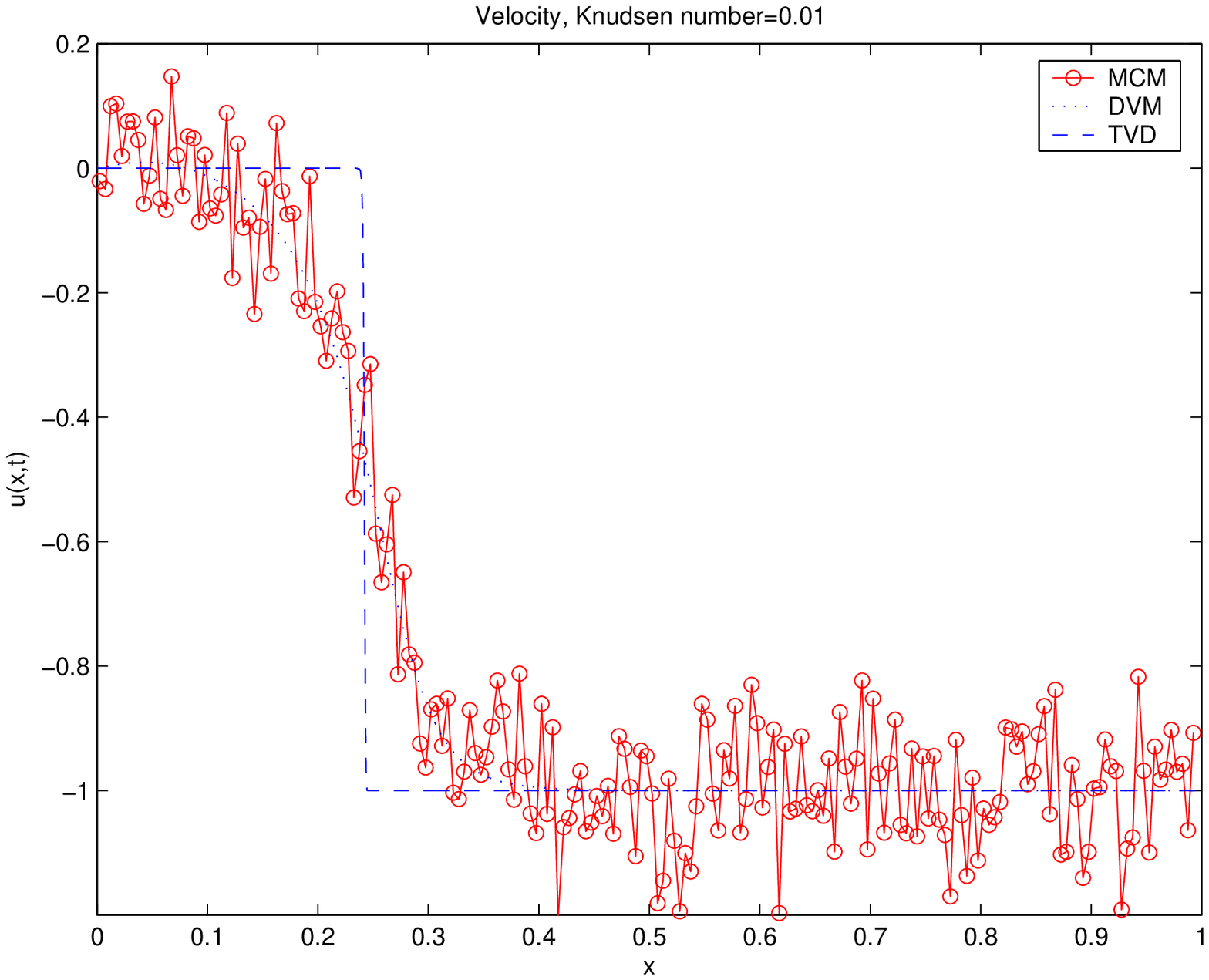}
\includegraphics[scale=0.40]{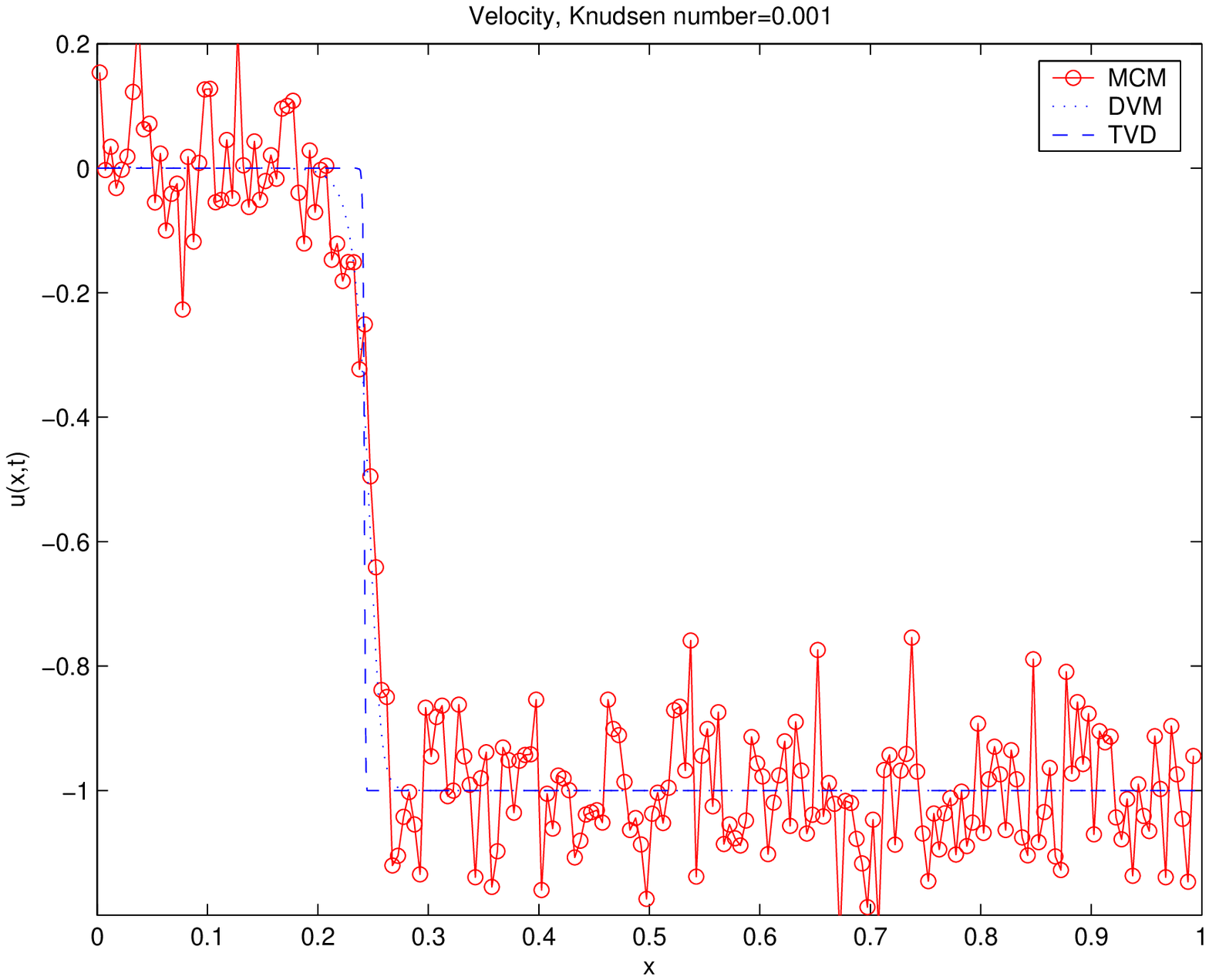}
\includegraphics[scale=0.40]{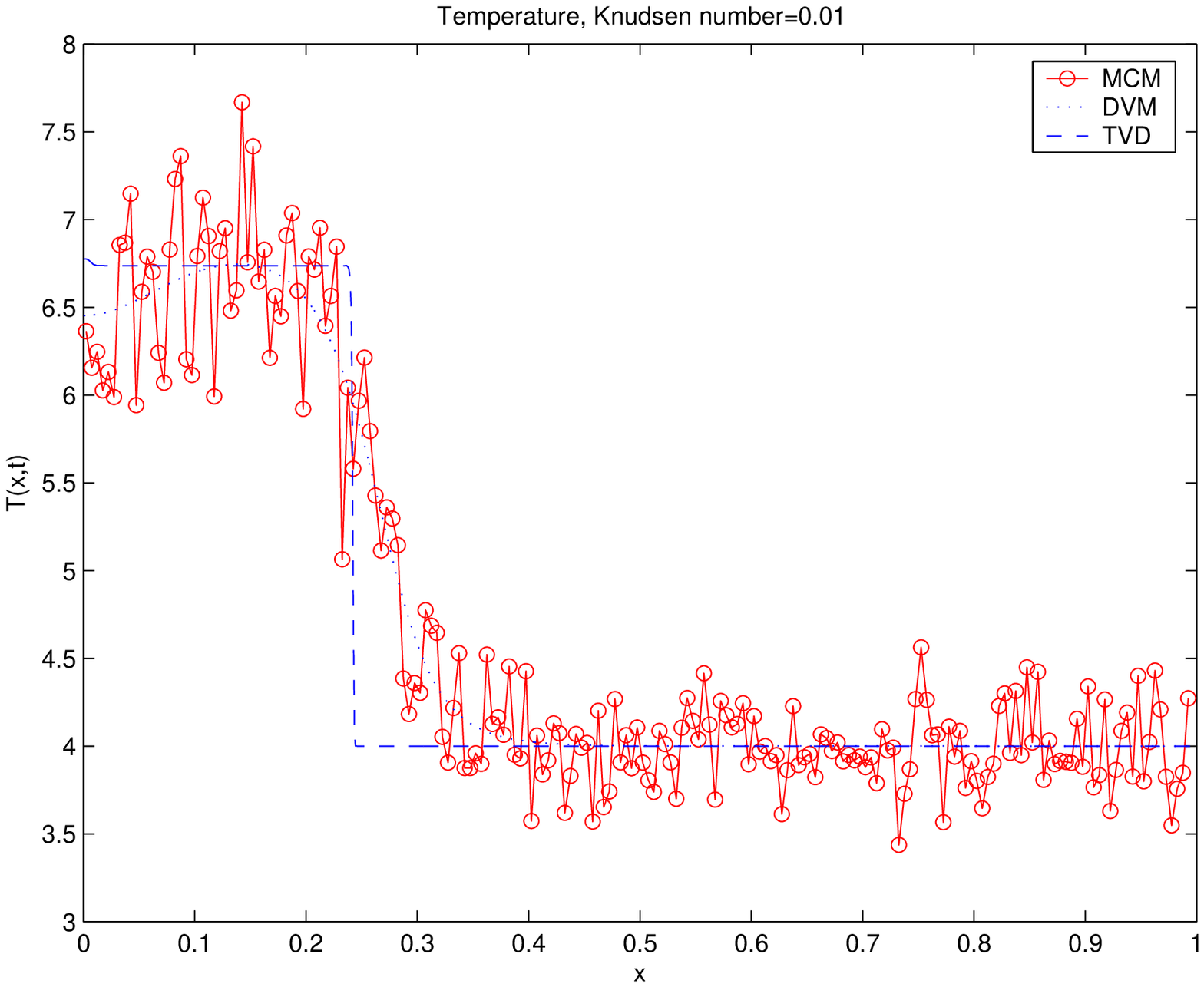}
\includegraphics[scale=0.40]{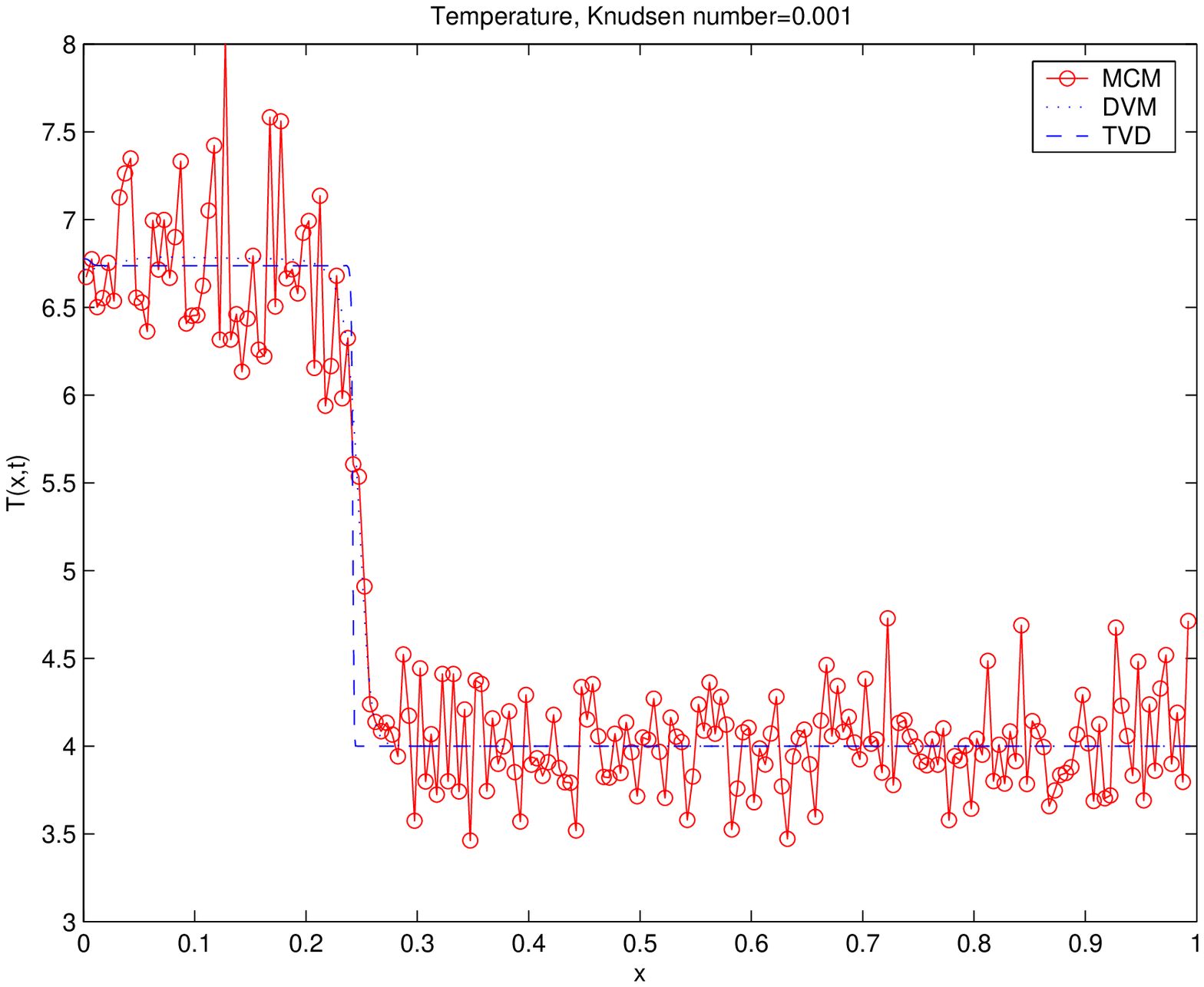}
\caption{Unsteady Shock. Solution at $t=0.065$ for MCM with
Knudsen numbers $\varepsilon=10^{-2}$ (left) and
$\varepsilon=10^{-3}$ (right). From top to bottom density, mean
velocity and temperature.} \label{US5}
\end{center}
\end{figure}

Next we consider an unsteady shock that propagates from left to
right. The shock is produced miming a specular wall on the left
boundary, thus for the stochastic component at each time step,
particles which escape from the computational domain on the left
side are put back in the first cell with opposite velocity and
opposite position respect to zero. On the other side particles are
injected with the initial mean velocity and temperature in a
number which corresponds to the initial density. For the
macroscopic part the usual specular and inflow boundary condition
are used. At the beginning the flow is uniform with mass
$\varrho=1$, mean velocity $u=-1$ and energy $E=2.5$. The
computations are stopped when $t=0.065$, the number of cells are
$200$ in space, while the initial number of particle are $500$ for
each space cell. In each Figure the solution computed with the
Euler scheme and the one computed with the DVM scheme is reported.
The FSI, FSI1 and MCM are respectively depicted for density, mean
velocity and temperature. We observe that for large Knudsen
numbers FSI (Figure \ref{US1} left) and FSI1 (Figure \ref{US1}
right) provide a small improvement with respect to MCM (Figure
\ref{US5} left) in term of fluctuations. When $\varepsilon$
decreases the non-equilibrium part becomes smaller and both FSI
and FSI1 (Figure \ref{US2} and \ref{US3}) contain less
fluctuations than MCM. Note that, since the time step here is
$O(\e)$, the reduction of fluctuations in FSI scheme is
essentially the same for $\e=0.001$ and $\e=0.0005$ whereas for
FSI1 scheme the solution shows a remarkable improvement as $\e$
diminishes. Finally for $\varepsilon=10^{-4}$ (Figure \ref{US4})
we are in an under-resolved regime and both hybrid methods yield
similar solutions at the same computational time.


\begin{figure}
\begin{center}
\includegraphics[scale=0.40]{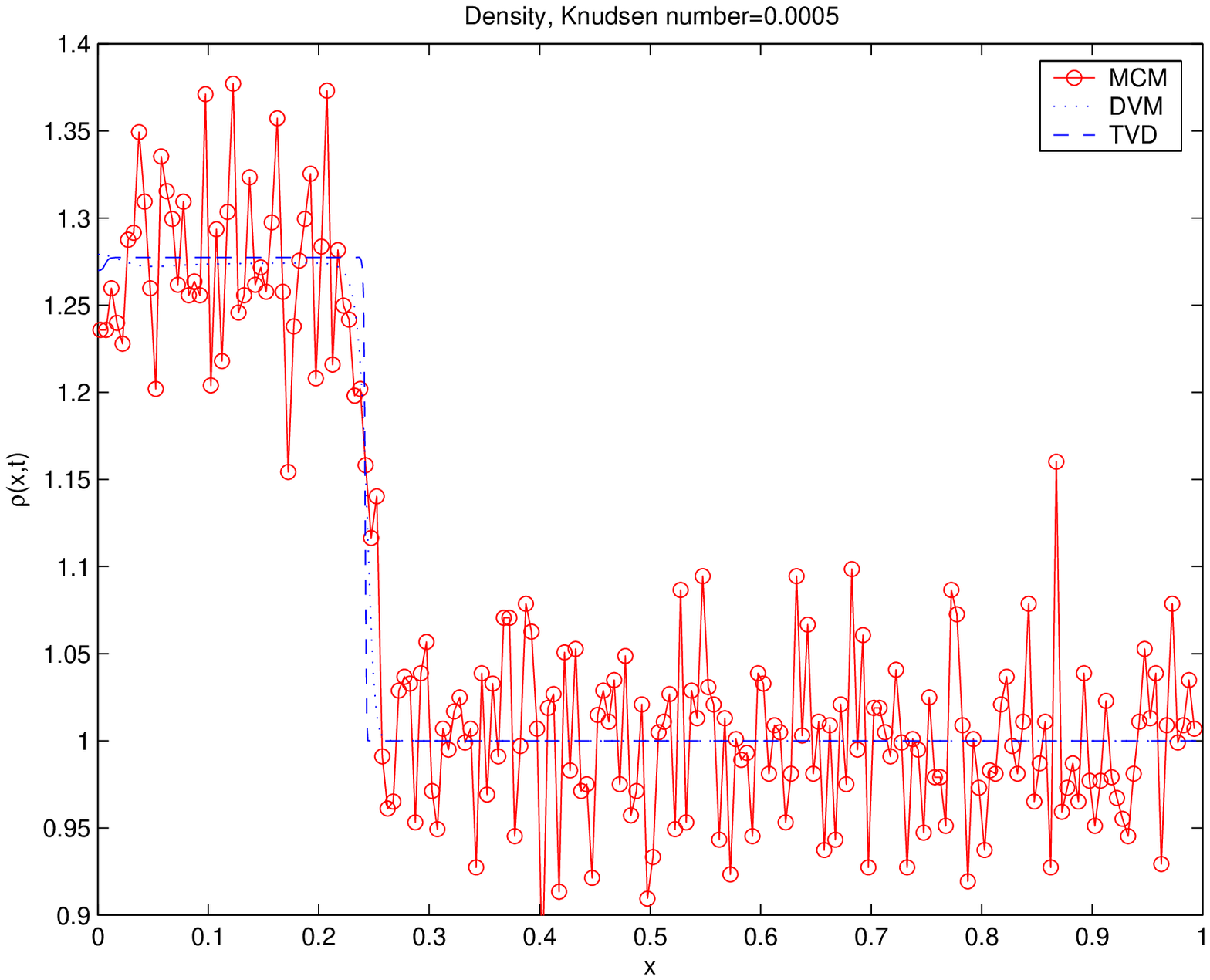}
\includegraphics[scale=0.40]{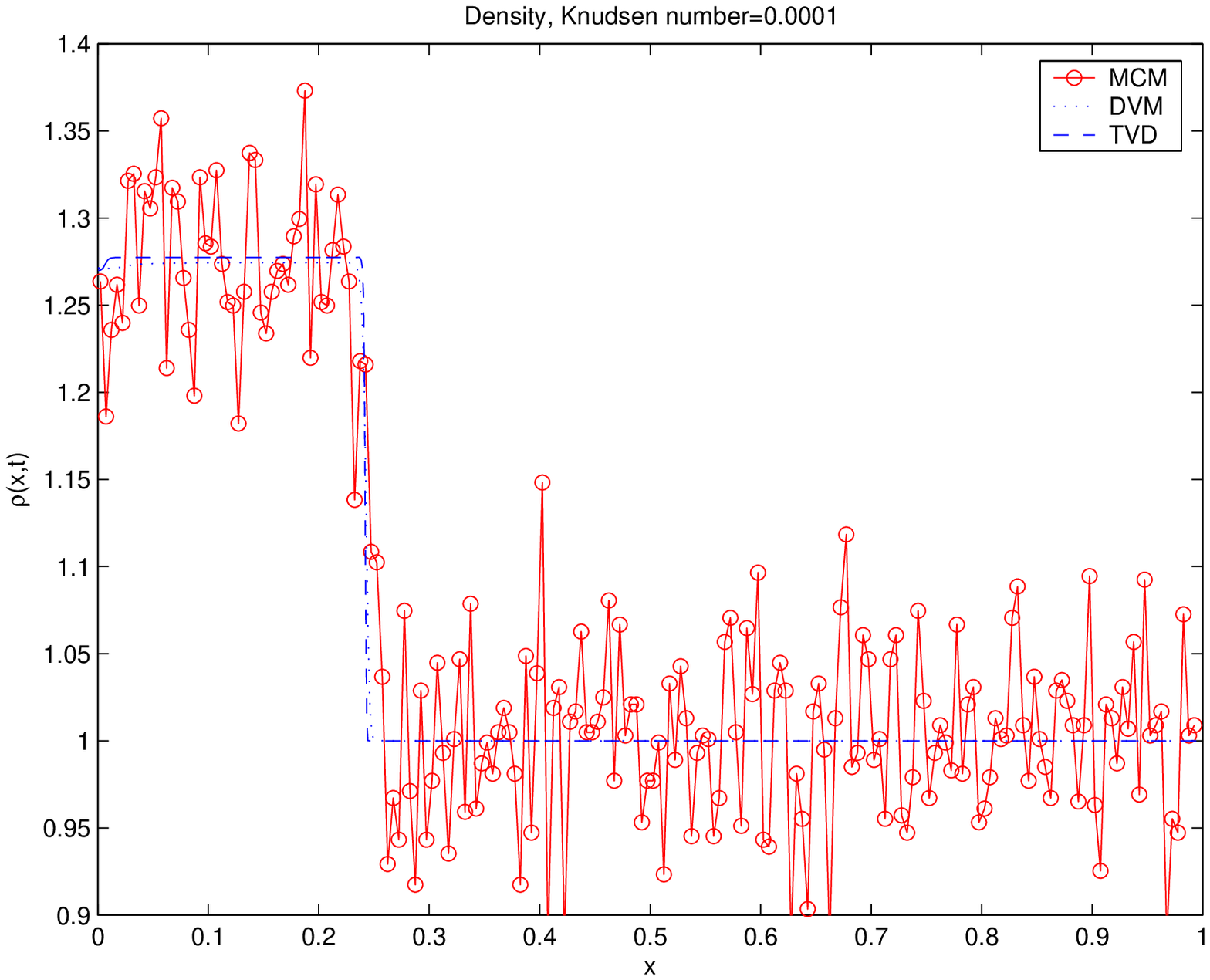}
\includegraphics[scale=0.40]{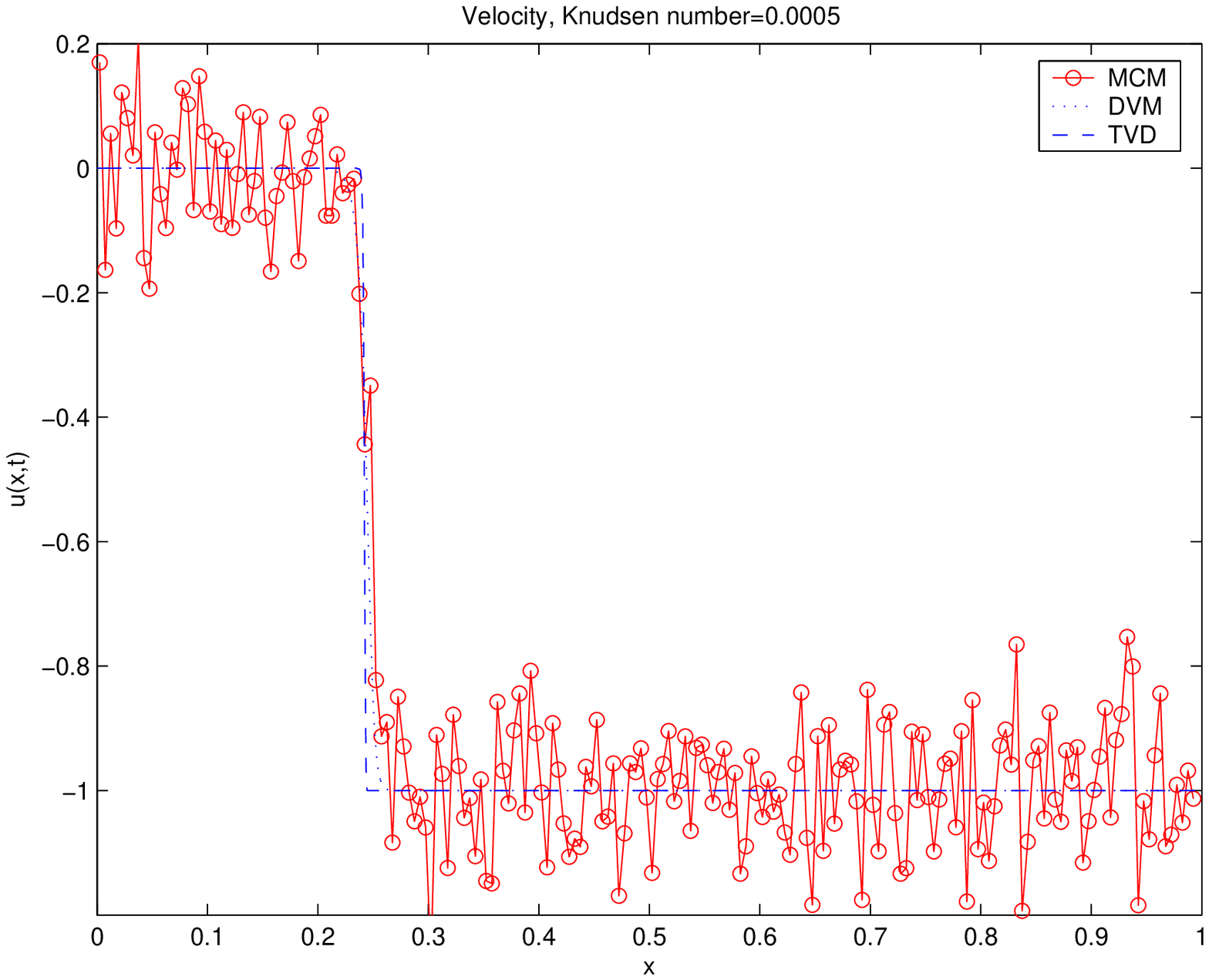}
\includegraphics[scale=0.40]{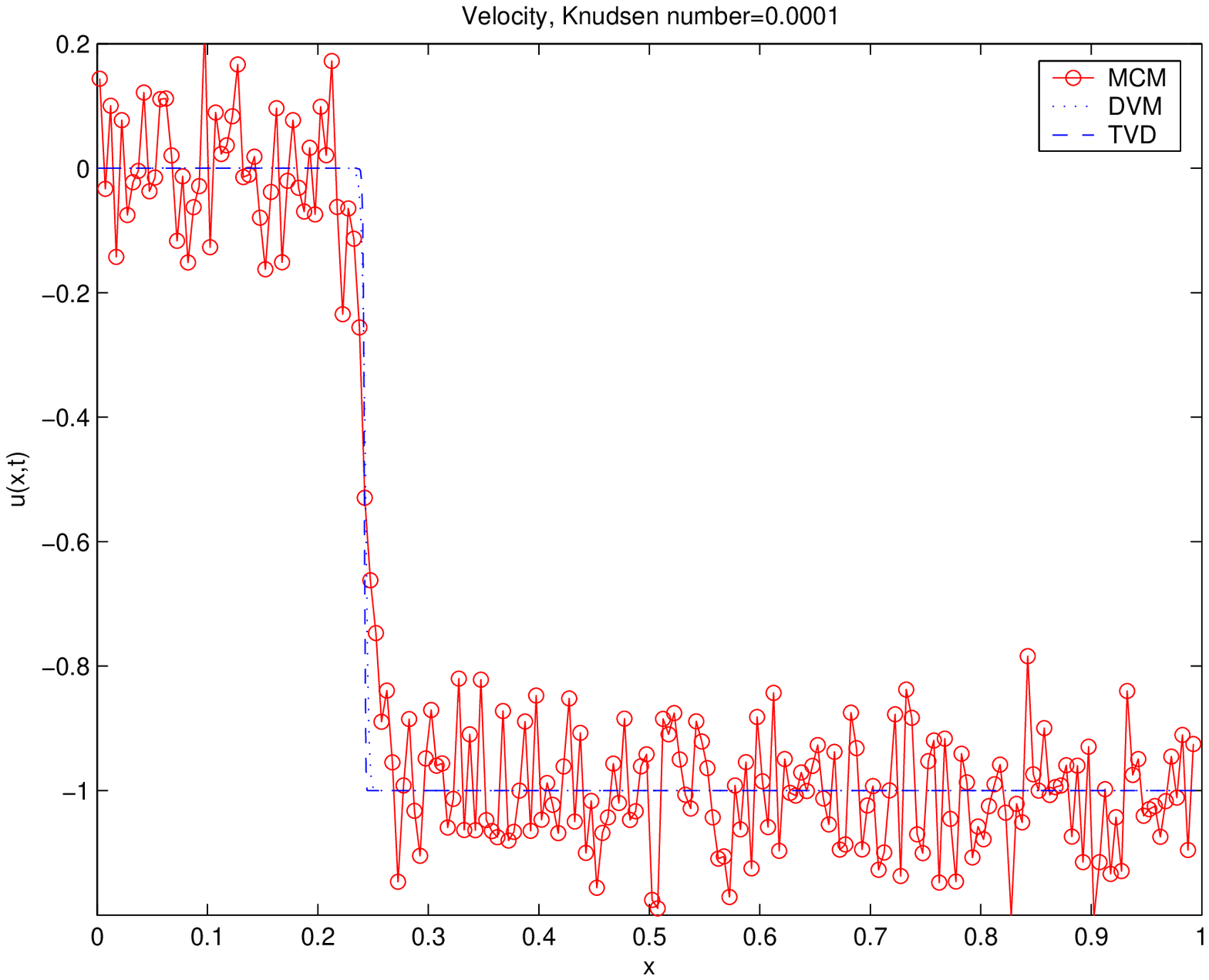}
\includegraphics[scale=0.40]{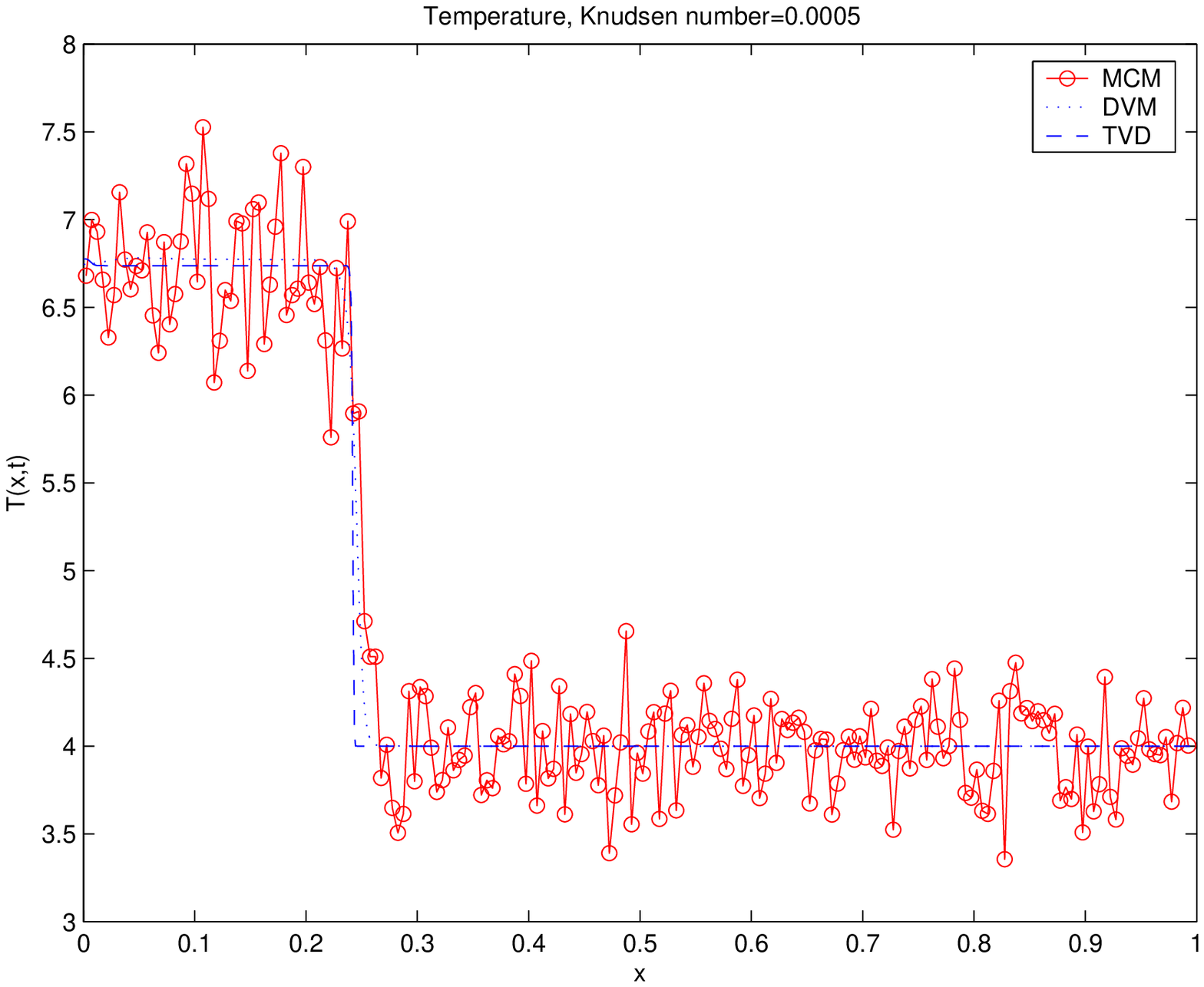}
\includegraphics[scale=0.40]{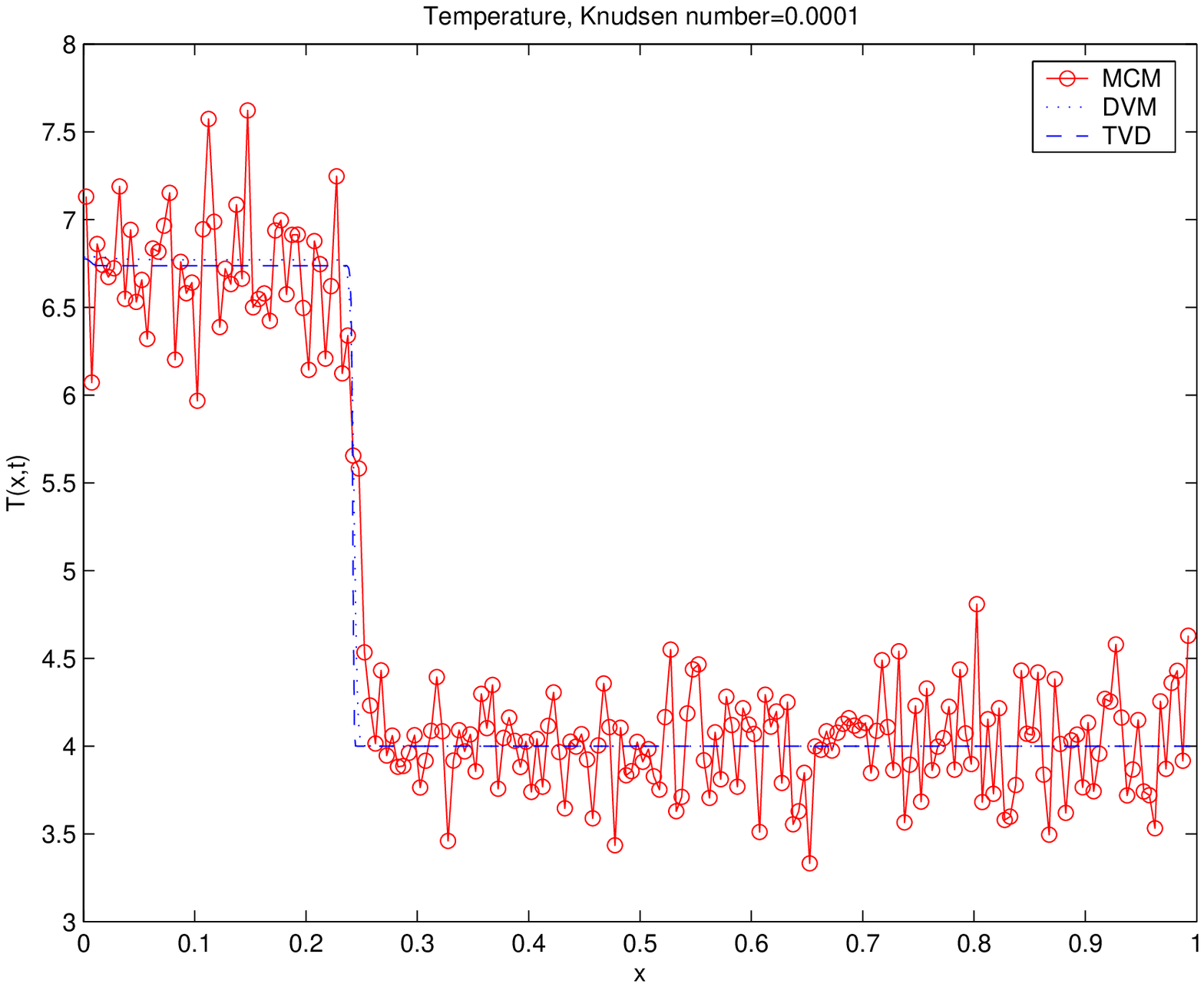}
\caption{Unsteady Shock. Solution at $t=0.065$ for MCM with
Knudsen numbers $\varepsilon=5\times 10^{-4}$ (left) and
$\varepsilon=10^{-4}$ (right). From top to bottom density, mean
velocity and temperature.} \label{US6}
\end{center}
\end{figure}

\subsection{1-D Lax Shock Tube Test}
\begin{figure}
\begin{center}
\includegraphics[scale=0.40]{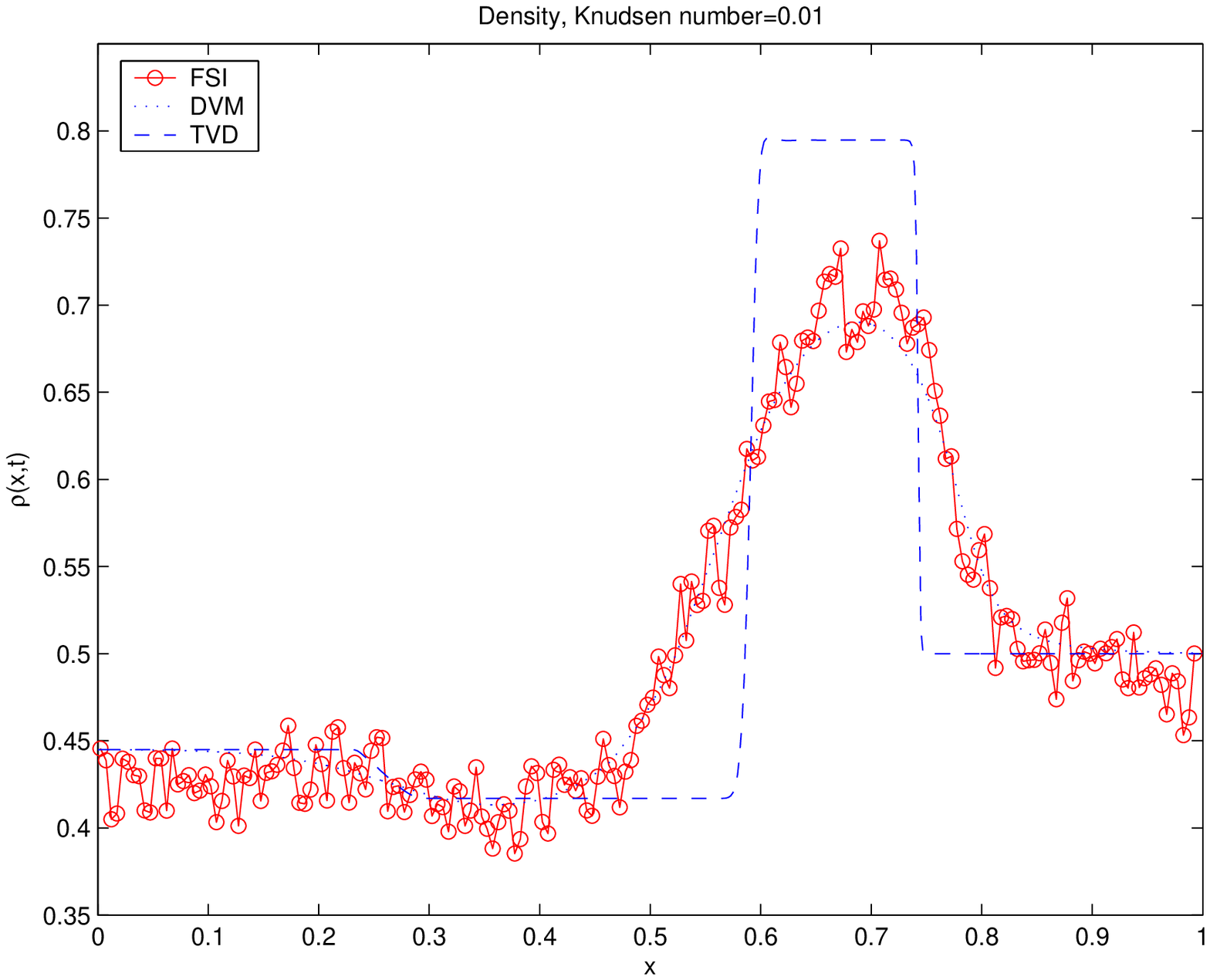}
\includegraphics[scale=0.40]{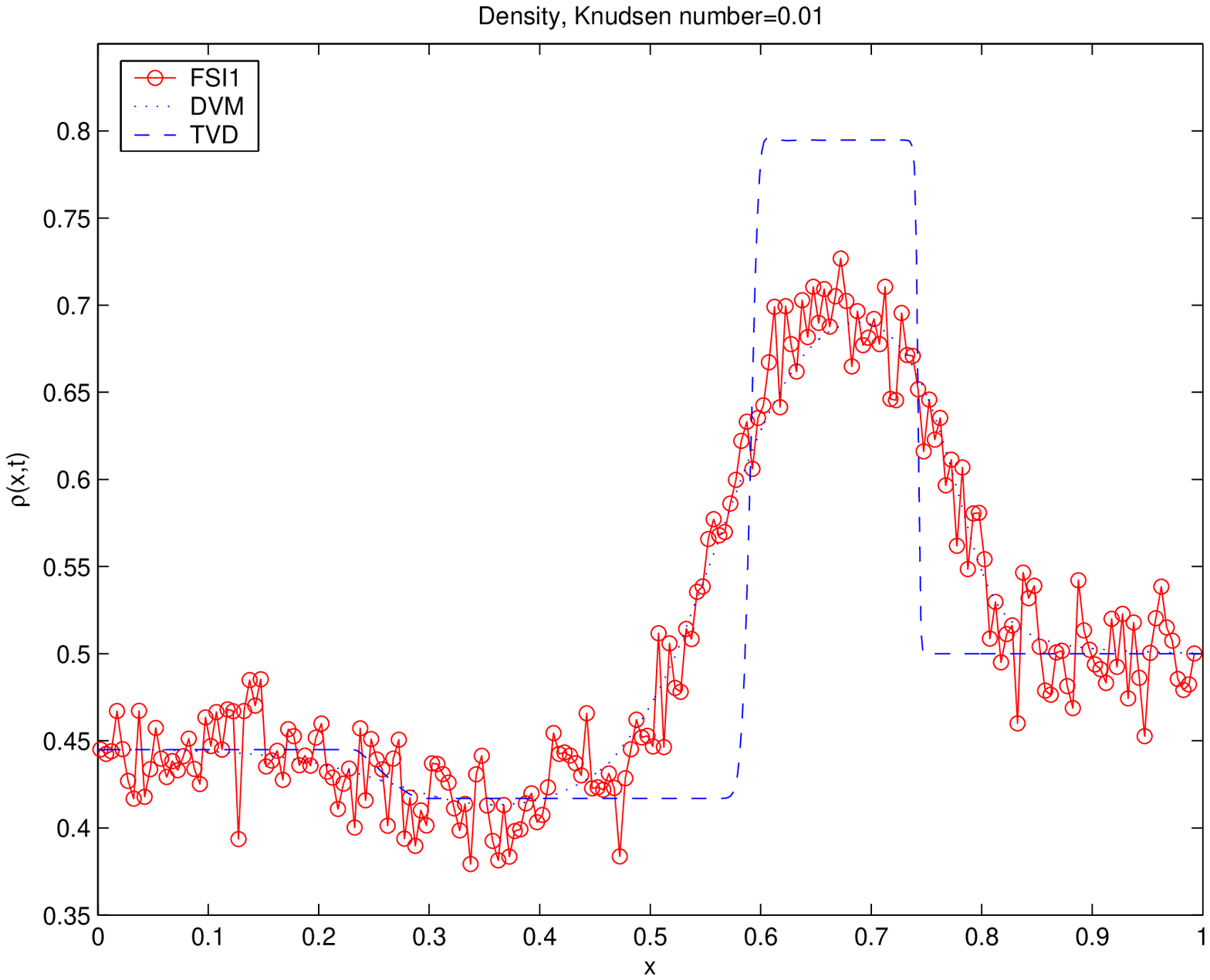}
\includegraphics[scale=0.40]{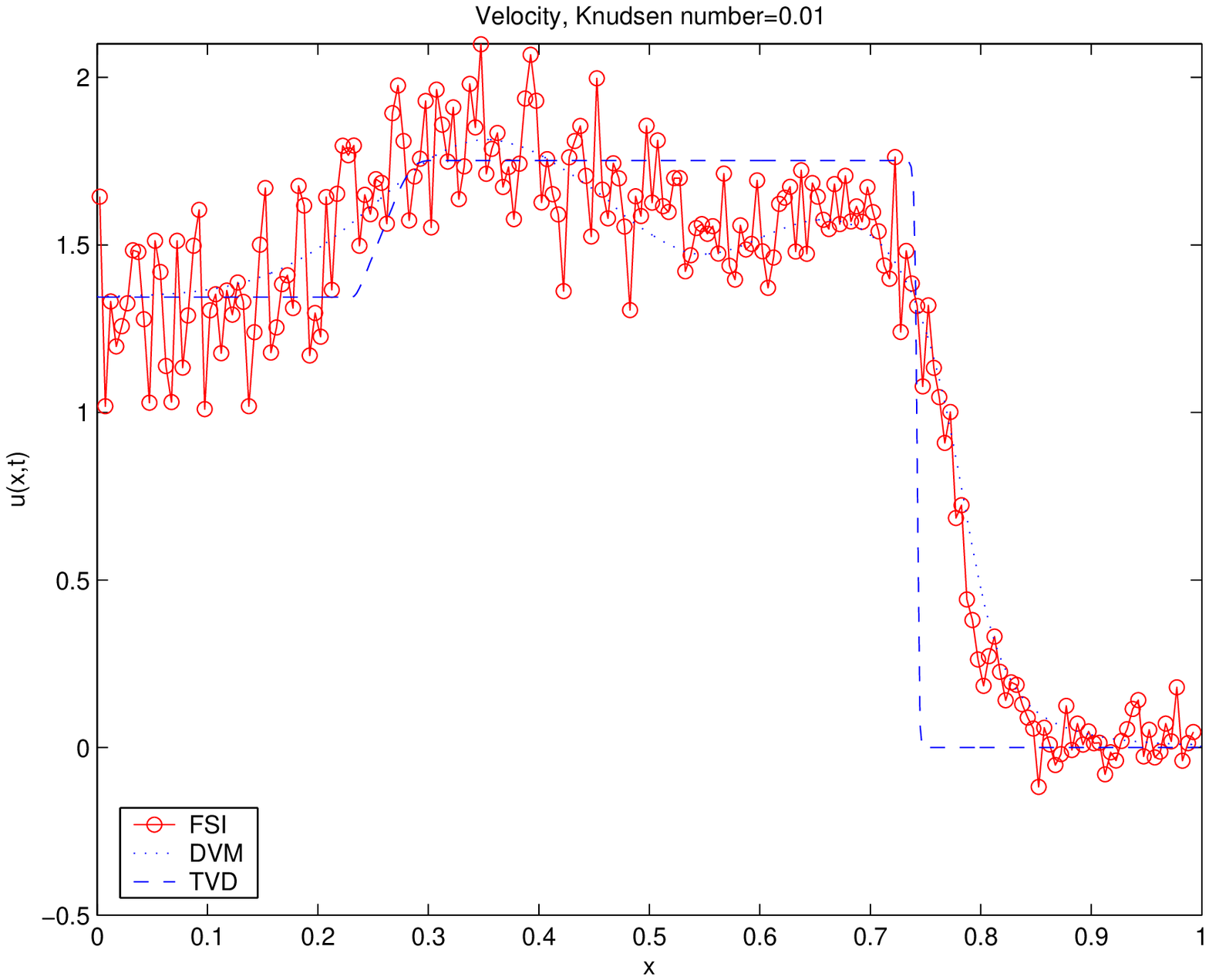}
\includegraphics[scale=0.40]{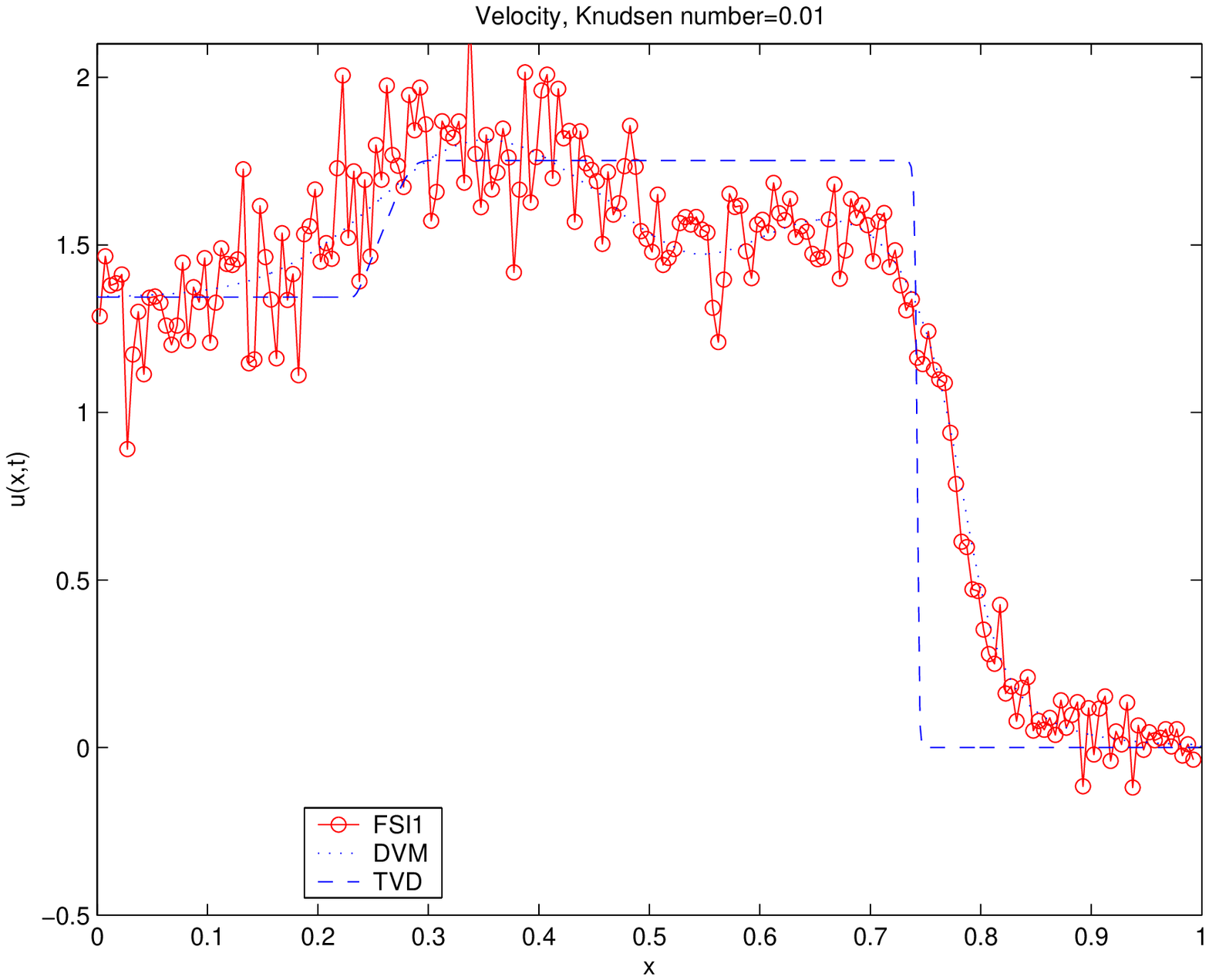}
\includegraphics[scale=0.40]{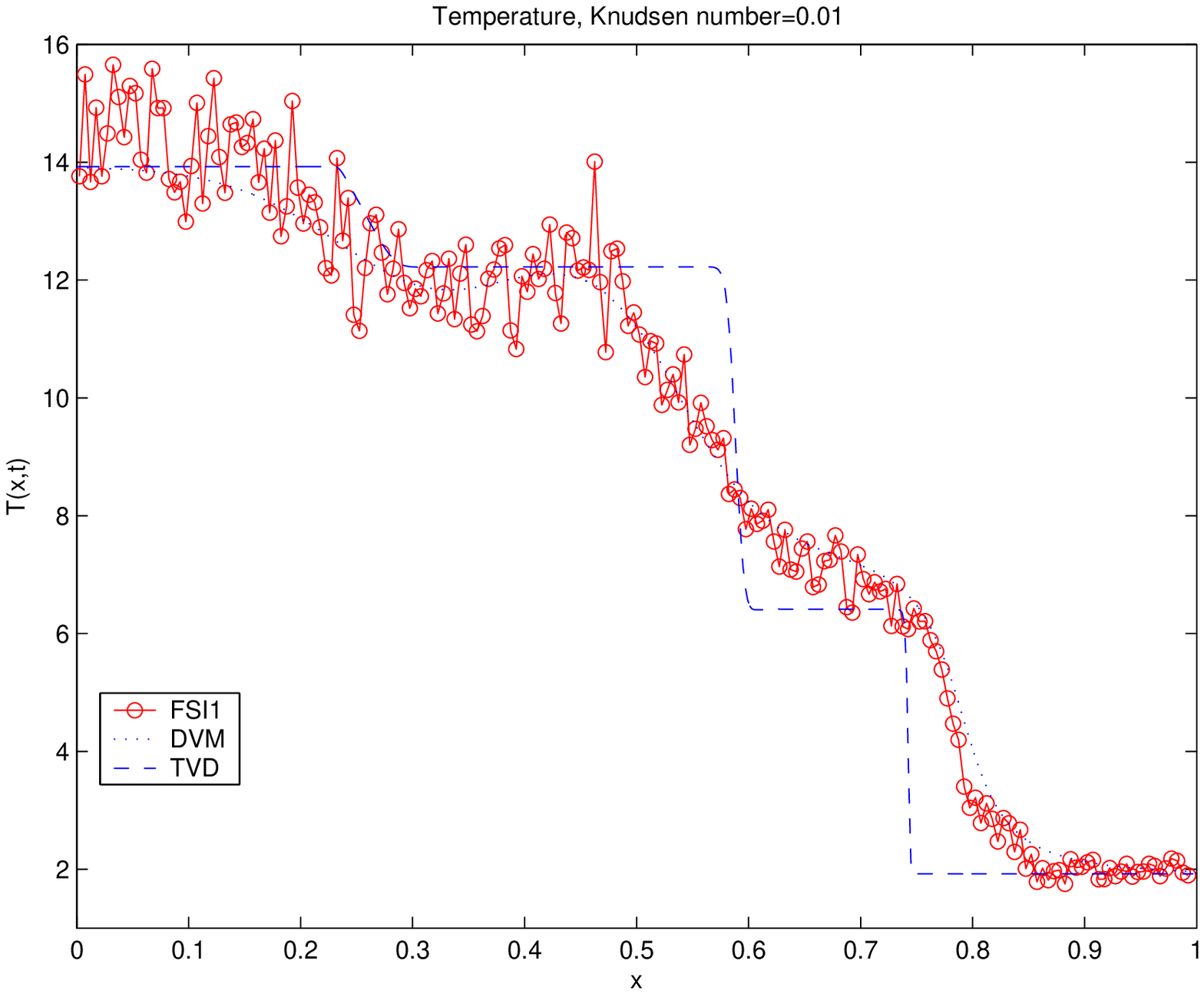}
\includegraphics[scale=0.40]{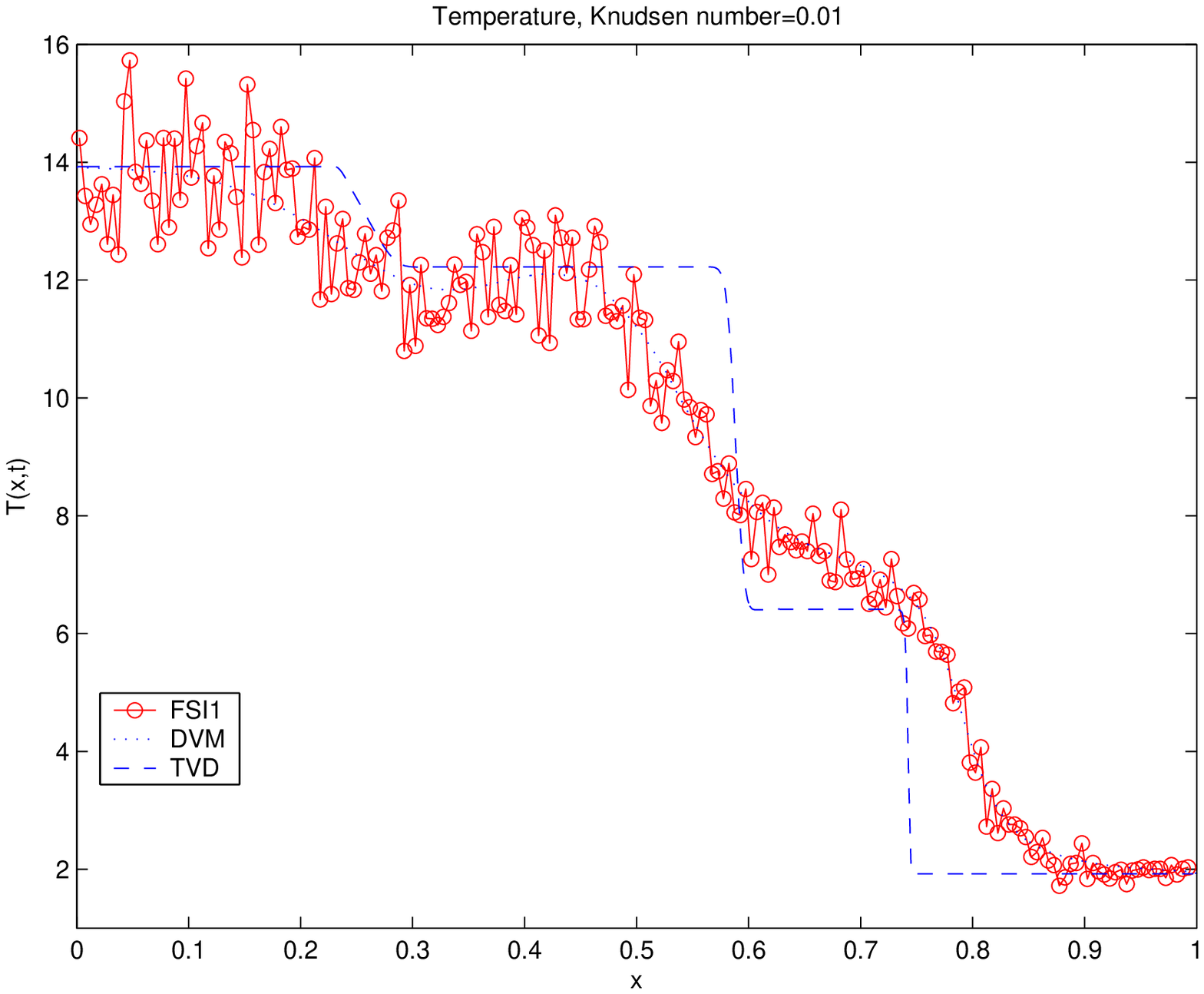}
\caption{Lax test: $\varepsilon=10^{-2}$. Solution at $t=0.05$ for
FSI (left) FSI1 (right). From top to bottom density, mean velocity
and temperature.} \label{L2}
\end{center}
\end{figure}

\begin{figure}
\begin{center}
\includegraphics[scale=0.40]{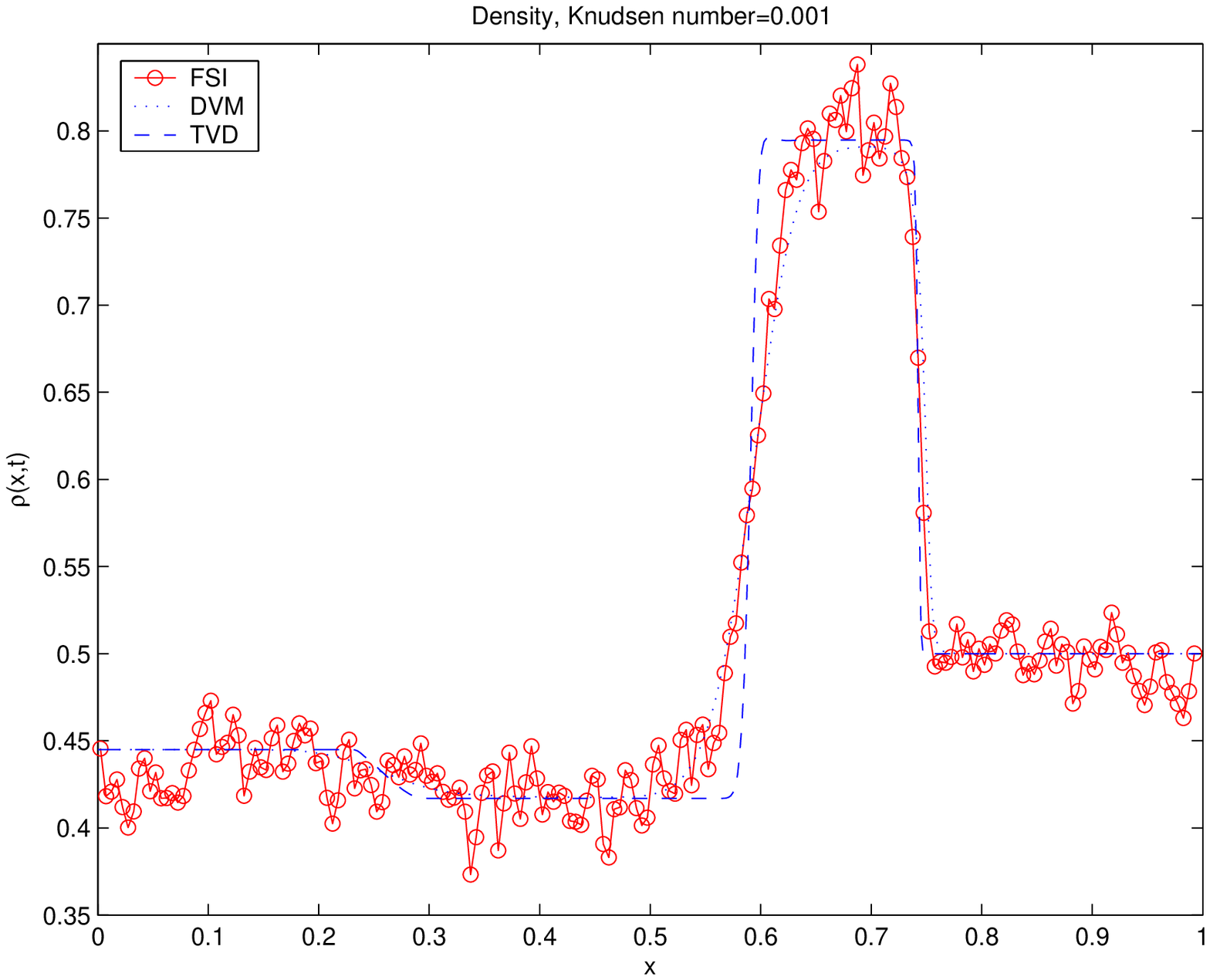}
\includegraphics[scale=0.40]{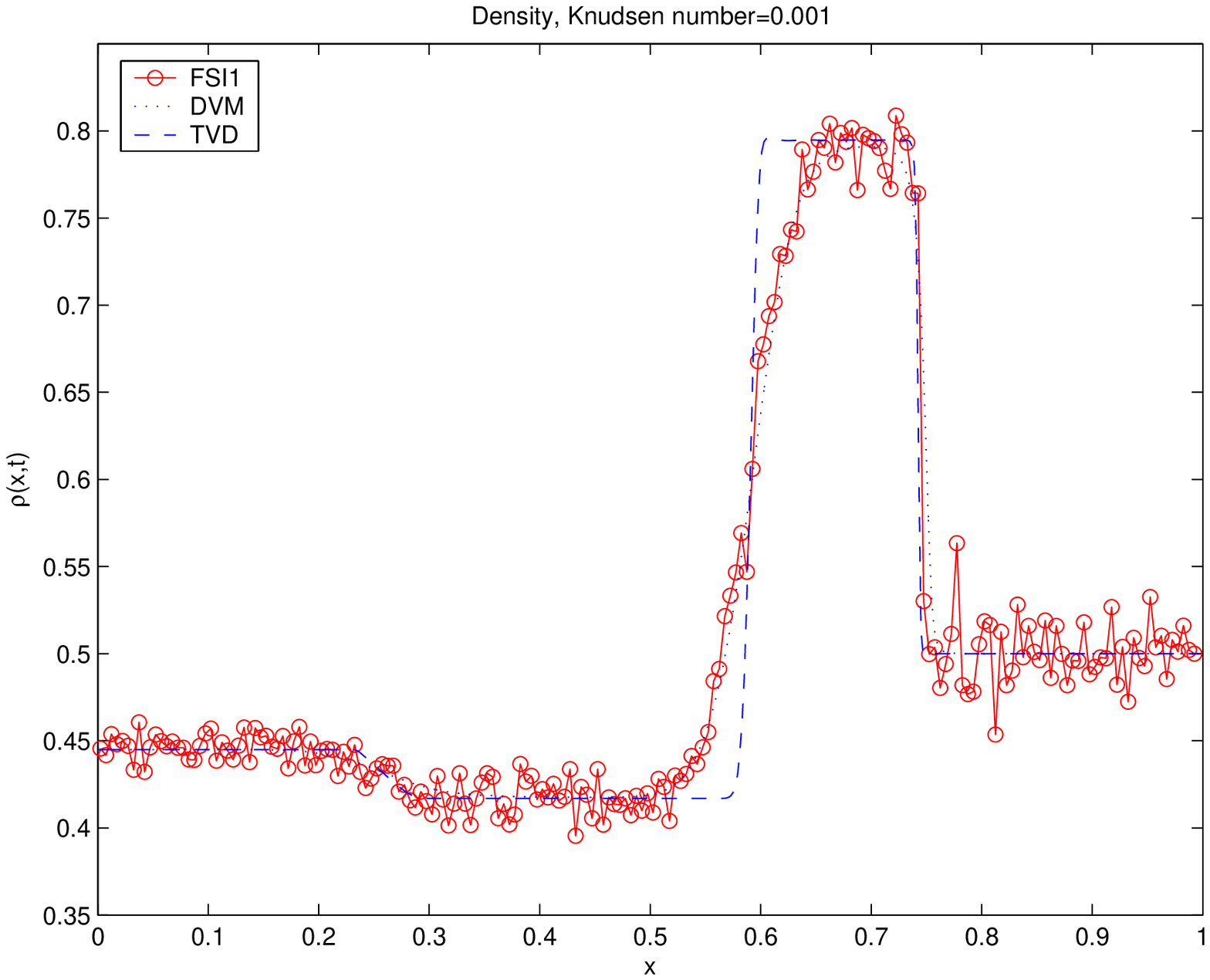}
\includegraphics[scale=0.40]{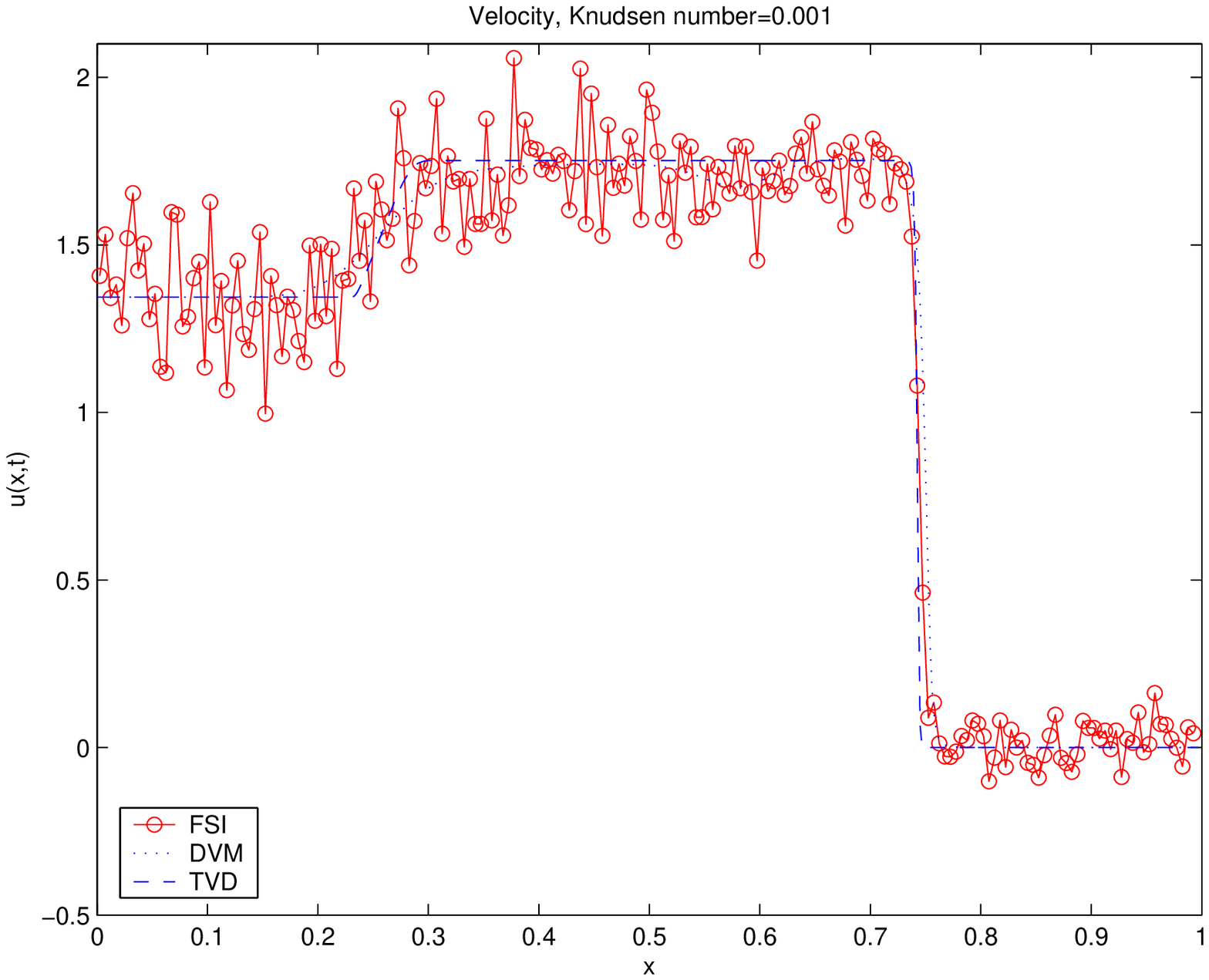}
\includegraphics[scale=0.40]{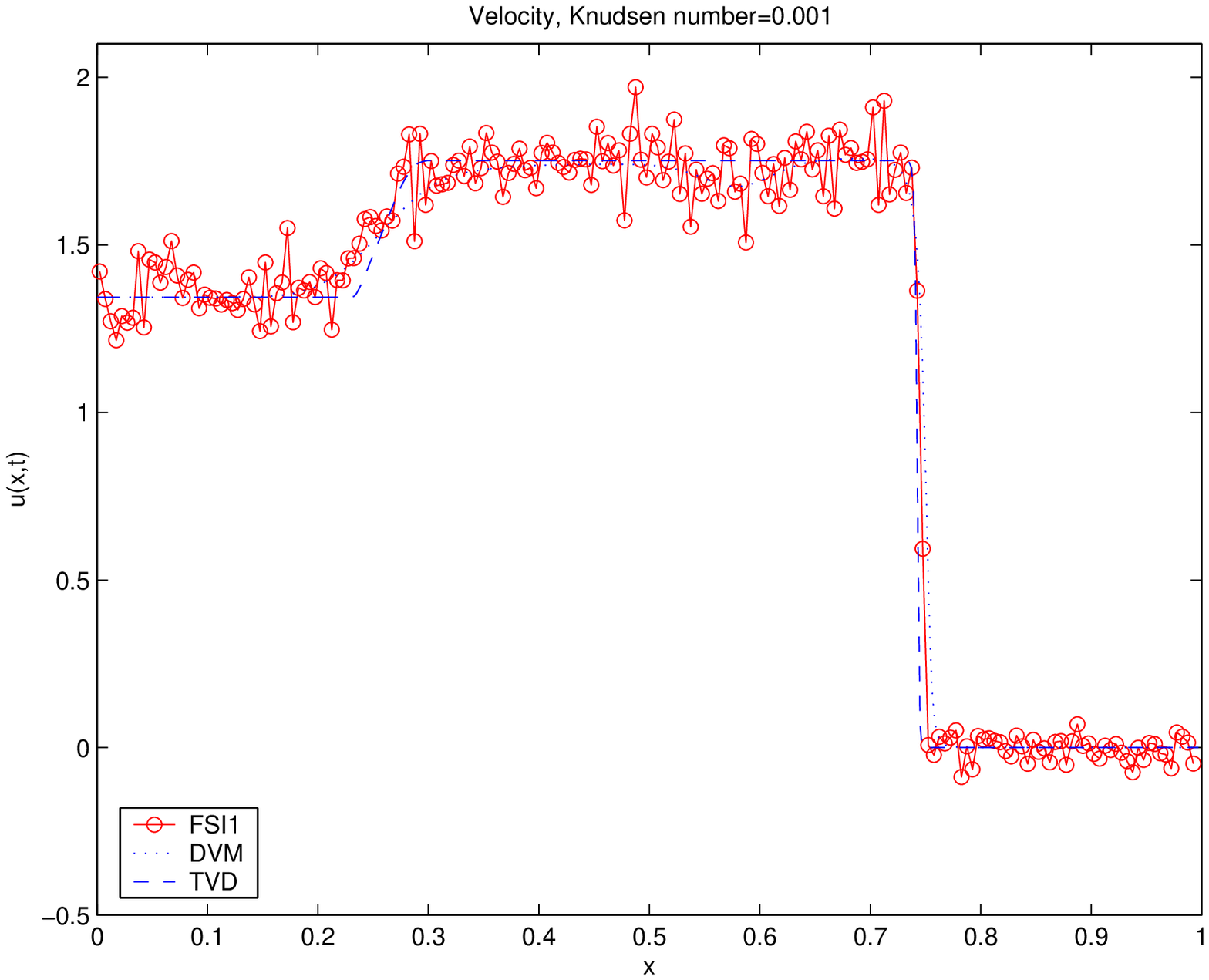}
\includegraphics[scale=0.40]{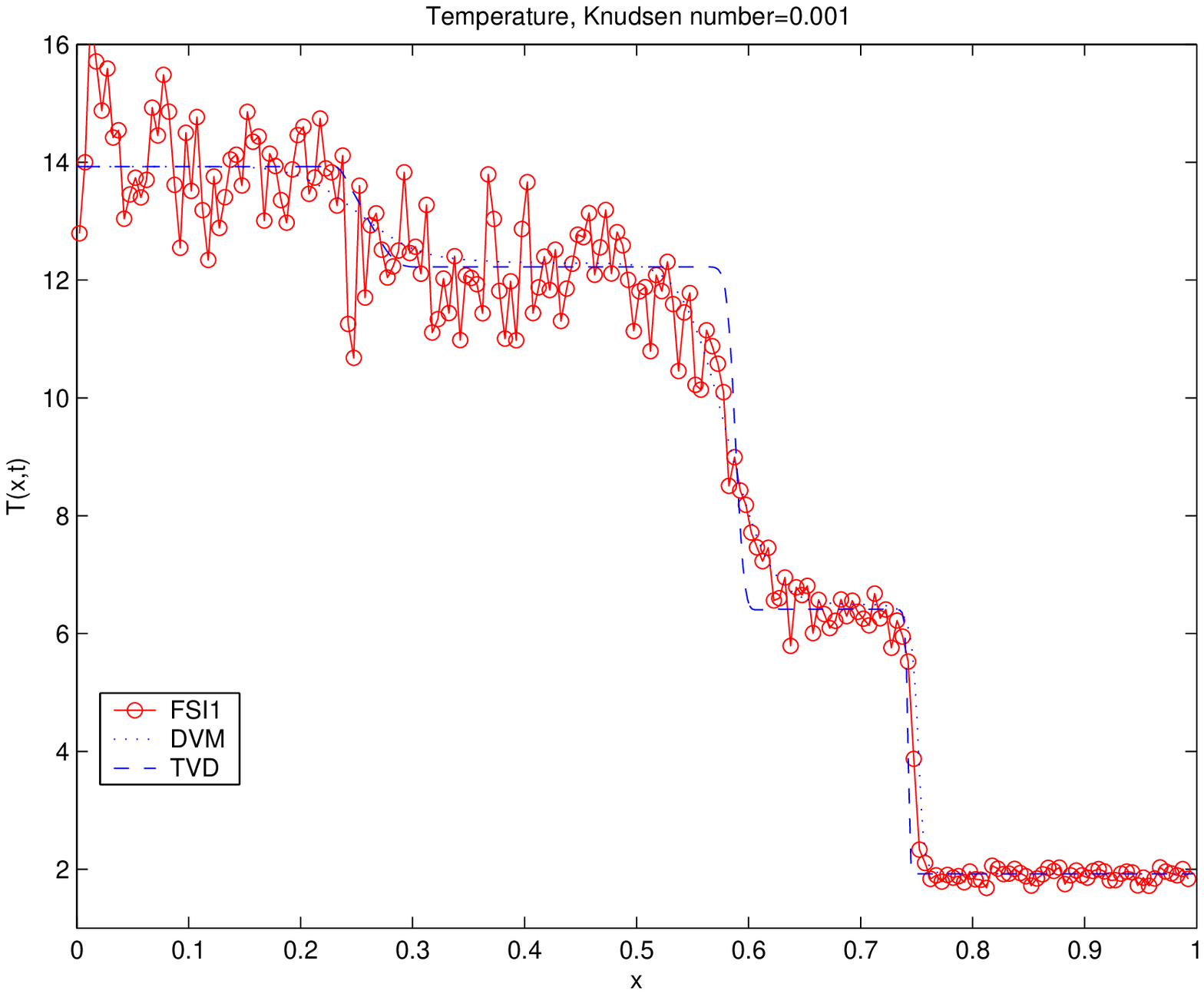}
\includegraphics[scale=0.40]{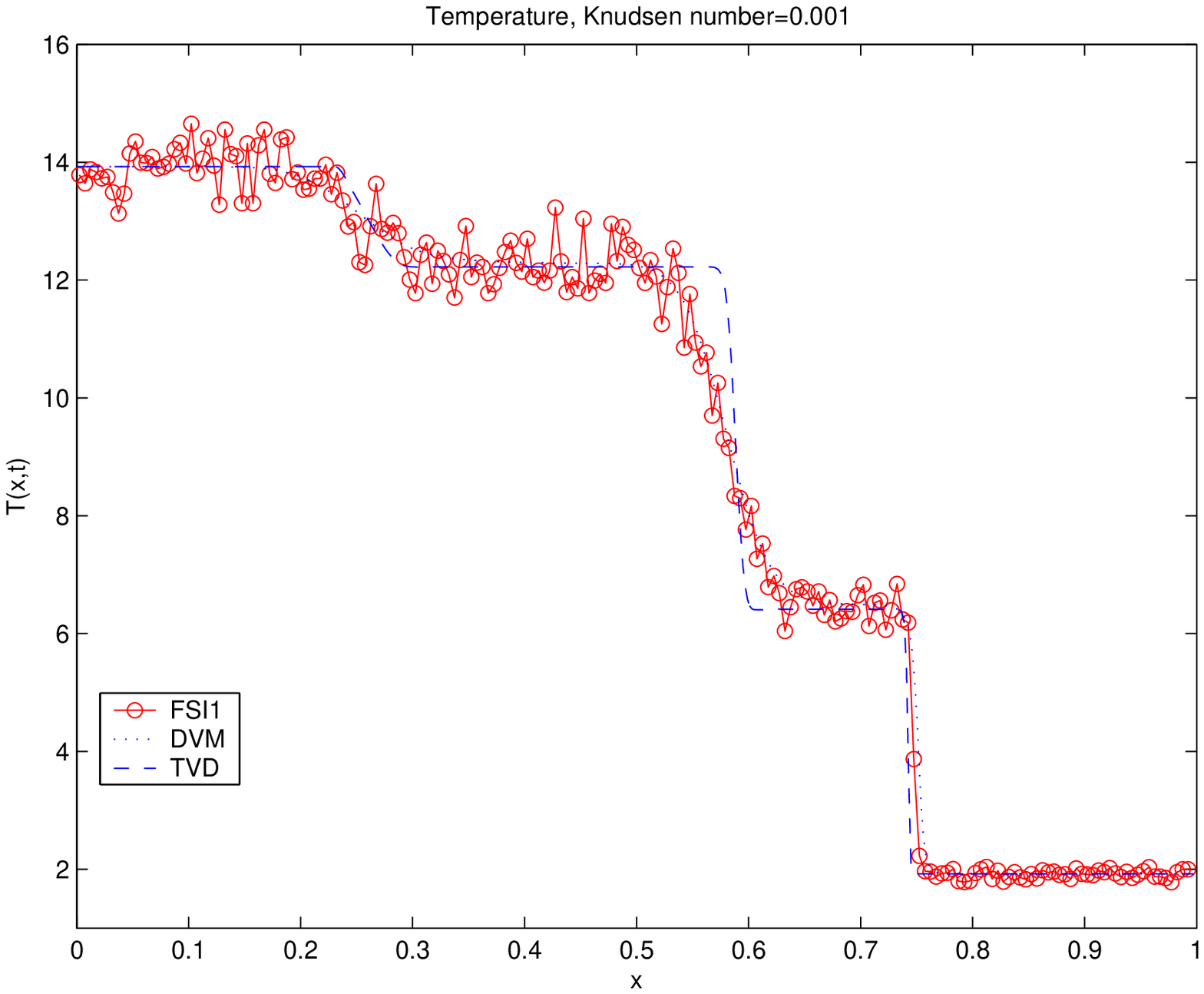}
\caption{Lax test: $\varepsilon=10^{-3}$. Solution at $t=0.05$ for
FSI (left) FSI1 (right). From top to bottom density, mean velocity
and temperature.} \label{L3}
\end{center}
\end{figure}

\begin{figure}
\begin{center}
\includegraphics[scale=0.40]{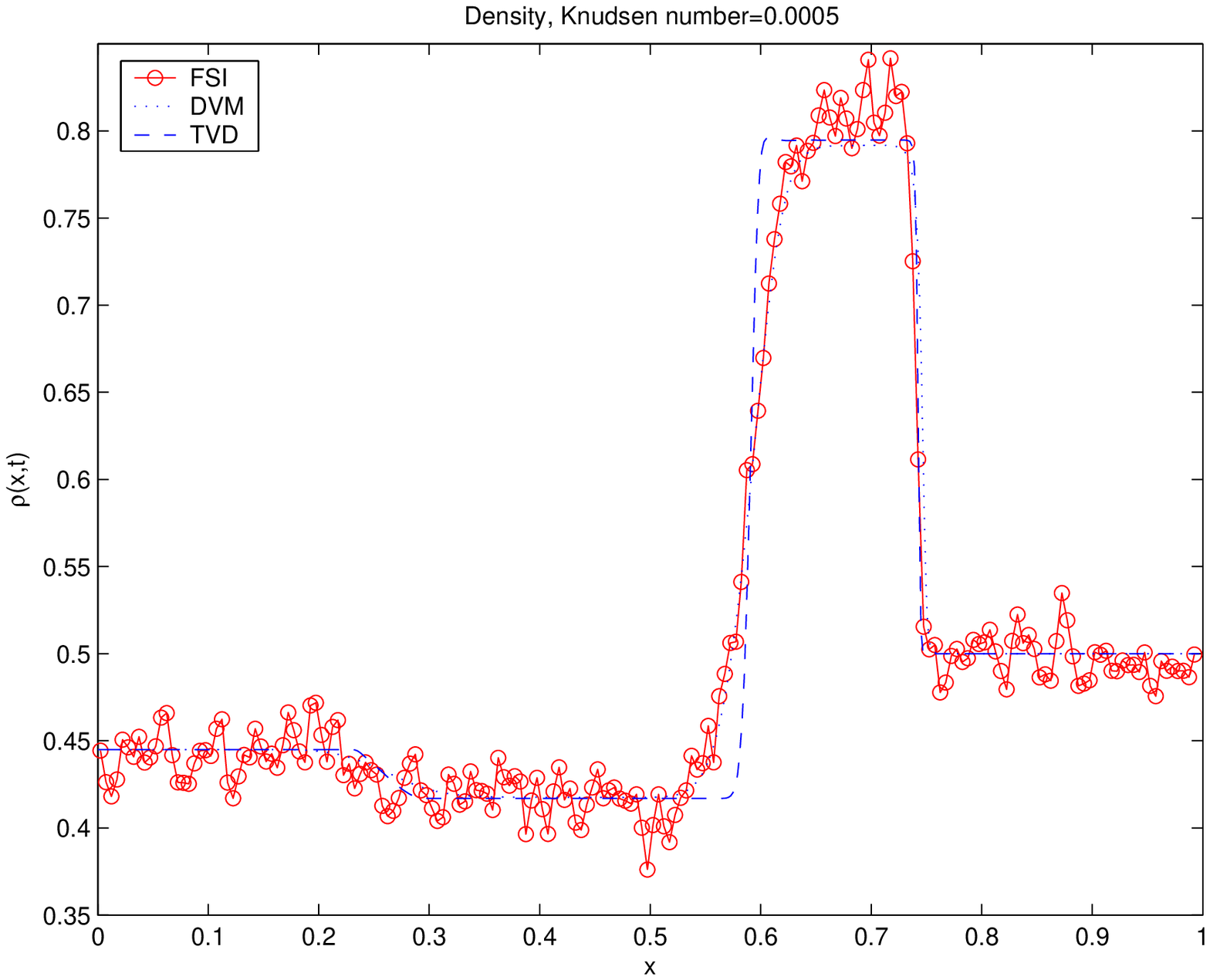}
\includegraphics[scale=0.40]{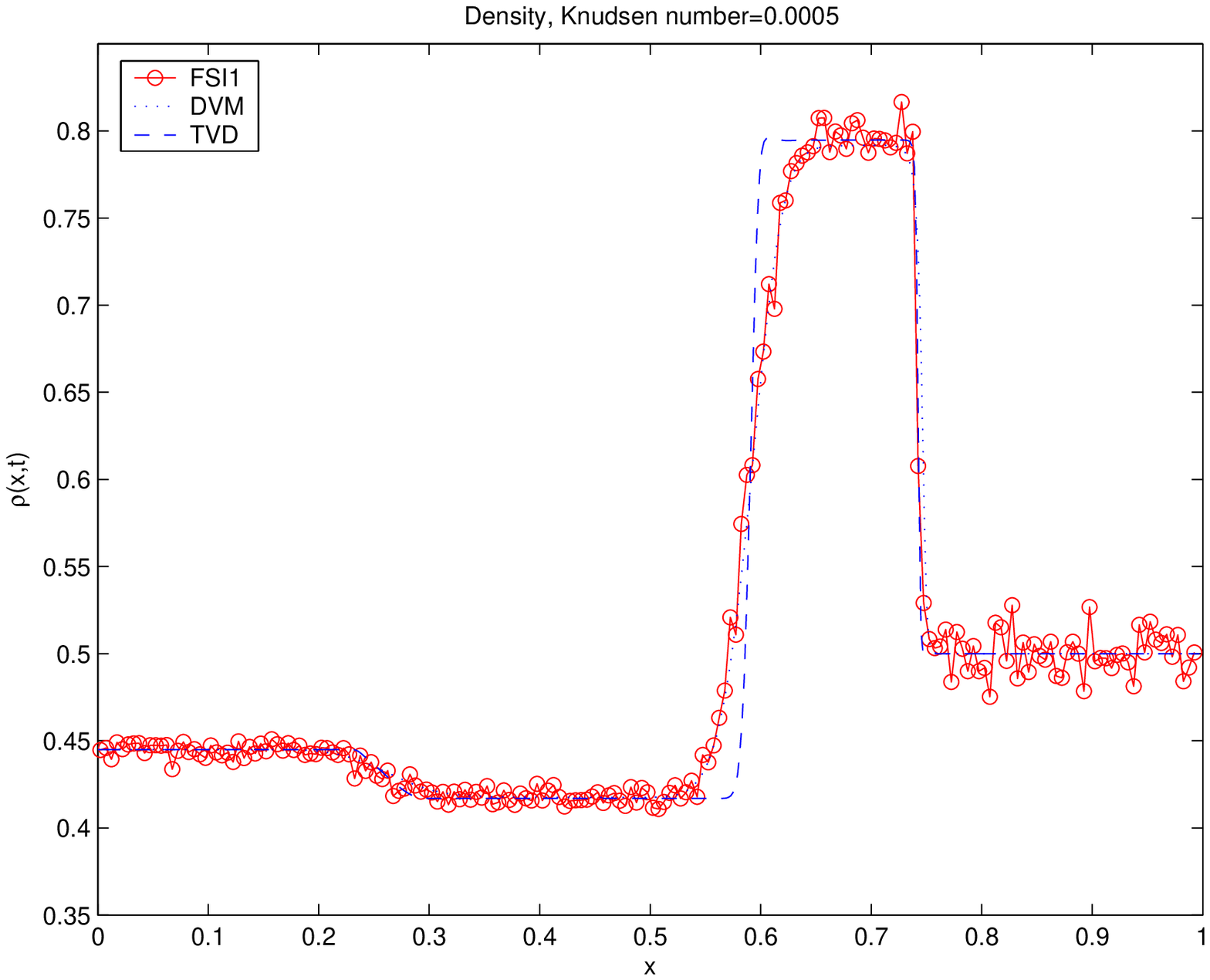}
\includegraphics[scale=0.40]{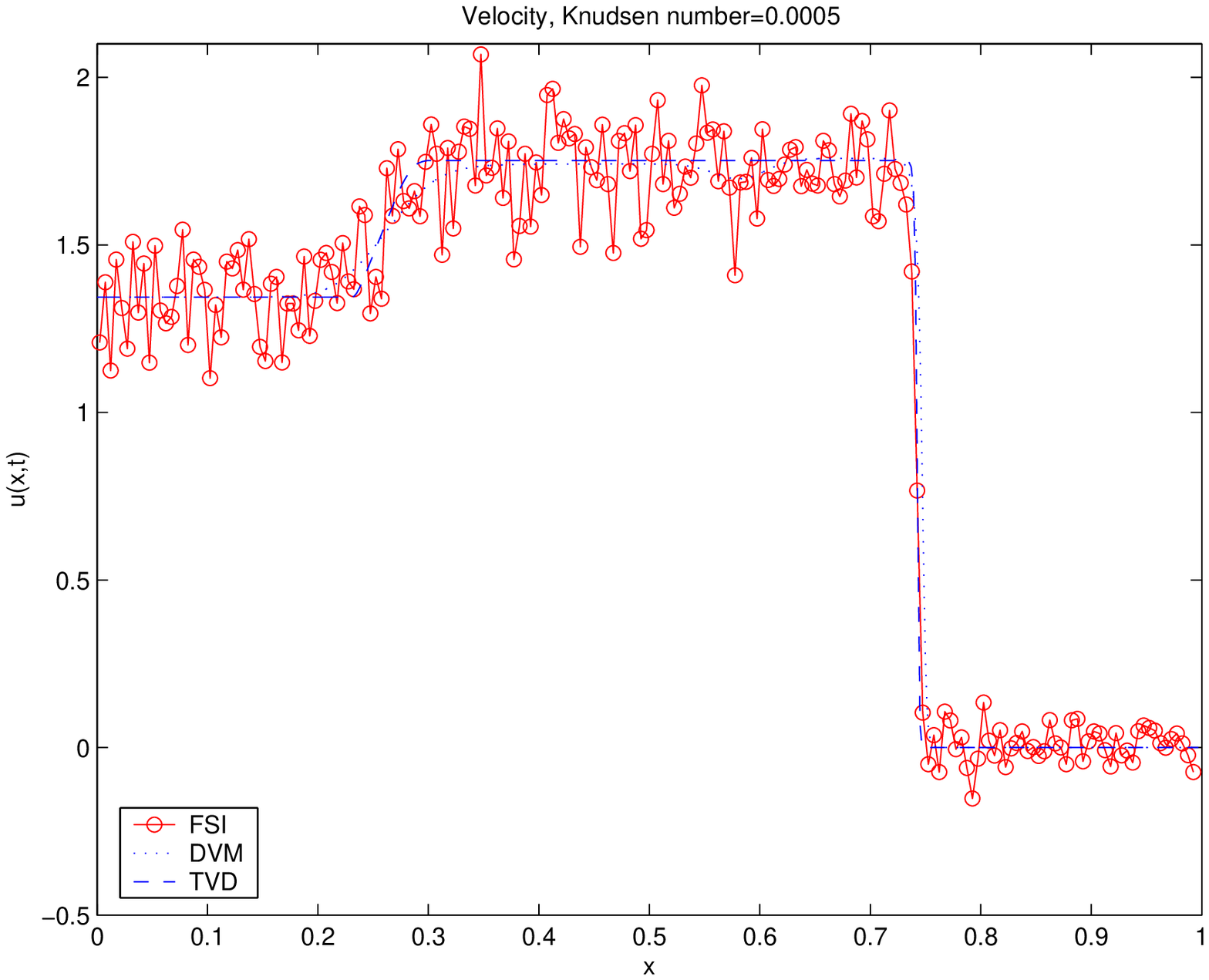}
\includegraphics[scale=0.40]{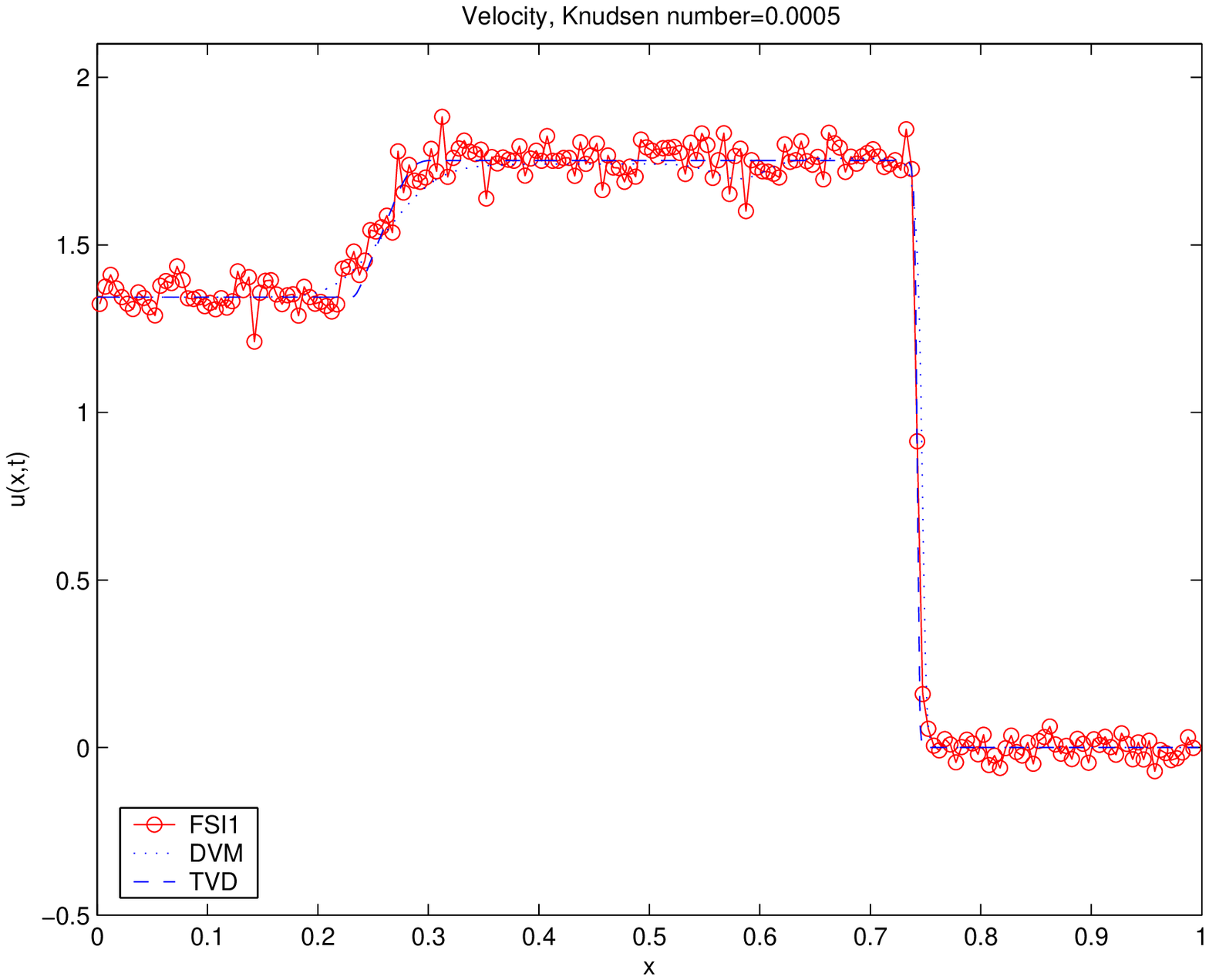}
\includegraphics[scale=0.40]{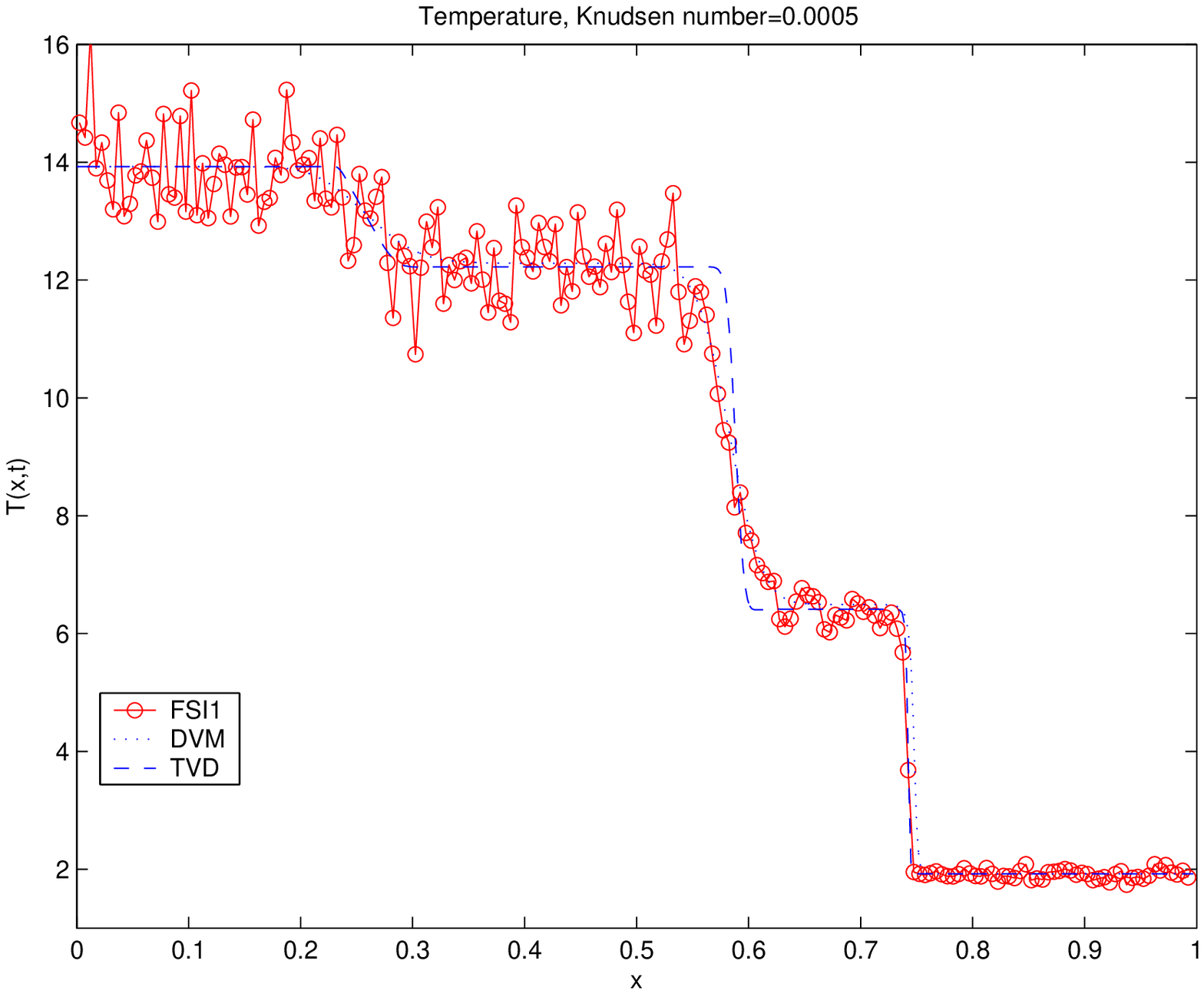}
\includegraphics[scale=0.40]{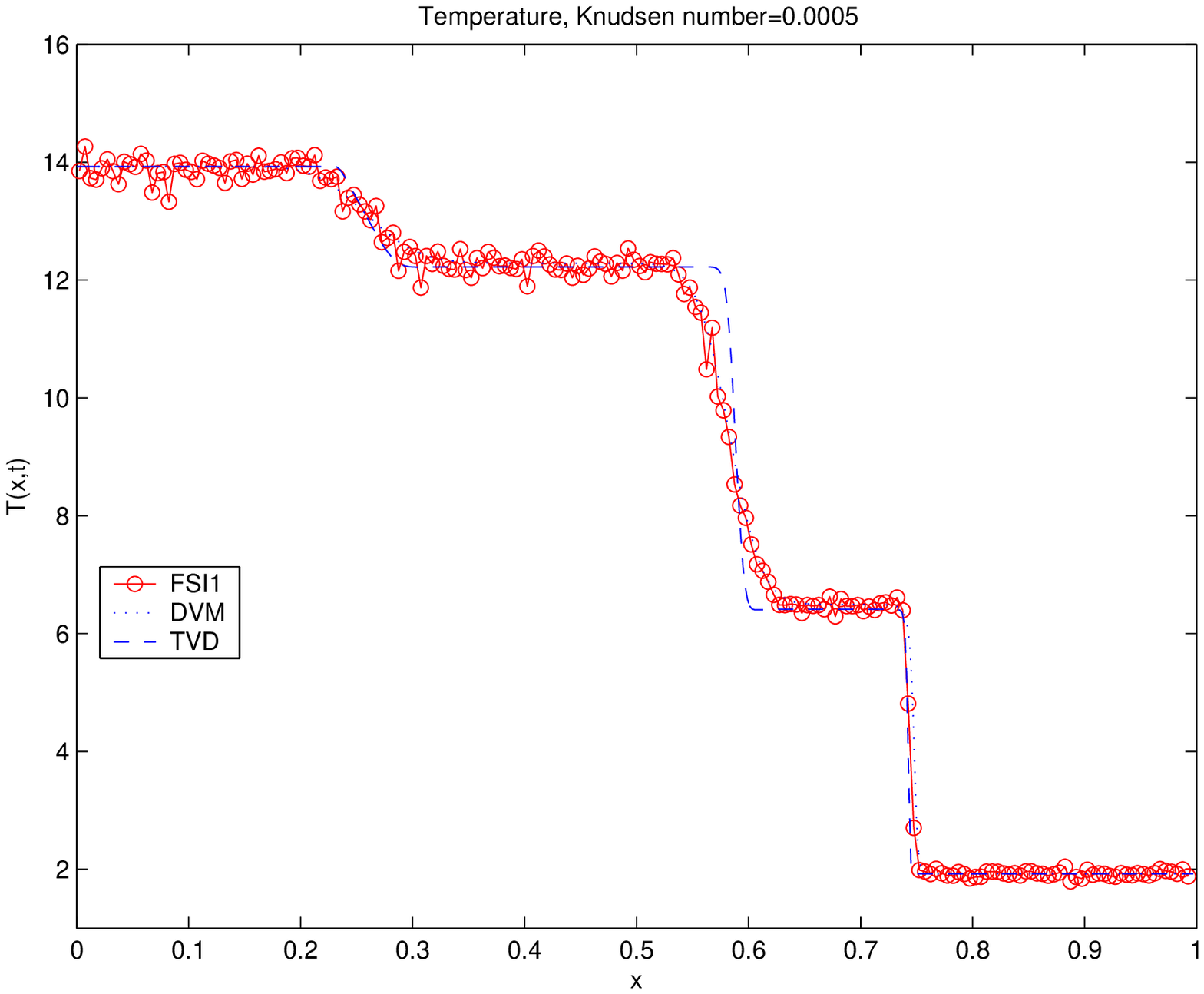}
\caption{Lax test: $\varepsilon=5\times 10^{-4}$. Solution at
$t=0.05$ for FSI (left) FSI1 (right). From top to bottom density,
mean velocity and temperature.} \label{L4}
\end{center}
\end{figure}

\begin{figure}
\begin{center}
\includegraphics[scale=0.40]{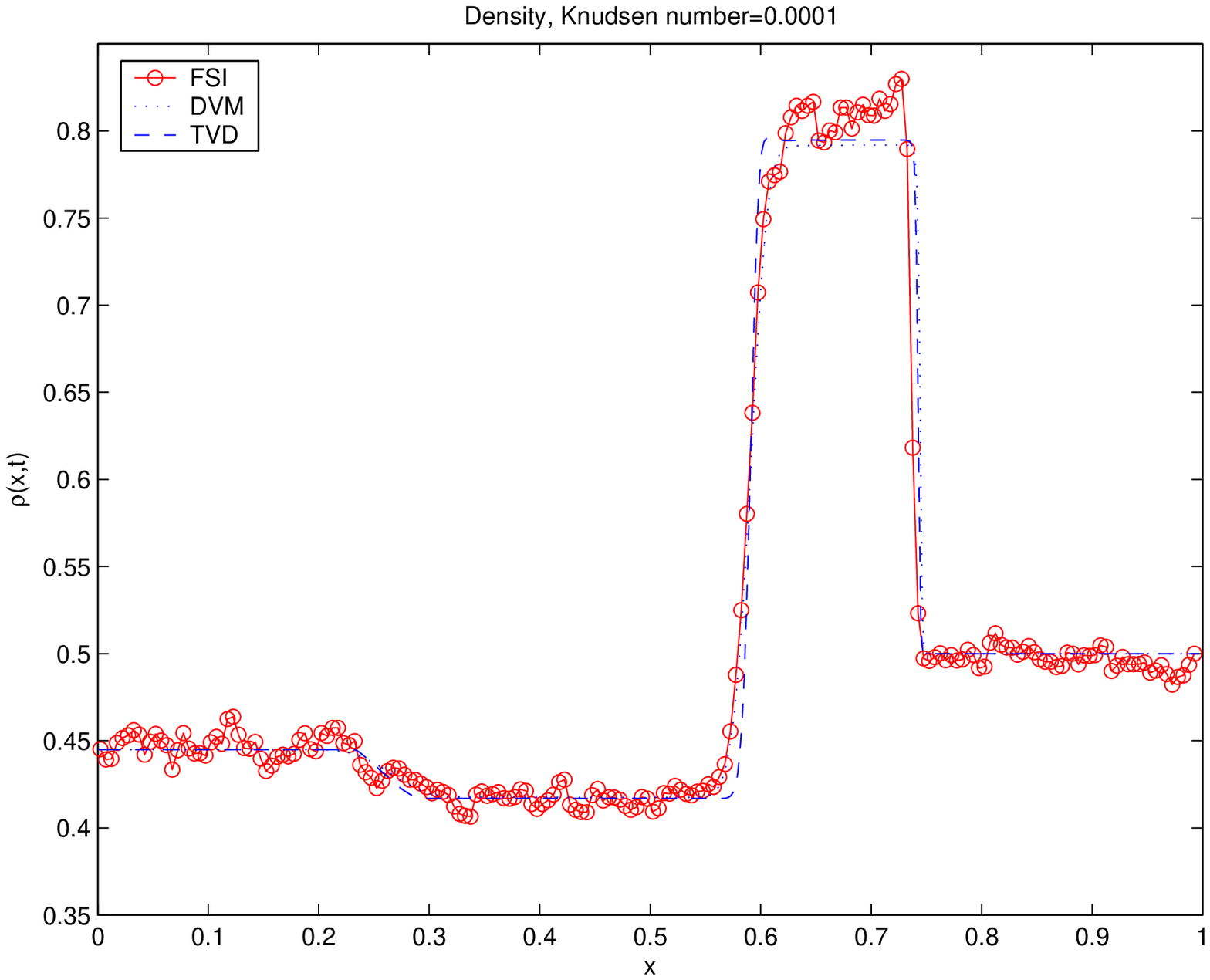}
\includegraphics[scale=0.40]{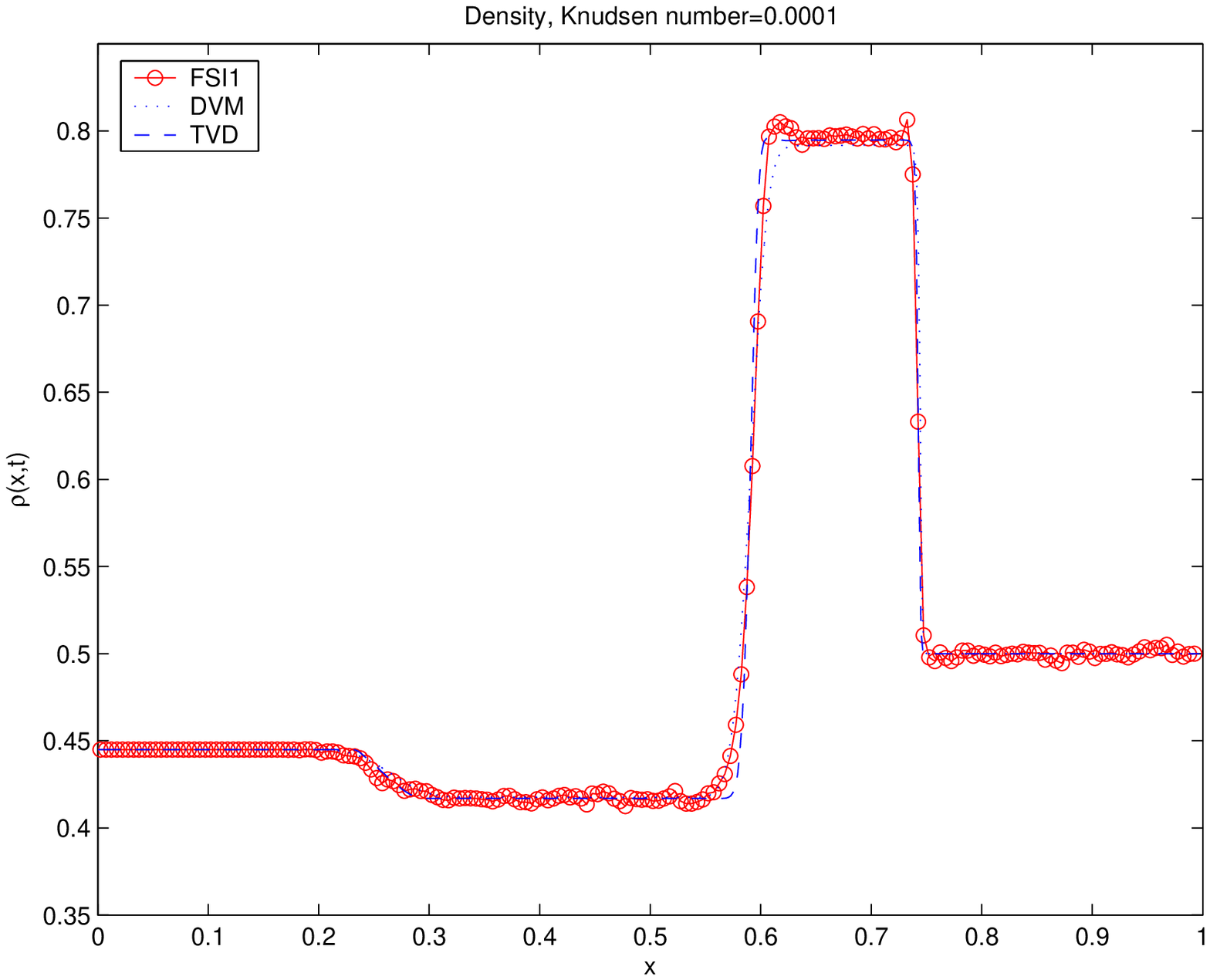}
\includegraphics[scale=0.40]{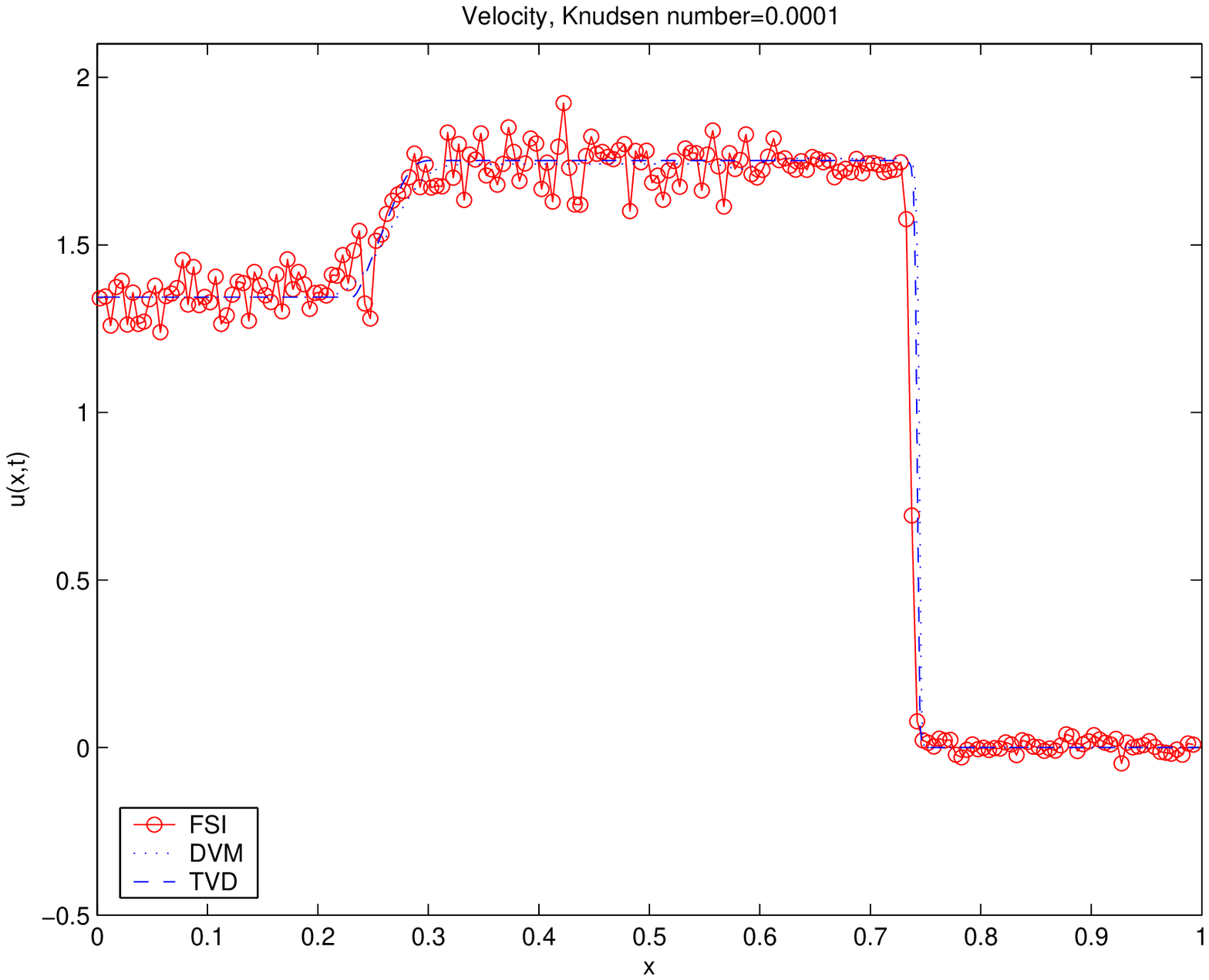}
\includegraphics[scale=0.40]{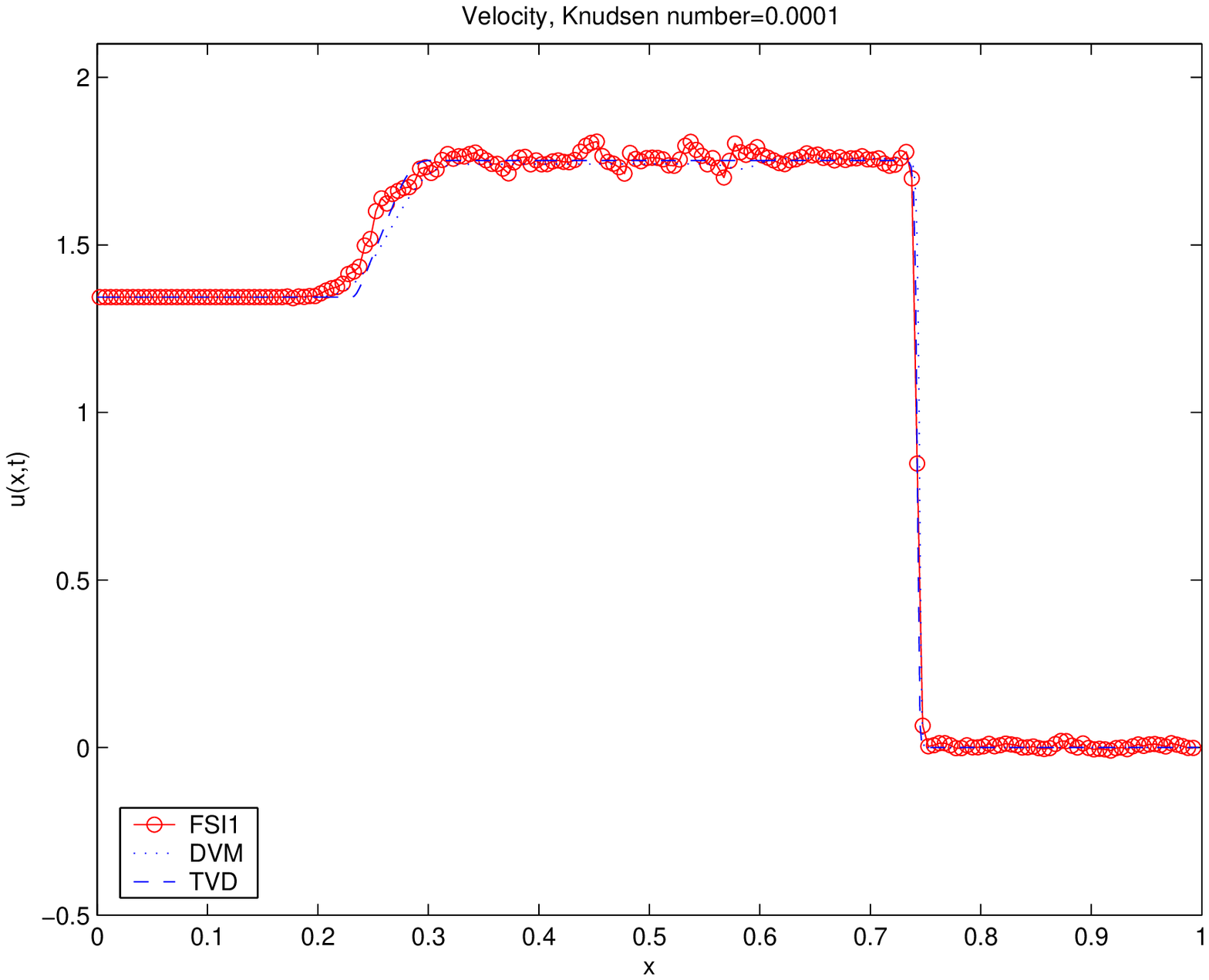}
\includegraphics[scale=0.40]{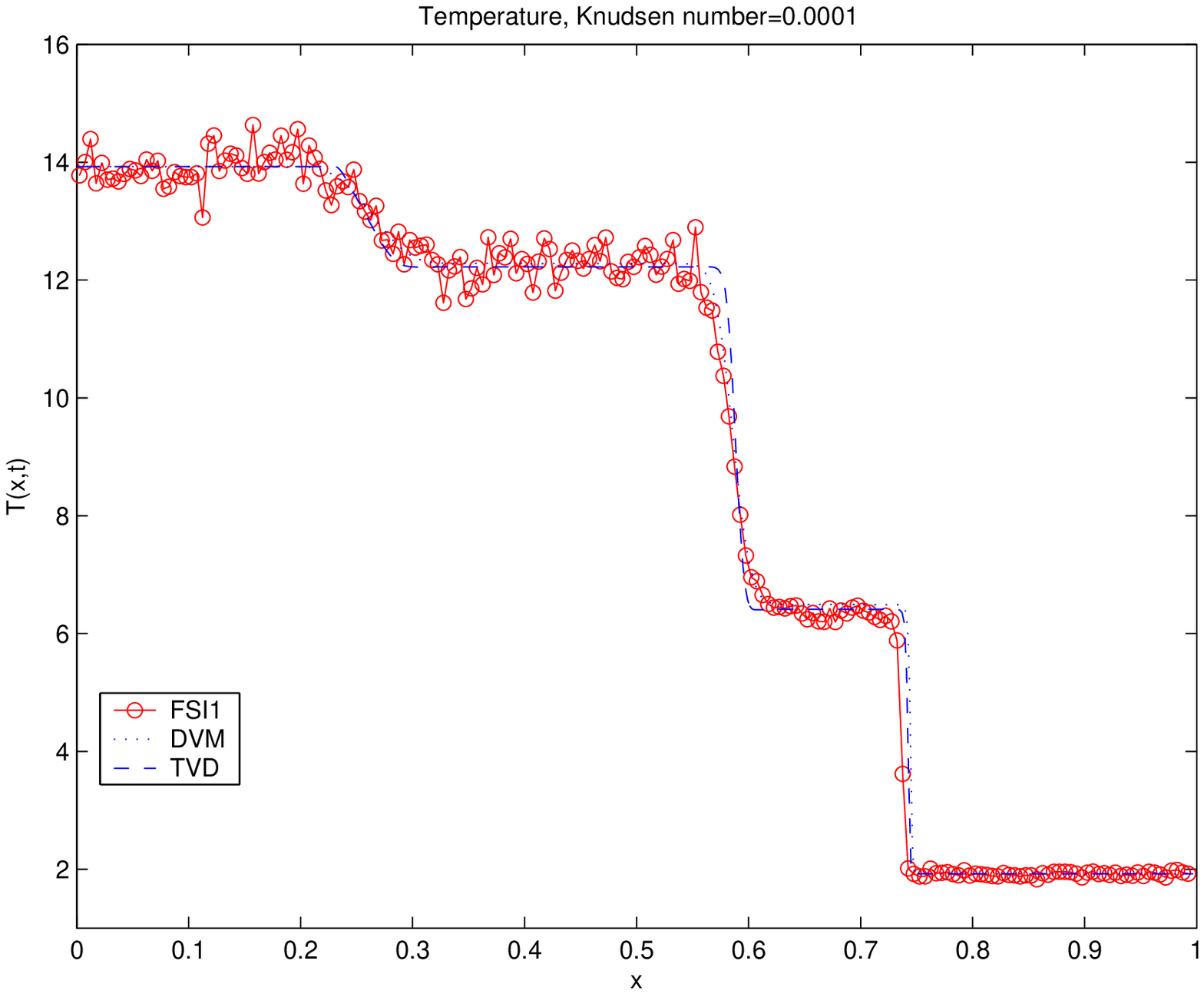}
\includegraphics[scale=0.40]{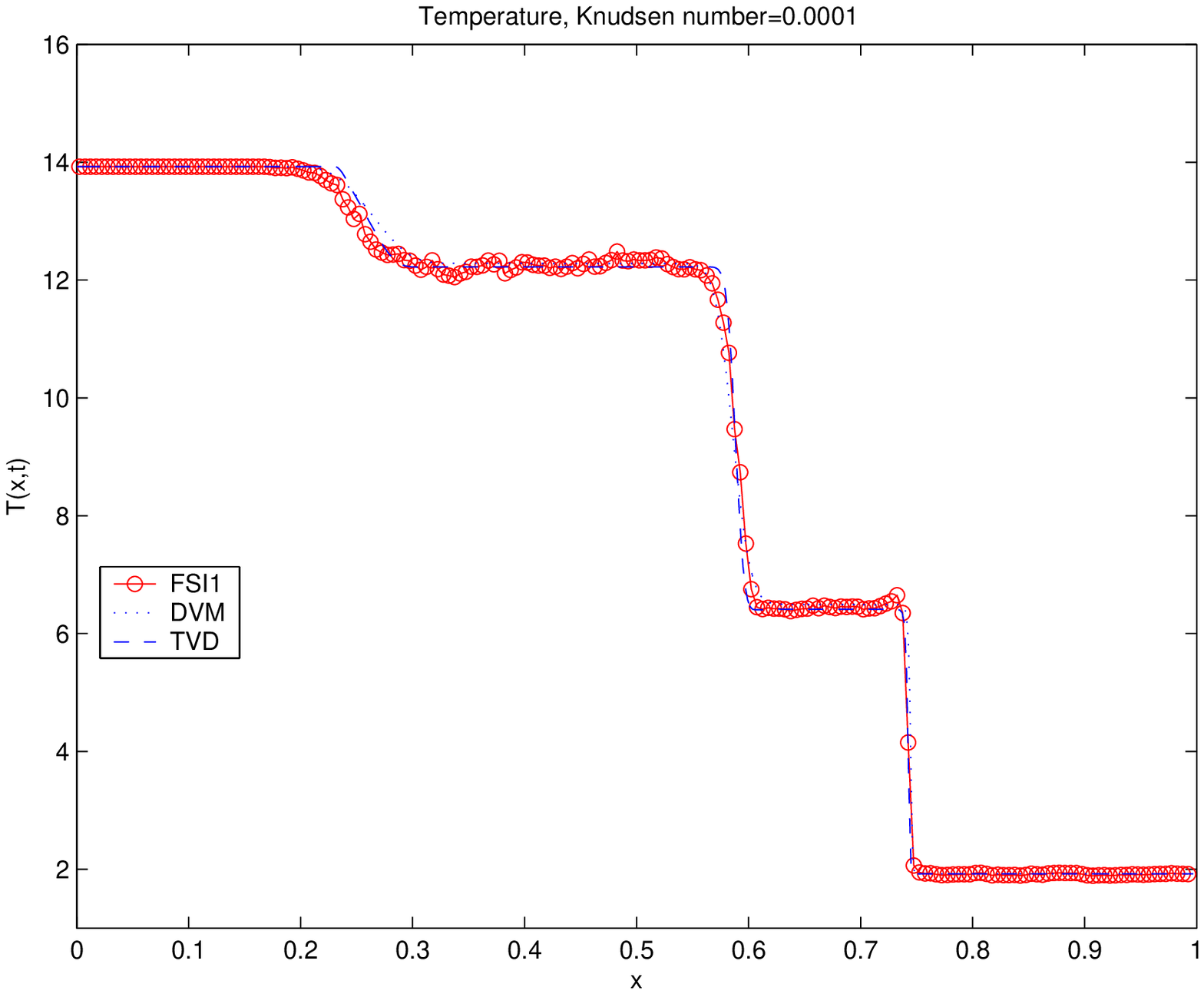}
\caption{Lax test: $\varepsilon=10^{-4}$. Solution at $t=0.05$ for
FSI (left) FSI1 (right). From top to bottom density, mean velocity
and temperature.} \label{L5}
\end{center}
\end{figure}

\begin{figure}
\begin{center}
\includegraphics[scale=0.40]{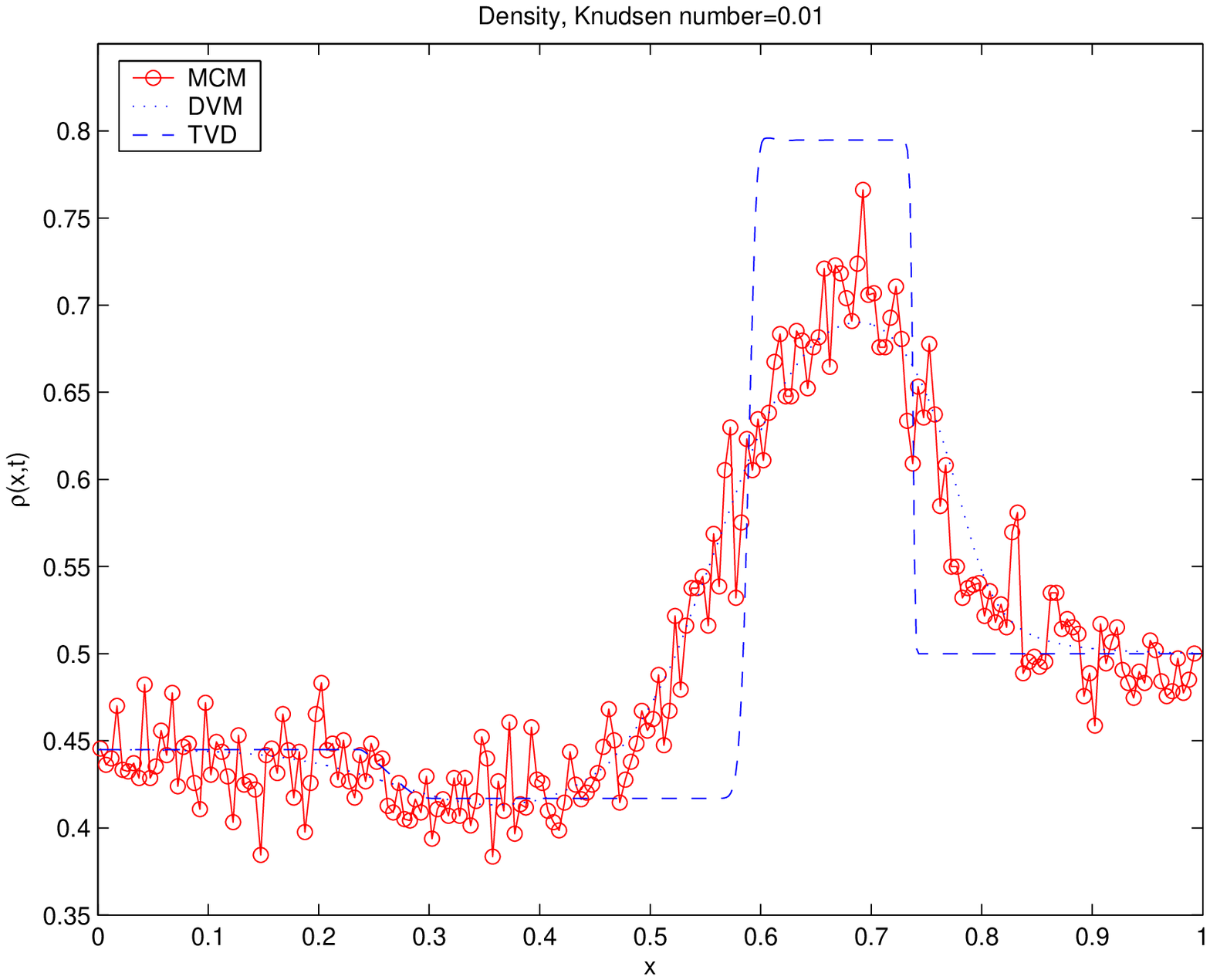}
\includegraphics[scale=0.40]{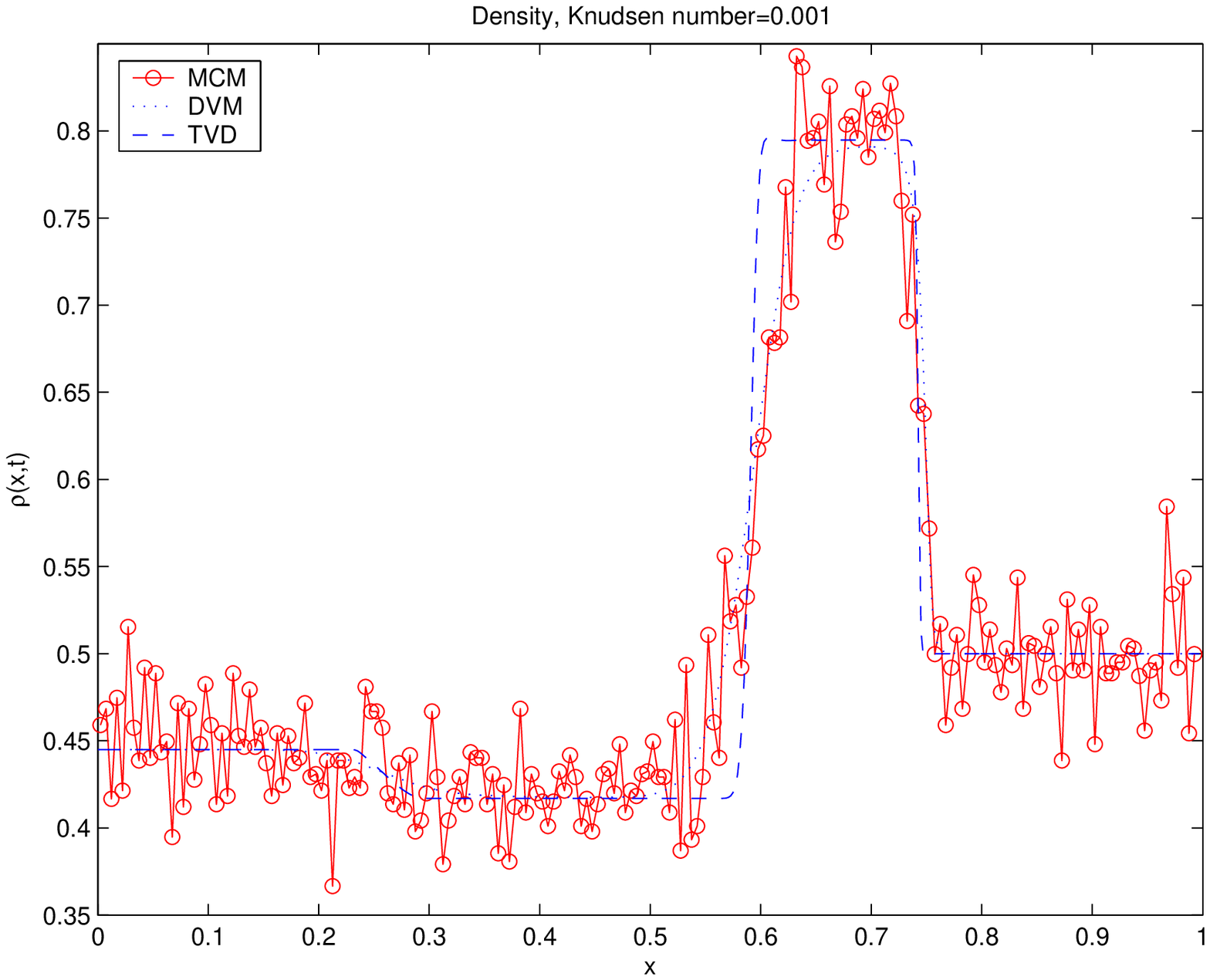}
\includegraphics[scale=0.40]{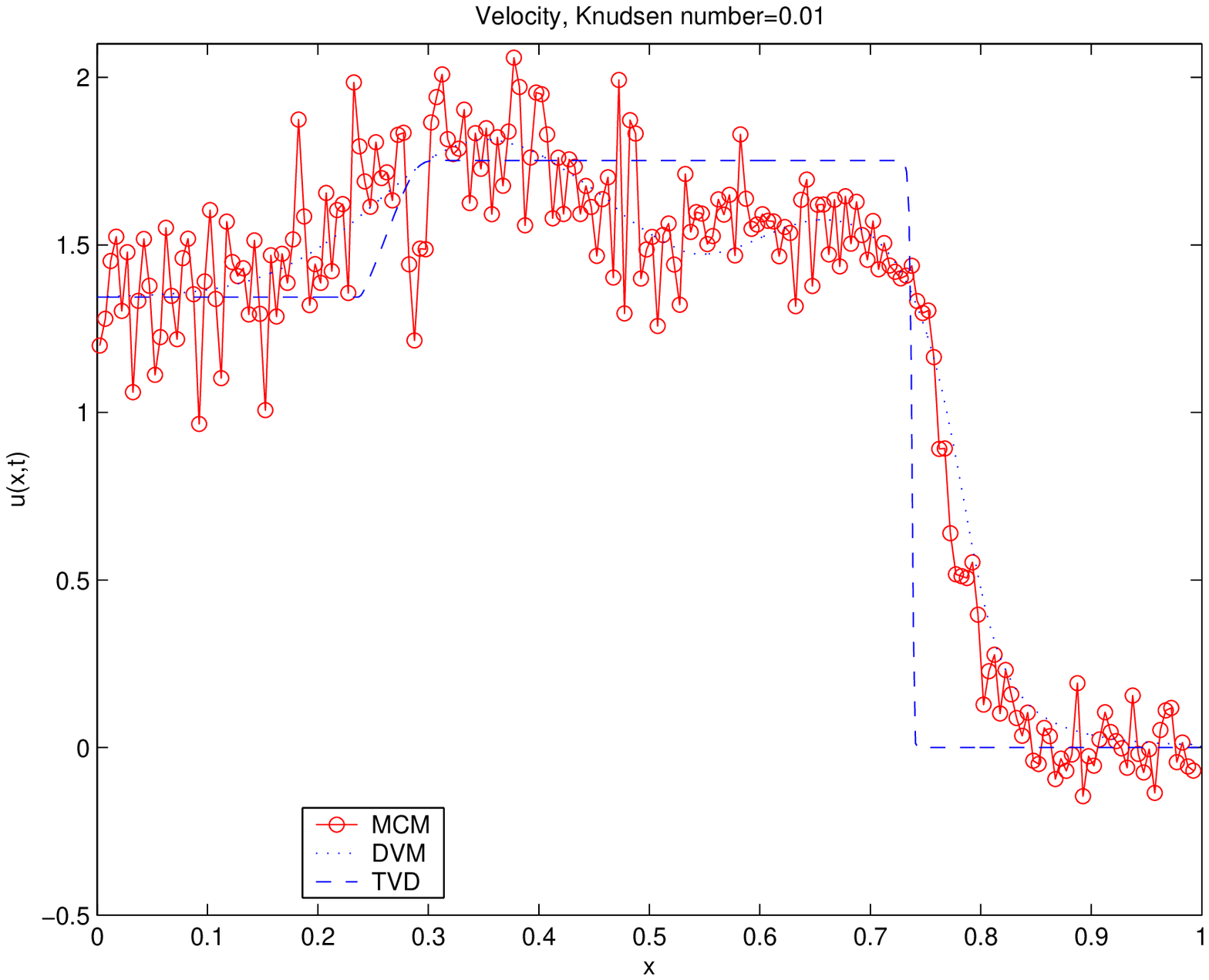}
\includegraphics[scale=0.40]{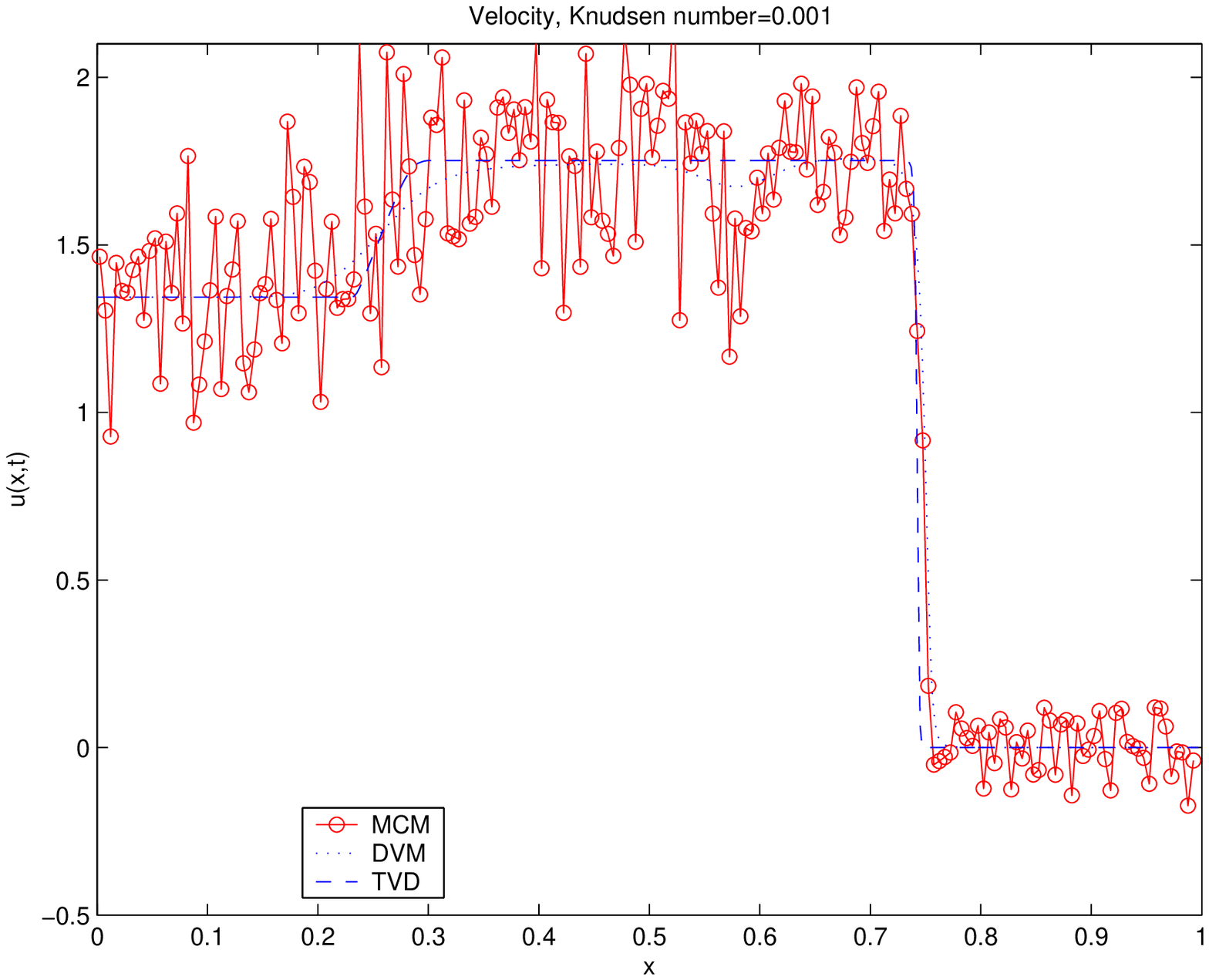}
\includegraphics[scale=0.40]{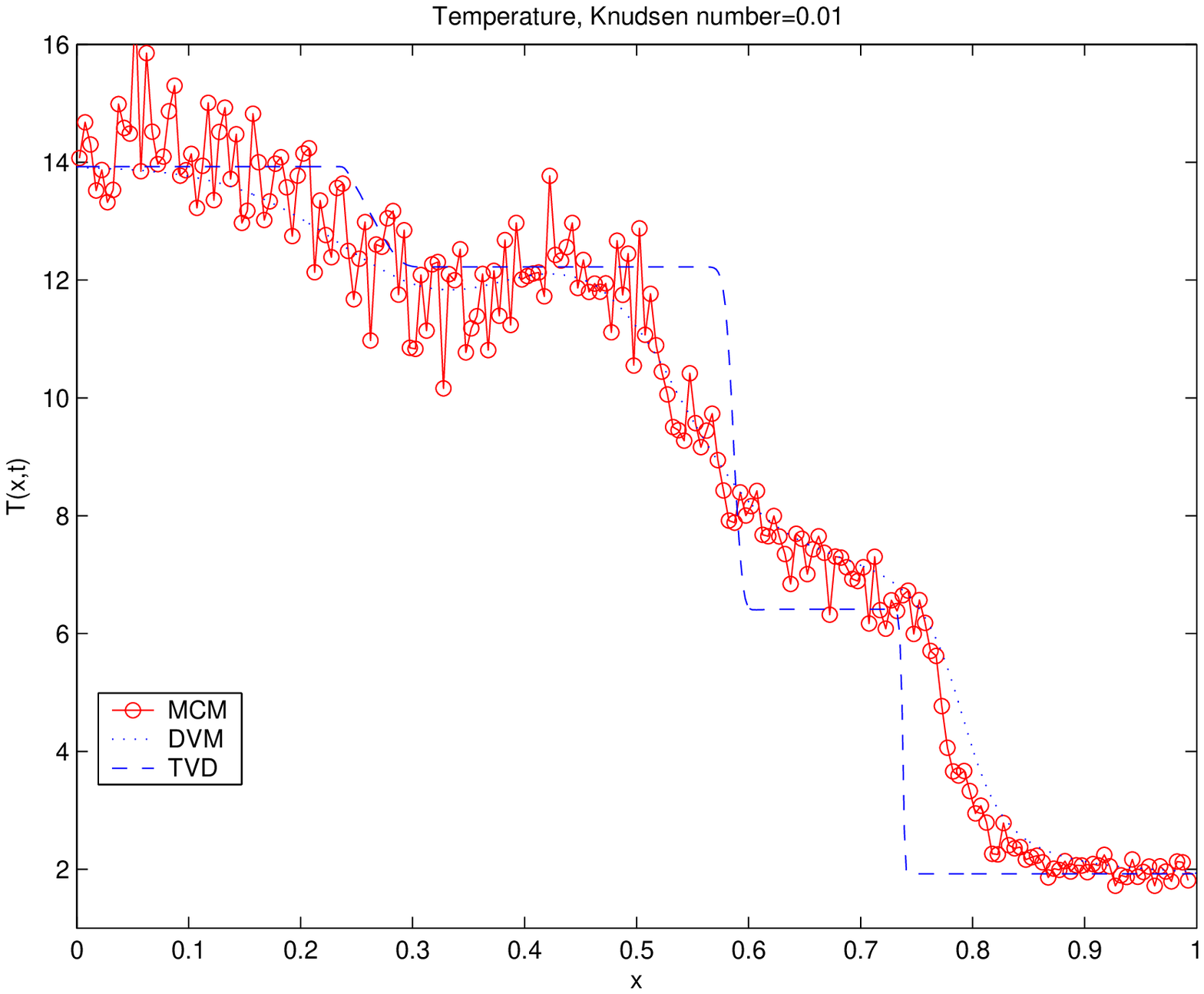}
\includegraphics[scale=0.40]{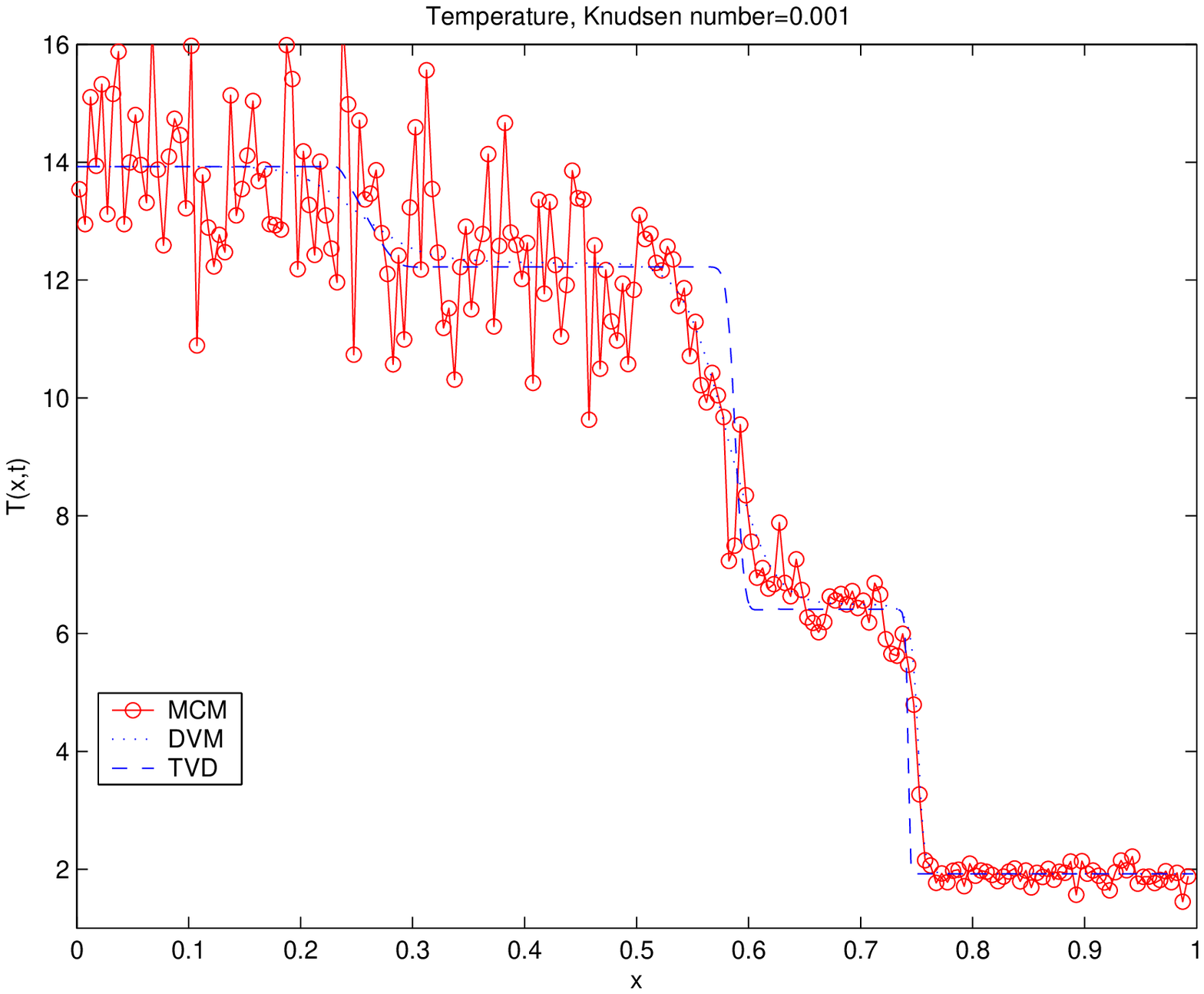}
\caption{Lax Test. Solution at $t=0.05$ for MCM with Knudsen
numbers $\varepsilon=5\times 10^{-2}$ (left) and
$\varepsilon=10^{-3}$ (right). From top to bottom density, mean
velocity and temperature.} \label{L6}
\end{center}
\end{figure}

Finally we consider a Lax shock tube test with initial values
\begin{equation}
\nonumber
\begin{array}{l}
\textbf{u}_{L}=\left(
\begin{array}{ll}
0.445 \\0.598\\3.5
\end{array}
\right), \ \hbox{if} \ 0  \leq x <0.5 \ \ \
\textbf{u}_{R}=\left(\begin{array}{l} 0.5 \\ 0 \\ 0.48
\end{array}\right), \ \hbox{if} \ 0.5\leq x \leq 1.
\end{array}
\label{eq:Sod}
\end{equation}
The solution is computed with $200$ grid points in space, the
final time is $t=0.05$. The initial number of particle is $500$
for each space cell.
Each Figure contains the DVM solution and the Euler solution as
reference. Similar considerations to those of the previous section
hold for this test case. Thus for large Knudsen numbers the
solutions computed with the hybrid methods (Figure
\ref{L2}-\ref{L3}) show small improvements compared to the Monte
Carlo scheme (Figure \ref{L6}). On the other hand when the Knudsen
number becomes smaller FSI and FSI1 schemes (Figure
\ref{L4}-\ref{L5}) give a considerable reduction of fluctuations.
This is especially true for FSI1 method which demonstrates the
importance of a good estimate of the equilibrium fraction
$\beta^c$ after the transport.


\begin{figure}
\begin{center}
\includegraphics[scale=0.40]{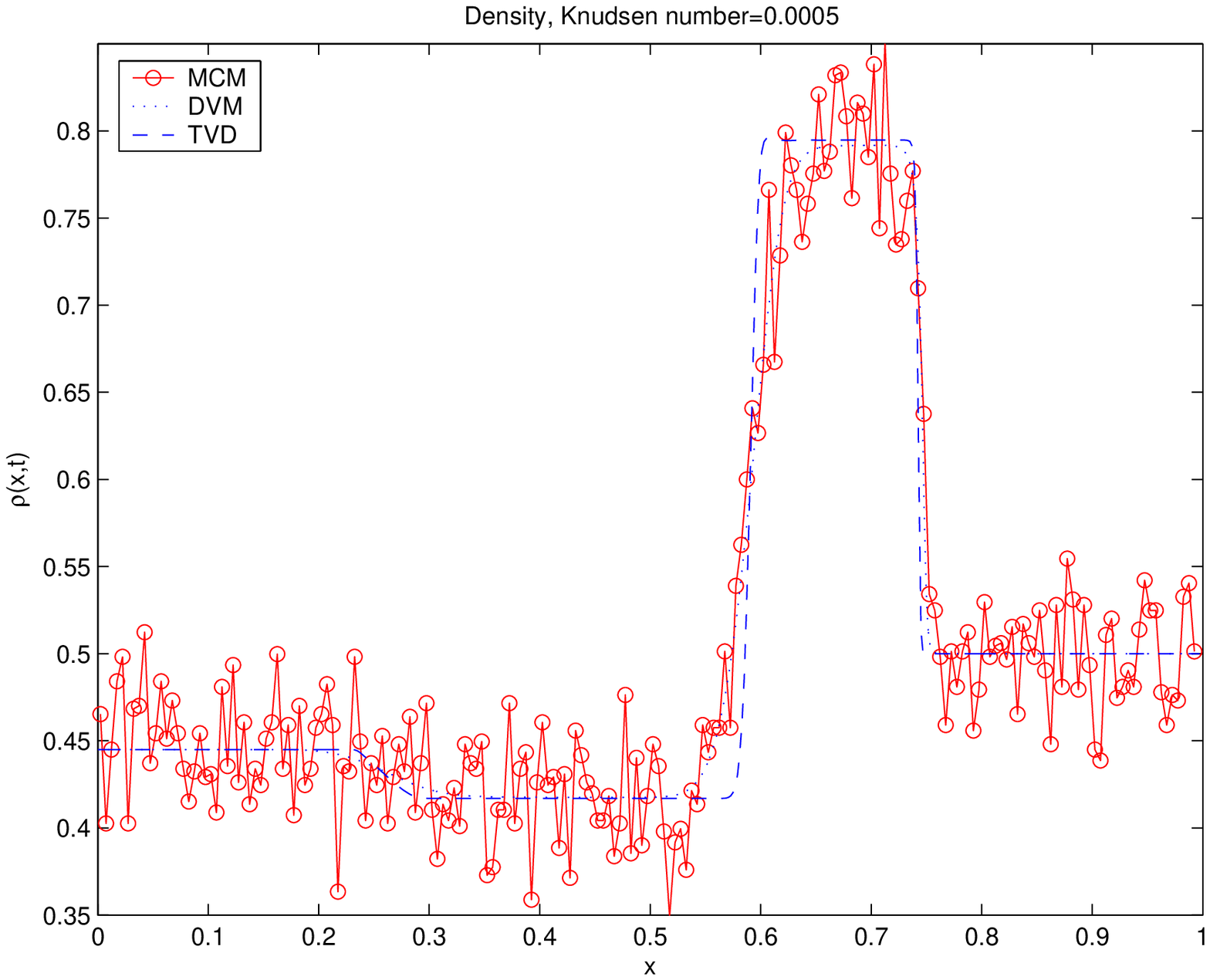}
\includegraphics[scale=0.40]{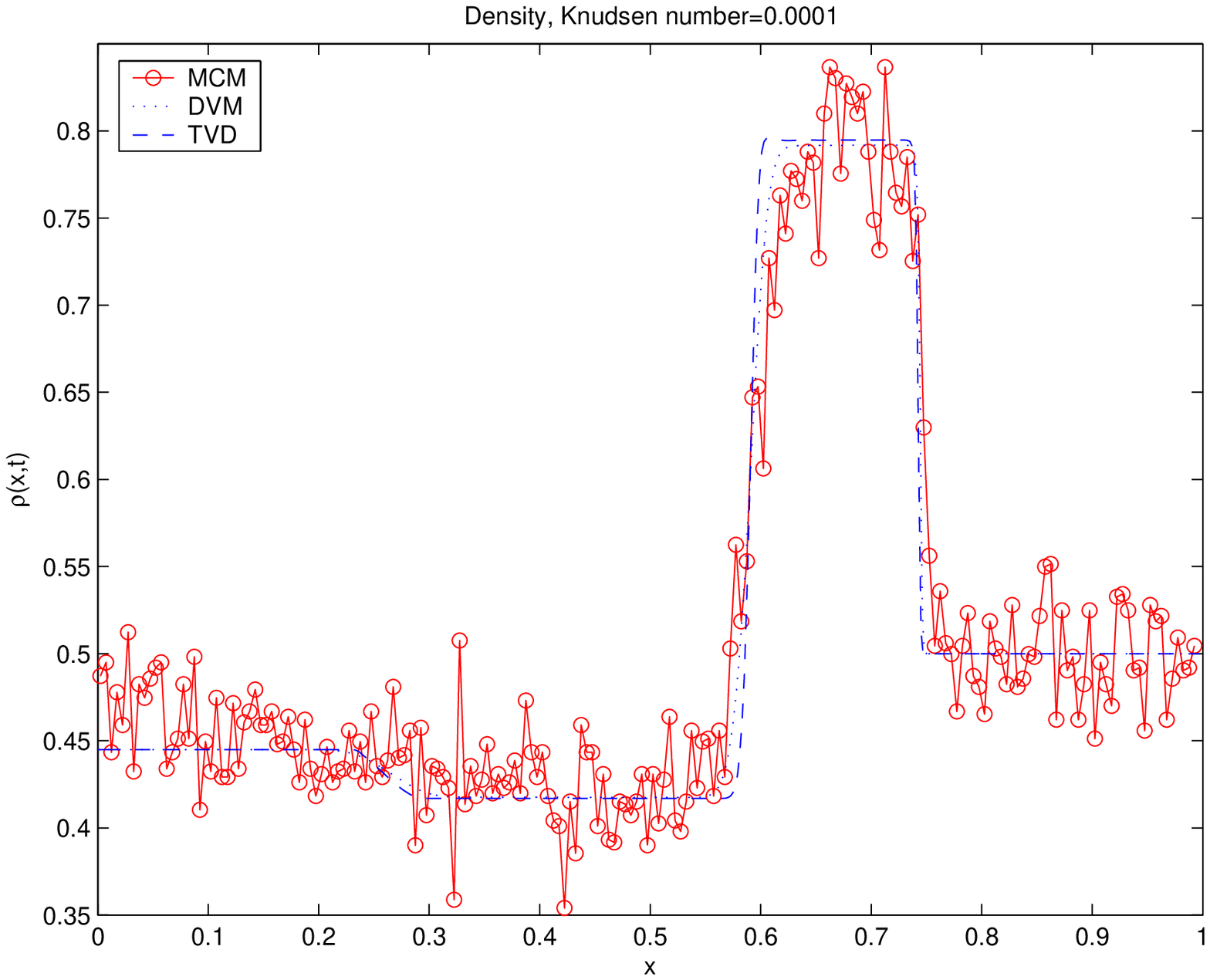}
\includegraphics[scale=0.40]{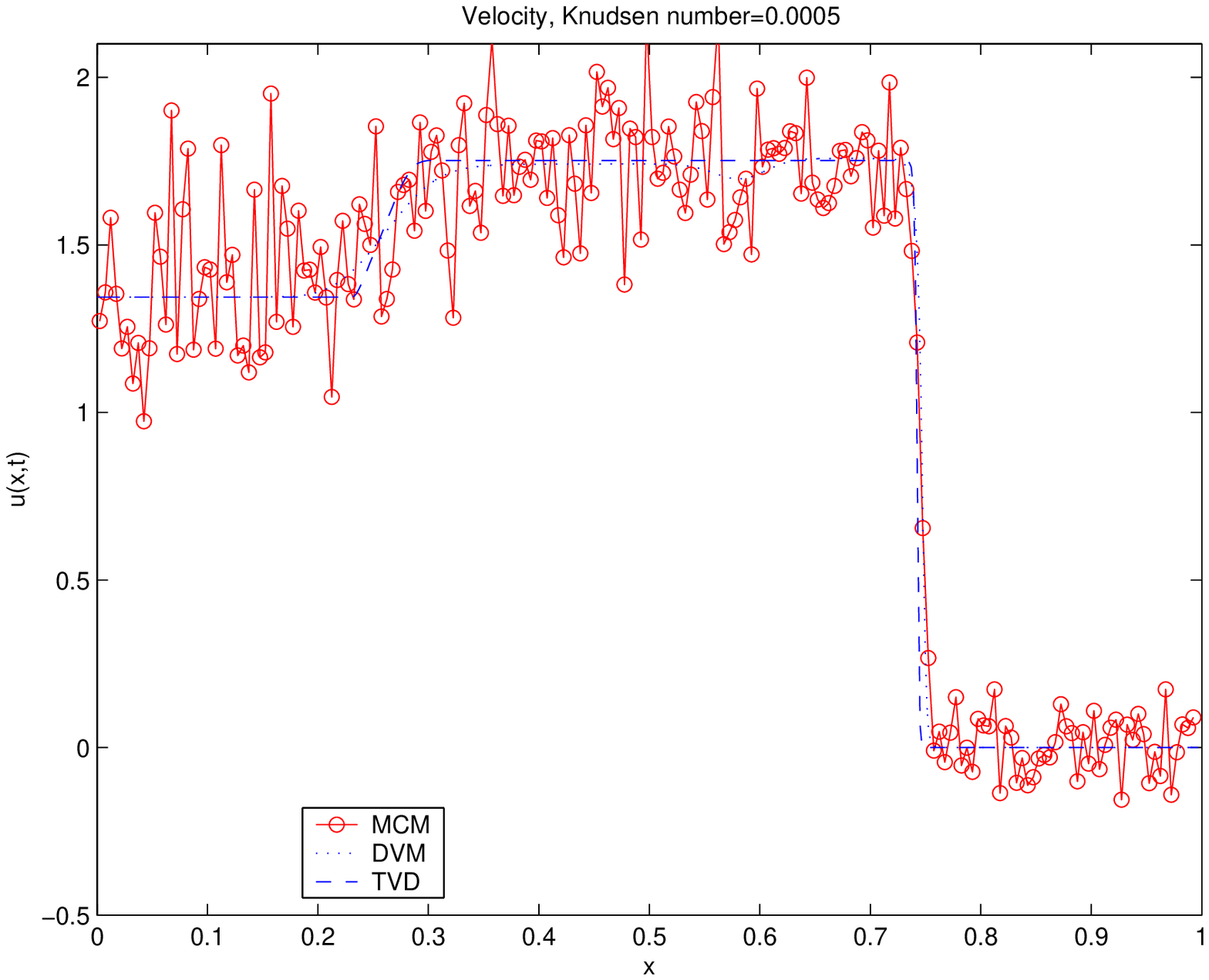}
\includegraphics[scale=0.40]{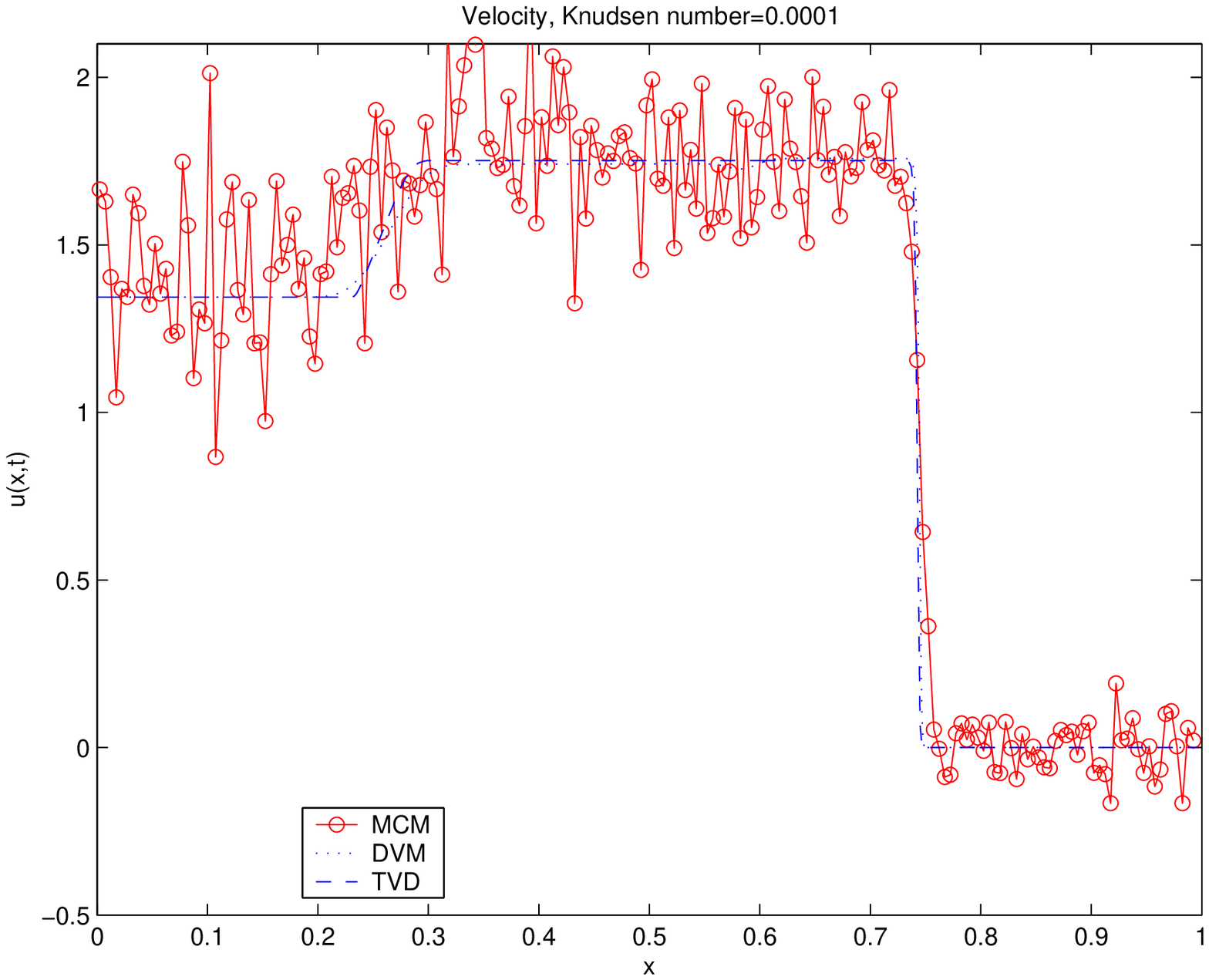}
\includegraphics[scale=0.40]{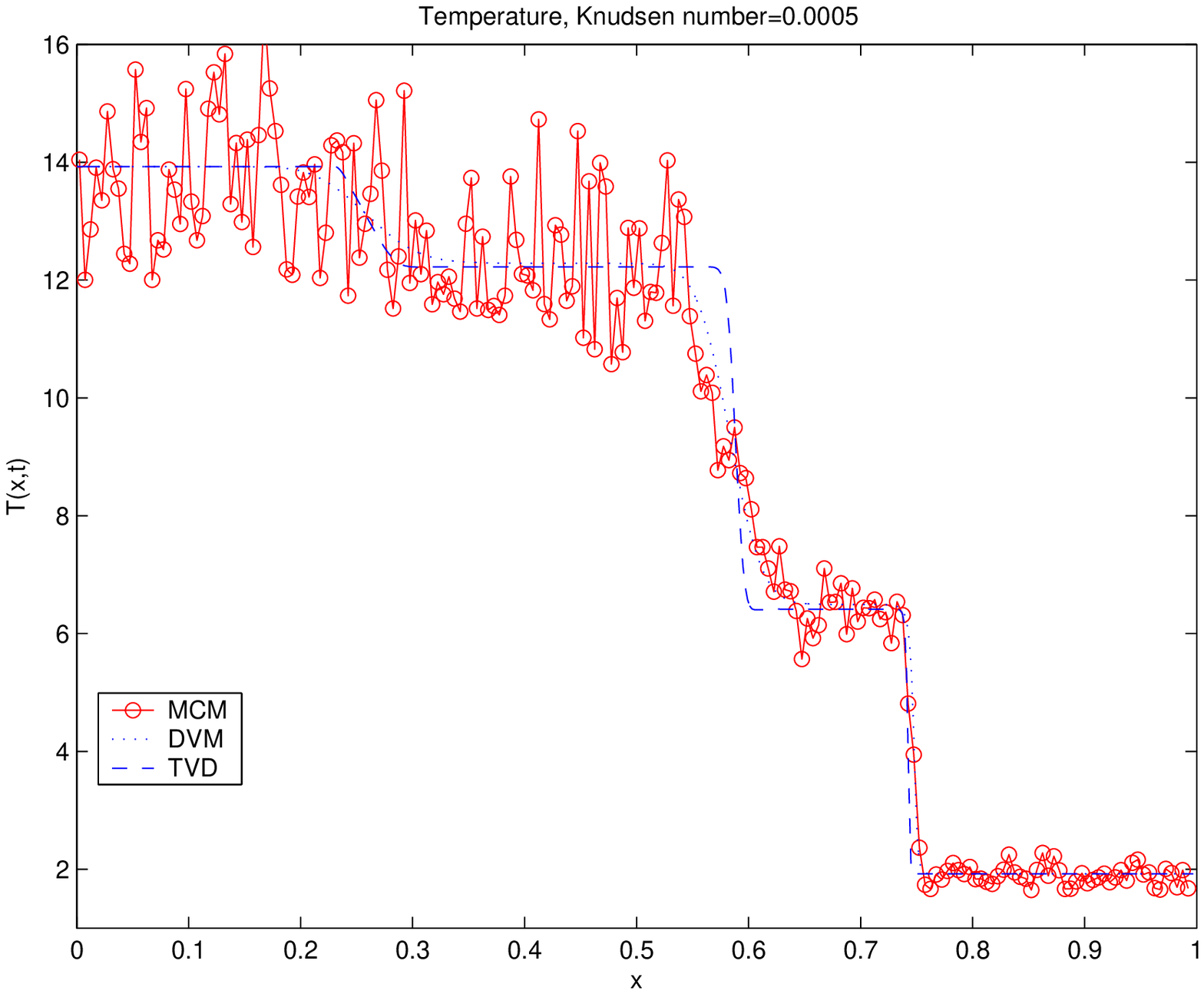}
\includegraphics[scale=0.40]{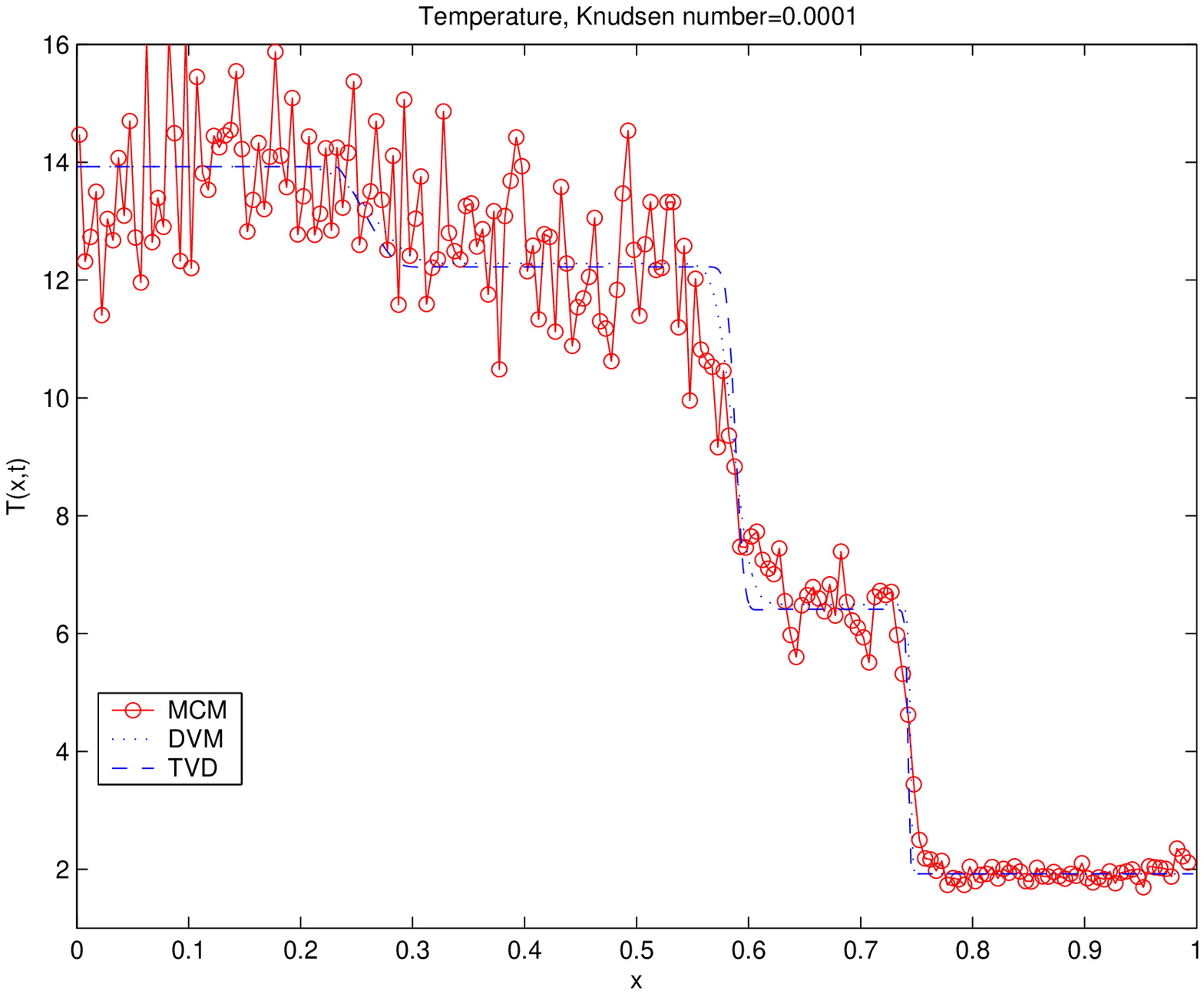}
\caption{Lax Test. Solution at $t=0.05$ for MCM with Knudsen
numbers $\varepsilon=5\times 10^{-4}$ (left) and
$\varepsilon=10^{-4}$ (right). From top to bottom density, mean
velocity and temperature.} \label{L7}
\end{center}
\end{figure}

\section{Conclusion}

In this paper we have extended the hybrid kinetic methods
developed in \cite{dimarco3} to the case of an arbitrary fluid
solver for the equilibrium component. Although, the simplified BGK
collision operator has been used to develop the schemes,
extensions to the full Boltzmann operator of rarefied gas dynamics
in principle are possible through the use of time relaxed methods
\cite{CPmc, PR}. We plan to extend the schemes to the full
Boltzmann equation in the nearby future.

The results obtained are very promising in terms of computational
cost with respect to traditional deterministic methods for kinetic
equations like discrete velocity model or spectral schemes
\cite{Mieussens, Russo1}. In addition the FSI hybrid algorithms
yield less fluctuations with respect to direct simulation Monte
Carlo methods and, close to the fluid regime, they permit to
compute results faster. A remarkable feature of the FSI1 scheme is
that the equilibrium fraction is essentially independent of the
choice of the time step and thus provides more accurate results
then Monte Carlo methods even in resolved regimes.

Some open questions remain on alternative ways to estimate and
increase the fraction of equilibrium in each space cell without
increasing the computational cost. It is also interesting, and
will be the subject of future works, to measure the response of
the FSI hybrid methods in others situations such as simulations of
nanosystem devices, plasma physics problems or turbulence.


\vskip 2ex {\bf Acknowledgements.} The authors would like to thank
Russ Caflisch and Pierre Degond for several stimulating
discussions.

\end{document}